\documentclass[review,onefignum,onetabnum]{siamart171218}

\usepackage{lipsum}
\usepackage{amsfonts}
\usepackage{graphicx}
\usepackage{epstopdf}
\usepackage{algorithmic}
\usepackage{array}

\ifpdf
  \DeclareGraphicsExtensions{.eps,.pdf,.png,.jpg}
\else
  \DeclareGraphicsExtensions{.eps}
\fi

\newcommand{\ie}{{\it i.e.}}
\newcommand{\eg}{{\it e.g.}}

\newcommand{\st}{{\it s.t.}}

\newsiamremark{remark}{Remark}
\newsiamremark{hypothesis}{Hypothesis}
\crefname{hypothesis}{Hypothesis}{Hypotheses}
\newsiamthm{claim}{Claim}

\headers{Weak Adversarial Networks for Constrained Optimization}{G. Bao, D. Wang, and B. Zou}

\title{WANCO: Weak Adversarial Networks for Constrained Optimization problems\thanks{Submitted to the editors DATE.
\funding{G. Bao was partially supported by the National Key Research and Development Program of China (2023YFA1009100) and the Key Project of Joint Funds for Regional Innovation and Development (U21A20425). D. Wang was partially supported by National Natural Science Foundation of China (Grant No. 12101524), Guangdong Basic and Applied Basic Research Foundation (Grant No. 2023A1515012199), Shenzhen Science and Technology Innovation Program (Grant No. JCYJ20220530143803007, RCYX20221008092843046), Guangdong Provincial Key Laboratory of Mathematical Foundations for Artificial Intelligence (2023B1212010001), and Hetao Shenzhen-Hong Kong Science and Technology  Innovation Cooperation Zone Project (No.HZQSWS-KCCYB-2024016).
}}}

\author{Gang Bao\thanks{School of Mathematical Sciences, Zhejiang University, Zhejiang, China 
  (\email{baog@zju.edu.cn}).}
\and Dong Wang\thanks{School of Science and Engineering, The Chinese University of Hong Kong, Shenzhen, Guangdong 518172, China \& Shenzhen International Center for Industrial and Applied Mathematics, Shenzhen Research Institute of Big Data, Guangdong 518172, China (\email{wangdong@cuhk.edu.cn}).} 
\and Boyi Zou\thanks{School of Science and Engineering, The Chinese University of Hong Kong, Shenzhen, Guangdong 518172, China (\email{boyizou@link.cuhk.edu.cn}).}}

\usepackage{amsopn}

\makeatletter
\newcommand*{\addFileDependency}[1]{
  \typeout{(#1)}
  \@addtofilelist{#1}
  \IfFileExists{#1}{}{\typeout{No file #1.}}
}
\makeatother

\ifpdf
\hypersetup{
  pdftitle={WANCO: Weak Adversarial Networks for Constrained Optimization problems},
  pdfauthor={G. Bao, D. Wang, and B. Zou}
}
\fi

\begin{document}

\nolinenumbers
\maketitle

\begin{abstract}
This paper focuses on integrating the networks and adversarial training into constrained optimization problems to develop a framework algorithm for constrained optimization problems. For such problems, we first transform them into minimax problems using the augmented Lagrangian method and then use two (or several) deep neural networks(DNNs) to represent the primal and dual variables respectively. The parameters in the neural networks are then trained by an adversarial process. The proposed architecture is relatively insensitive to the scale of values of different constraints when compared to penalty based deep learning methods. Through this type of training, the constraints are imposed better based on the augmented Lagrangian multipliers. Extensive examples for optimization problems with scalar constraints, nonlinear constraints, partial differential equation constraints, and inequality constraints are considered to show the capability and robustness of the proposed method, with applications ranging from Ginzburg--Landau energy minimization problems, partition problems, fluid-solid topology optimization, to obstacle problems.

\end{abstract}



\section{Introduction}

Constrained optimization problems arise from an intricate and practical set of issues, with the objective of identifying solutions to an optimization problem while adhering to specific constraints. These problems play a crucial role in various domains such as finance, economics, engineering, operation research, machine learning, and others. For a comprehensive exploration of constrained optimization problems, we refer to~\cite{bertsekas2014constrained} and the references therein. A wide range of problems in scientific computing can be classified as constrained optimization problems, for example, inverse scattering problems~\cite{Bao_2015,bao2022mathematical}, optical sciences \cite{bao2001mathematical}, topology optimization problems~\cite{bendsoe2013topology}, operation research, optimal control, and obstacle problems.

In this paper, we focus on constrained optimization problems related to differential operators that can appear in either the objective functional or the constraints. There are many efficient classical methods for solving constrained optimization problems~\cite{bertsekas2014constrained, glowinski1989augmented, wright1997primal}, such as the penalty method, the Lagrange multiplier method, the augmented Lagrangian method, and the interior point method. When solving differential operator related problems with constraints, traditional numerical methods such as the finite difference method, the finite element method, or the spectral method can be applied with careful construction, discretization, and post-processing for satisfying the constraints.

In recent years, there has been rapid development in artificial intelligence (AI) for science, with particular attention being drawn to physics-informed machine learning in scientific computing. The universal approximation of continuous functions using neural networks, as proven in~\cite{Cybenko_1989}, has gained significant attention in the past decade, leading to an increasing trend of utilizing deep neural networks to solve differential operator related problems. A lot of studies have demonstrated the effectiveness of deep learning approaches in solving partial differential equations, optimization problems with differential equation constraints, or inverse problems, including physics-informed neural network (PINN)~\cite{Raissi_2019}, deep Galerkin method (DGM)~\cite{Sirignano_2018}, deep Ritz method (DRM)~\cite{E_2018}, weak adversarial network (WAN)~\cite{Zang_2020, Bao_2020}, random feature method (RFM)~\cite{Chen_2022_RF}, local extreme learning machines (locELM)~\cite{Dong_2021}, particle weak-form based neural networks (ParticleWNN) \cite{zang2023particlewnn} and many other mechanism driven methods~\cite{kharazmi2019variational,Yu_2021,Linka_2022,McClenny_2022,Yu_2022,Wu_2023,Wang_2024}. In contrast with traditional methods, deep learning methods are usually mesh-free. Therefore, they are less affected by the curse of dimensionality, making them efficient for computing high-dimensional problems. For a comprehensive overview of using DNNs to solve differential operator related problems, we refer  to~\cite{Lu_2021, Karniadakis_2021,Cuomo_2022} and references therein.

The framework for using DNNs to solve differential operator constrained optimization problems involves constructing a problem-specific loss function, followed by solving an optimization problem with boundary conditions and other constraints. One common approach is to utilize penalty methods to transform the constrained optimization problem into an unconstrained optimization problem, as used in PINN, DGM, and DRM. The results are highly dependent on the selection of the weights of different terms in the loss function. In particular, one needs to carefully consider the scale of the values of constraints to guarantee that the training process can minimize each term in the loss function equally. The construction method (\eg, \cite{Lyu_2021, Sukumar_2022}) is another way to deal with specific boundary constraints, while it is only efficient to particular constraints and computational domains.

Except for the penalty or construction based methods, one may consider the Lagrange method or augmented Lagrangian method(ALM) to transform the constrained optimization problem into a minimax problem.  It has been widely used for solving constrained optimization problems involving neural networks and been successfully applied in various fields, such as natural language processing (NLP)~\cite{nandwani2019primal}, classification~\cite{sangalli2021constrained} and some others~\cite{Li_2021, Qi_2022, Burman_2023}. As for problems related to differential operators, neural networks are used to represent the decision variables, and then the problem is discretized, with the introduction of Lagrange multipliers on some fixed grid points. For example, in~\cite{Hwang_2022}, ALM was used to impose equality constraints in addressing the physical conservation laws of kinetic PDE problems, which resulted in more accurate approximations of the solutions in terms of errors and the conservation laws. In~\cite{Lu_2021hard}, the applications of ALM, penalty method, and adaptive penalty method for inverse design problems were compared under the framework of physics-informed neural networks (PINNs). The results demonstrated that ALM is relatively more robust and less sensitive to hyperparameters. Enhanced physics-informed neural networks with augmented Lagrangian relaxation method (AL-PINNs)~\cite{Son_2023} treated the initial and boundary conditions as constraints for the optimization problem of the PDE residual. Various numerical examples showed that AL-PINNs yield a much smaller relative error compared to state-of-the-art adaptive loss-balancing algorithms. The augmented Lagrangian deep learning (ALDL)~\cite{Huang_2022} constructed the Lagrange multipliers in the sense of adversarial network to impose the essential boundary conditions of variational problems, avoiding the discretization of the Lagrange multipliers. However, the papers mentioned above either rely on the ALM form, with the Lagrange multipliers defined and updated on a batch of fixed grid points, or only handle boundary conditions as constraints in variational problems. These methods are restricted  by the selection of grid points and may be inefficient for high-dimensional problems due to the curse of dimensionality. Additionally, they are difficult to integrate with other improved frameworks, such as adaptive methods that need to adaptively sample the collocation points.

Inspired by the framework of WAN~\cite{Zang_2020, Bao_2020}, we propose the weak adversarial networks for constrained optimization (WANCO) to solve constrained optimization problems.  WANCO has a wide range of applications and can handle various constraints, including for example, scalar constraints, linear and nonlinear constraints, PDE constraints, and inequality constraints. It is based on the augmented Lagrangian form of the constrained optimization problems by representing decision variables and Lagrange multipliers using individual networks. The problems are then trained through an adversarial process. The proposed WANCO is less sensitive to the weights of different terms in the loss function and can be easily integrated with other techniques to improve the performance of the training result. In the numerical examples, we demonstrate the capability and robustness of the parameter setting of WANCO on various constrained optimization problems, including the Ginzburg--Landau energy minimization problems, Dirichlet partition problems in different dimensions, fluid-solid optimization problems, and obstacle problems. 

The rest of this paper is organized as follows. In Section~\ref{Sec: Framework of WANCO}, we introduce the proposed WANCO with its network architectures. In Section~\ref{Subsec: Ginzburg--Landau Energy}, we use the mass preserving Ginzburg--Landau energy minimization problem as an example to demonstrate the performance of WANCO. In particular, Section~\ref{subsection:GLcomparison} presents the comparison results to other deep learning based approaches, demonstrating the insensitivity of WANCO to parameters, and Section~\ref{subsection:activation} shows comparison results to shed light on the choice of activation functions. In Section~\ref{Subsec: Dirichlet partition}, we investigate the high dimensional performance of WANCO by considering Dirichlet partition problems up to four dimensions. The PDE constrained optimization problems are considered in Section~\ref{Subsec: Fluid-solid optimization} for fluid-solid optimization problems as the example. In Section~\ref{Subsec: Obstacle problem}, WANCO solves the obstacle problems involving inequality constraints. We draw a conclusion and discussion in Section~\ref{Sec: Conclusions}.

\section{Framework of WANCO}\label{Sec: Framework of WANCO}

To introduce the weak adversarial networks for constrained optimization problems (WANCO), we consider a general constrained optimization problem in the following form:
\begin{equation}\label{Eq: Constrained Optimization Original Form}
\begin{aligned}
     &\min \limits_{u}\: L(u)\\
    &\st \:    C(u) = 0,
\end{aligned}
\end{equation}
where $u = u(x)$ represents the decision variable that needs to be optimized, $L$ denotes the objective functional of $u$, and $C(u)=0$ represents the constraint that can describe various types of constraints, such as scalar constraints, boundary constraints, PDE constraints, and even inequality constraints by introducing slack variables. This problem finds applications in diverse fields, including fluid mechanics, material science, quantum mechanics, industrial design, and topology optimization. 

For constrained optimization problems~\eqref{Eq: Constrained Optimization Original Form}, we first transform it into a minimax problem with the Lagrange multiplier $\lambda$ and a penalty term as follows,
\begin{equation}\label{Eq: Augmented Lagrangian}
    \min\limits_{u}\max\limits_{\lambda}\:L_{\beta}(u,\lambda)=  C_0 L(u)-\left\langle\lambda,C(u)\right\rangle + \frac{\beta}{2}\|C(u)\|^2, 
\end{equation}
in which $C_0$ and $\beta$ are scalar coefficients. $C_0$ is a fixed parameter used to avoid the objective functional value being too small in specific problems. $\beta$ gradually increases during the optimization process, similar to that in the augmented Lagrangian method.

\begin{center}
\begin{figure}[htb!]
\centering
\includegraphics[width=0.85\textwidth]{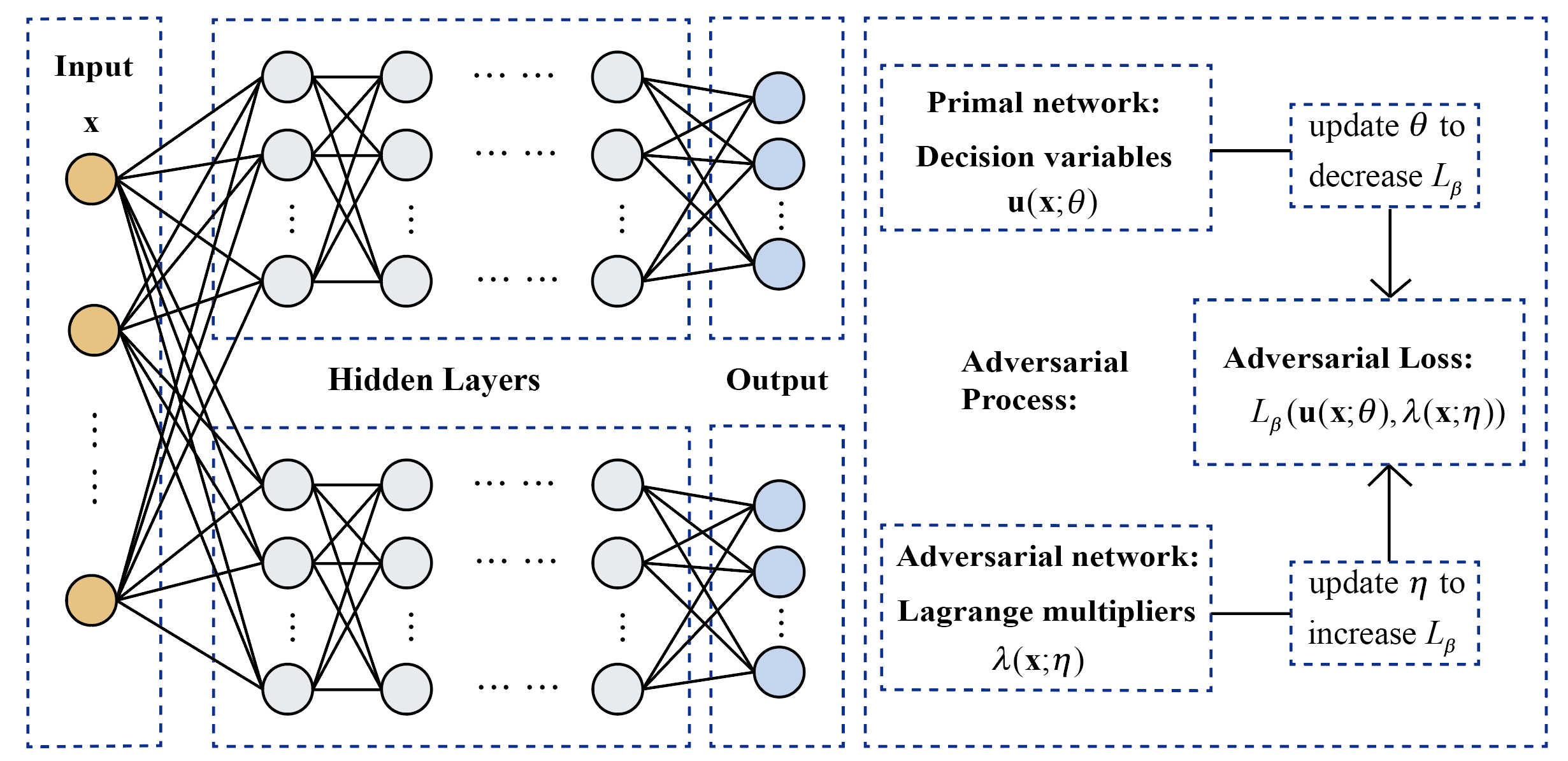} 
\caption{ Network structure of WANCO, see Section~\ref{Sec: Framework of WANCO}.}
\label{Fig: WANCO network structure}
\end{figure} 
\end{center}

In WANCO, we propose to use two(or several) individual deep neural networks (similar to WAN and IWAN) to represent $u = u(x;\theta)$, and $\lambda = \lambda(x;\eta)$ with parameters $\theta$ and $\eta$, respectively. Problem~\eqref{Eq: Augmented Lagrangian} is then transformed into a minimax problem with respect to $\theta$ and $\eta$,
\begin{equation}\label{Eq: WANCO method}
    \min\limits_{\theta}\max\limits_{\eta}\:L_{\beta}(u(x;\theta),\lambda(x;\eta))= C_0 L(u(x;\theta))-
    \langle\lambda(x;\eta),C(u(x;\theta))\rangle + \frac{\beta}{2}\|C(u(x;\theta))\|^2. 
\end{equation}
The parameters $\theta$ and $\eta$ are then alternatively updated by stochastic gradient descent and stochastic gradient ascent of $L_{\beta}(u(x;\theta),\lambda(x;\eta))$, with a gradually increasing value of $\beta$. The algorithmic framework of WANCO is summarized in Algorithm~\ref{alg: WANCO framework}, and the neural network diagram of WANCO is illustrated in Figure~\ref{Fig: WANCO network structure}.
\begin{algorithm}[htb]
  \caption{WANCO: weak adversarial networks for constrained optimization problems.}  
  \label{alg: WANCO framework}
  \begin{algorithmic}
  \STATE{\bf Inputs:}\\
  \STATE{Two (or several) individual neural networks to represent the decision variable $u = u(x;\theta)$ and Lagrange multiplier $\lambda = \lambda(x;\eta)$.\\
  $C_0$, $\beta$: weight parameter of the objective functional and constraint.\\
$\alpha$: the amplifier of $\beta$ in the augmented Lagrangian method.\\ 
$N$: number of iterations.\\ 
$N_u$, $N_{\lambda}$: number of primal and adversarial network parameter updates in each iteration.\\
$\tau_\theta$, $\tau_\eta$: learning rate of primal and adversarial network.}\\
    \FOR{i = $1:N$}
    \STATE{Update the Loss function $L_{\beta}(u(x;\theta),\lambda(x;\eta))$ in Equation~\eqref{Eq: WANCO method}.} 
        \FOR{j = $1:N_u$}
            \STATE{{Update $\theta \leftarrow \theta-\tau_\theta \nabla_\theta L_\beta$.}}
            \ENDFOR
        \FOR{k = $1:N_\lambda$} 
            \STATE{{Update $\eta \leftarrow \eta+\tau_\eta \nabla_\eta L_\beta$.}} 
            \ENDFOR
        \STATE{{Update $\beta\leftarrow \alpha*\beta$.}}
    \ENDFOR
    \STATE{{\bf Outputs:} Neural networks $u(x;\theta)$ and $\lambda(x;\eta)$ with arbitrary input $x$ and trained parameters $\theta$ and $\eta$.}
\end{algorithmic}  
\end{algorithm}

In terms of the network architecture, we represent the decision variables by the residual neural networks (ResNet). The Lagrange multiplier is represented by different networks based on their properties. For example,\\
1. if the constraint is scalar, $\lambda(x;\eta)$ can be simply set as $\lambda(0;\eta)$ with constant input $0$ and a shallow neural network,\\
2. if the constraint is defined over the entire computational domain ({\it e.g.}, PDE constraints), $\lambda(x;\eta)$ is set to be a function in the domain represented by a ResNet,\\ 
3. if the constraint is an inequality, structured neural networks can be used to enforce the non-positivity of $\lambda(x;\eta)$.

To be specific, we use the ResNet~\cite{He_2016} with $N_d$ hidden layers and width $N_w$, expressed as follows,
\begin{equation}\label{Eq: u-net RNN}
    \begin{aligned}
        &\mathbf{u}(\mathbf{x}_0;\theta) = \mathbf{W}_{out}\mathbf{x}_{N_d +1} + \mathbf{b}_{out},  \\
        &\mathbf{x}_{k+1} = f_k(\mathbf{x}_k) = \Phi_k(\mathbf{W}_k\mathbf{x}_k+\mathbf{b}_k)+\mathbf{x}_k ,\quad k = N_d,...,1,\quad \\
        &\mathbf{x}_1 =  \mathbf{W}_{in}\mathbf{x}_0+\mathbf{b}_{in}.\quad 
    \end{aligned}
\end{equation}
Here, $\mathbf{x}_0\in \Omega$ is the input of dimension $d$, $\{\mathbf{W}_{in}\in\mathbb{R}^{N_w\times d}, \mathbf{W}_{k}\in\mathbb{R}^{N_w\times N_w}, \mathbf{W}_{out}\in\mathbb{R}^{n\times N_w}\}$ are the weight matrices of the input layer, $k$-th hidden layer, and output layer, $\{\mathbf{b}_{in}\in\mathbb{R}^{N_w}, \mathbf{b}_{k}\in\mathbb{R}^{N_w}, \mathbf{b}_{out}\in\mathbb{R}^{n}\}$ are the corresponding bias vectors, and $\Phi_k$ represents the activation function. The parameters $N_d$ and $N_w$ are used to define the size of the weight matrices and bias vectors in the neural network, which ultimately determine the network's performance and capabilities from many empirical observations and the universal approximation theorem. ResNet effectively mitigates the problems of gradient vanishing and gradient explosion by incorporating a residual connection or skip connection, represented by the $+\mathbf{x}_k$ term in the second equation of~\eqref{Eq: u-net RNN}. As for the activation function $\Phi_k$, we will use $\tanh^3$ instead of common activation functions~\cite{Sharma_2020} (\eg, $\mathrm{sigmoid}$, $\tanh$, $\mathrm{softmax}$, $\mathrm{ReLU}$, and $\mathrm{ReLU}^3$). A detailed numerical comparison will be provided in Section~\ref{Subsec: Ginzburg--Landau Energy}.

After constructing the saddle point problem~\eqref{Eq: WANCO method}, the loss function is optimized iteratively by the standard optimizer {\tt Adam}~\cite{kingma2014adam} with the stochastic gradient descent(SGD) and ascent method. When computing different terms of the loss function, we employ several standard DNN techniques, including the Monte Carlo (MC) method for integrals and auto-differentiation for derivatives. To avoid getting stuck in undesirable local minima, we incorporate shuffling of training points during the training process, as shown in \cite{E_2018} for example. To simplify the optimization process, we design specific networks that automatically satisfy boundary or positivity constraints. In cases where these constraints are not straightforward, they can be imposed as additional adversarial networks in the WANCO framework. The numerical results presented in Sections~\ref{Subsec: Ginzburg--Landau Energy}-\ref{Subsec: Obstacle problem} demonstrate that WANCO is a versatile algorithm that can be combined with various DNN techniques to enhance accuracy and simplify the optimization process. The specific techniques will be discussed when they are used in different examples.

The proposed WANCO can be applied to general constrained optimization problems. In this paper, we specifically focus on problems that involve differential operators either in the objective functional or the constraints, with applications in free interface related optimization problems. In the follows, we consider several typical examples to demonstrate the performance of the proposed methods:

1. linear scalar constraints: Ginzburg--Landau energy minimization with mass preserving (see Section~\ref{Subsec: Ginzburg--Landau Energy});

2. nonlinear scalar constraints: Dirichlet partition problems (see Section~\ref{Subsec: Dirichlet partition});

3. partial differential equation constraints: fluid-solid optimization (see Section~\ref{Subsec: Fluid-solid optimization}); and

4. inequality constraints: obstacle problems (see Section~\ref{Subsec: Obstacle problem}).

In the follows, we first study the parameter effect and compare the performance of different activation functions using the mass preserving Ginzburg--Landau energy minimization problem. We then consider the Dirichlet partition problems to show the high-dimensional performance of WANCO. Then, we extend to consider fluid-solid optimization with PDE constraints. Furthermore, we explore the performance of WANCO on problems with inequality constraints through obstacle problems.

\section{Minimization of the Ginzburg--Landau energy with mass conservation}\label{Subsec: Ginzburg--Landau Energy}

In this section, we focus on the empirical study of the sensitivity of the parameters, the choice of activation functions, and the comparison to the penalty based methods and adversarial network with only Lagrange multipliers.  We illustrate them through the minimization of the Ginzburg--Landau energy with mass conservation as an example.

Consider the minimization of the Ginzburg--Landau energy with mass conservation and homogeneous Dirichlet boundary condition in $\Omega=[0,1]\times[0,1]$ as follows,
\begin{equation}
    \begin{aligned}
     &\min\limits_{u}\: \int_\Omega\frac{\epsilon}{2}|\nabla u|^2+\frac{1}{\epsilon}(u^2-1)^2 \ d\mathbf{x} \\
    \st  &\: \int_\Omega u \ d\mathbf{x} = V ,\quad \textrm{and}\quad  u=-1\quad \textrm{on}\quad \partial \Omega,\end{aligned}\label{Eq: Allen Cahn 2d}
\end{equation}
where $u:(x,y)\xrightarrow{}\mathbb{R}$ is the phase order parameter, and $V$ is a constant used to enforce the integral of $u$, which is equivalent to the total mass in a two-phase system if $u$ represents rescaled density.

In this case, we use ResNet to represent $u$. Regarding the Lagrange multiplier, we represent $\lambda$ using a shallow neural network with one hidden layer, which includes a small number of neurons with the activation function $\Phi$. We set the input to be $0$, and thus $\lambda(\mathbf{x};\eta)$ can be viewed as a scalar variable $\lambda(0;\eta)$.

As for the boundary conditions of the primal network $\mathbf{u}(\mathbf{x};\theta)$, we enforce them explicitly at the output layer. Specifically, we impose the following boundary conditions at the output layer:
\begin{equation}
u(\mathbf{x}_0;\theta) =  ({W}_{out}f_{N_d}(\mathbf{x}_{N_d}) + {b}_{out})\times(x_0^{(1)})\times(1-x_0^{(1)})\times(x_0^{(2)})\times(1-x_0^{(2)}) - 1,\label{formula: output modification of u-net G-L energy}
\end{equation}
in which $x_0^{(1)}$ and $x_0^{(2)}$ correspond to two coordinates of the input $\mathbf{x}_0$ and $u(\mathbf{x}_0;\theta)=-1$ for $\mathbf{x}_0\in \partial\Omega$. This approach is suitable for relatively regular geometries, such as the unit square we considered in this example. For more general boundary conditions, we can treat them as constraints and use the framework of WANCO to convert them into another adversarial network, as demonstrated in Section~\ref{Subsec: Fluid-solid optimization}.

\subsection{Sensitivity to parameters and comparison to other deep learning based approaches}\label{subsection:GLcomparison}
For problem~\eqref{Eq: Allen Cahn 2d}, we obtain the following minimax problem:
\begin{equation}\label{WAN_Augmented Lagrangian eq_2d parameter version--A-C equation}
    \min_{\theta}\max_{\eta}\: \mathcal{L}^{GL}(\theta,\eta)
\end{equation}
where 
\begin{align*}
\mathcal{L}^{GL}(\theta,\eta) = & C_0\int_{\Omega}\left[\frac{\epsilon}{2}|\nabla u(\mathbf{x};\theta)|^2+\frac{1}{\epsilon}(u(\mathbf{x};\theta)^2-1)^2\right]d\mathbf{x} \\ 
&- \lambda(0;\eta)\left(\int_\Omega u(\mathbf{x};\theta) d\mathbf{x}-V\right) +\frac{\beta}{2}\left(\int_\Omega u(x;\theta) d\mathbf{x}-V\right)^2.
\end{align*}
We then solve this problem with the proposed WANCO in contrast to: 

1. DRM combined with penalty method (DRM-P),  

2. DRM combined with adaptive penalty method(DRM-AP), and

3. DRM combined with Lagrange multipliers.

From the classical analysis of problem~\eqref{Eq: Allen Cahn 2d}, the zero level set of the solution should be a circle with the enclosed area conserved. In this example, we set $\epsilon=0.05$, $V=-0.5$, $C_0=400$, and $\alpha=1.0003$. For the primal network $u(\mathbf{x};\theta)$, the number of hidden layers is $N_d=4$ and the width is $N_w=50$. To represent the adversarial network, we employ a shallow neural network $\lambda(0;\eta)$ with $1$ hidden layer and width $10$ since the Lagrange multiplier $\lambda$ is a scalar in this case. The training step is set to $N=5000$ with inner iteration steps $N_u=N_\lambda=2$. In each iteration, $N_r=40,000$ points are used for training. The initial learning rate is $\tau_\theta=\tau_\eta=0.016$ and halved at certain training milestones. This learning rate setting is consistent throughout the subsequent numerical examples, and will not be reiterated later. The activation function used here is $\tanh^3(x)$. Both DRM-P and DRM-AP have the same structure as the primal network mentioned above.

\begin{figure}[ht!]
	\centering
	\begin{tabular}{|m{1.5cm}|m{3.0cm}|m{3.0cm}|m{3.0cm}|}
	\hline 
		 &   \centering{WANCO} & \centering{DRM-AP} &  \qquad\; DRM-P \\
	 \hline    \footnotesize{$\beta=100000$}  &  \includegraphics[width = 0.22\textwidth]{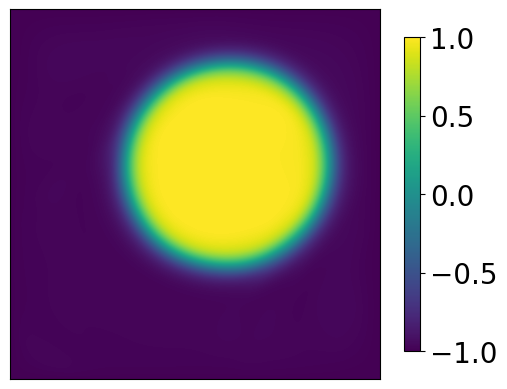}  &  
  \includegraphics[width = 0.22\textwidth]{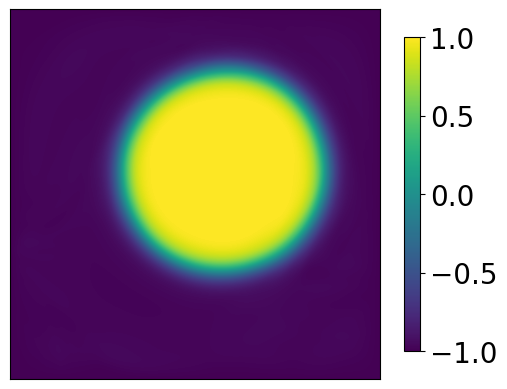} &
  \includegraphics[width = 0.22\textwidth]{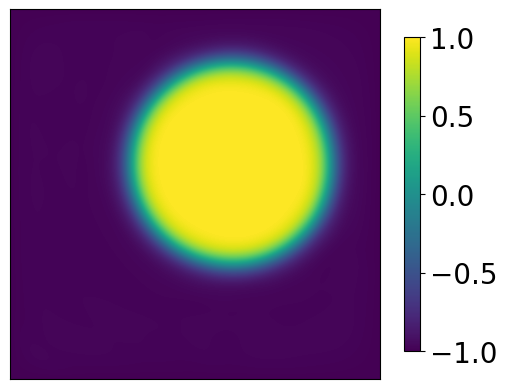}\\
  
    \hline  \footnotesize{$\beta=10000$} & \includegraphics[width = 0.22\textwidth]{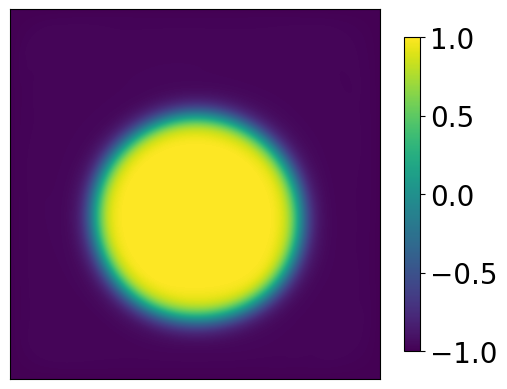}&  
  \includegraphics[width = 0.22\textwidth]{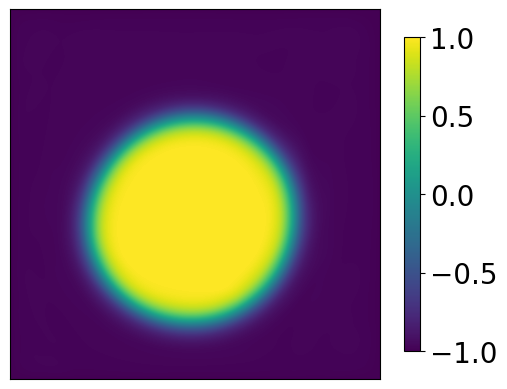}&
  \includegraphics[width = 0.22\textwidth]{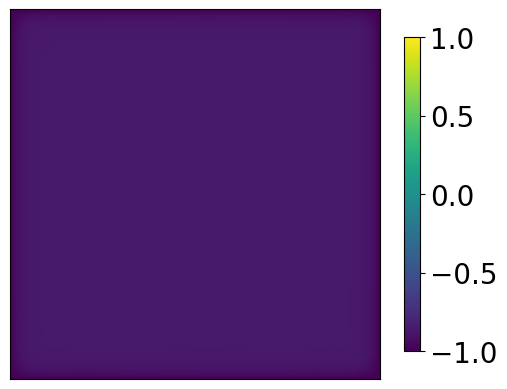}\\
  
    \hline \footnotesize{$\beta=1000$} & \includegraphics[width = 0.22\textwidth]{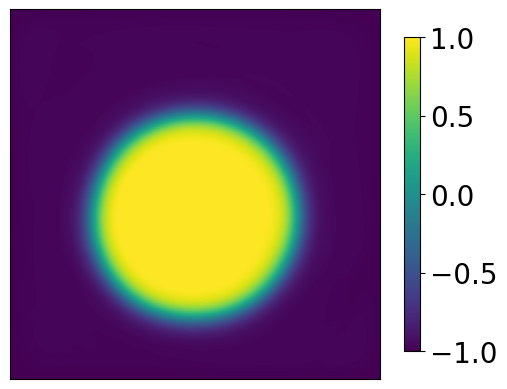}&  
  \includegraphics[width = 0.22\textwidth]{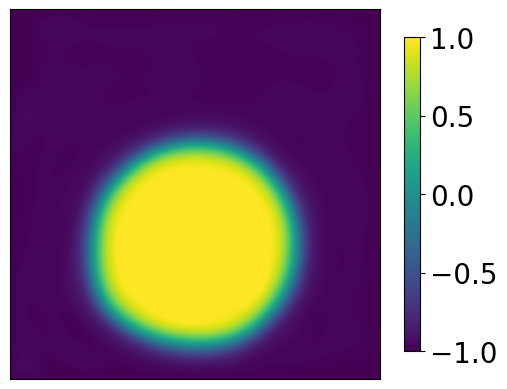}&
  \includegraphics[width = 0.22\textwidth]{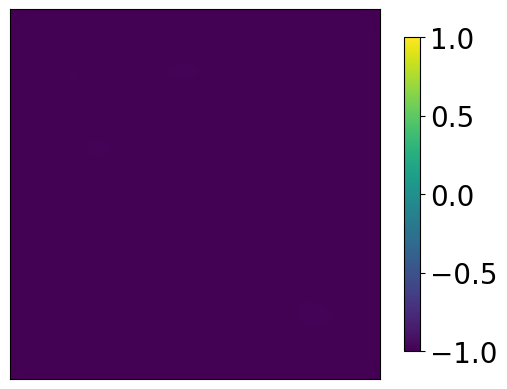}\\
		\hline
	\end{tabular}
\caption{Comparison of the training results among WANCO, DRM-AP, and DRM-P under different $\beta$. See Section~\ref{subsection:GLcomparison}.}
\label{Fig: G-L model  WANCO-AP-P}
\end{figure} 

\begin{figure}[htb!]
\centering
\includegraphics[width = 0.32 \textwidth] {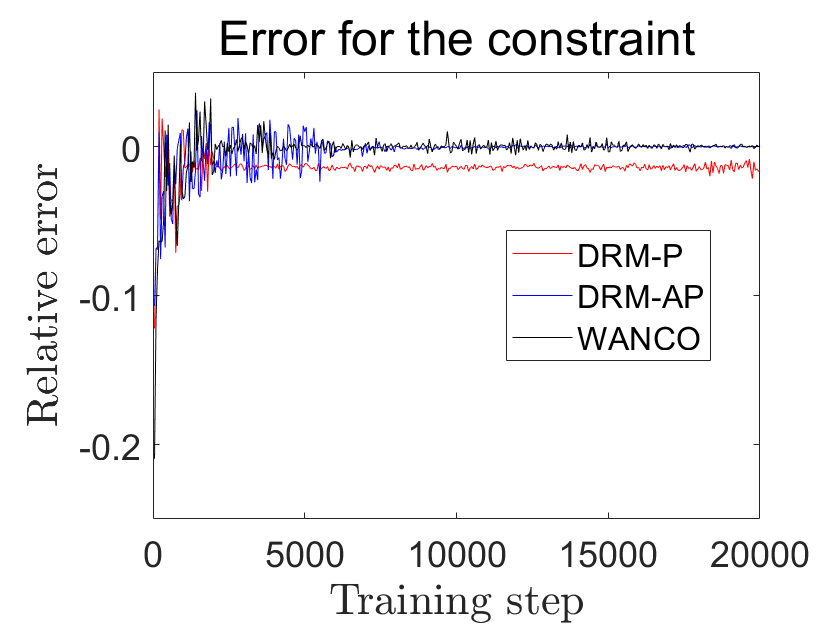}\ \
\includegraphics[width = 0.32 \textwidth]{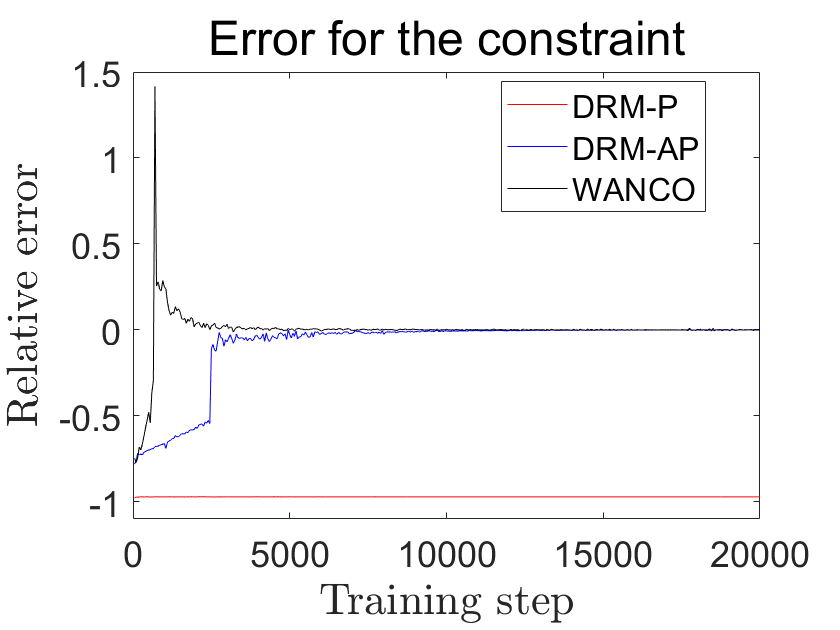}\ \
\includegraphics[width = 0.32 \textwidth]{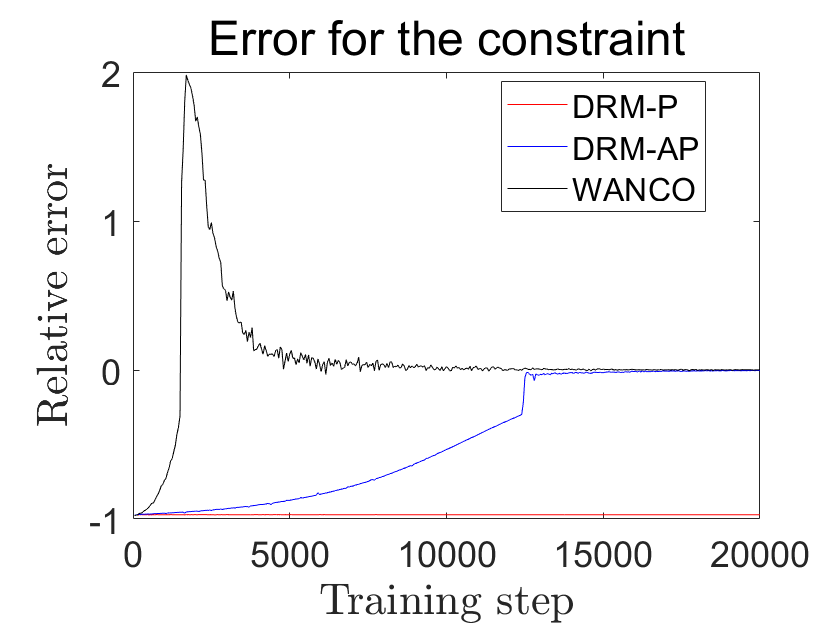}
\caption{The relative error of the constraint for three schemes with corresponding $\beta$. The first, second, and third pictures correspond to $\beta = 100000, 10000, 1000$, respectively. The relative error of the constraint, defined as $\frac{V-\int_\Omega u \ d\mathbf{x} }{V} $, is recorded every $50$ step. See Section~\ref{subsection:GLcomparison}.}
\label{Fig: G-L model error}
\end{figure}

In Figure~\ref{Fig: G-L model  WANCO-AP-P}, we display the results from WANCO, DRM-P, and DRM-AP for different choices of $\beta = 100000, 10000, 1000$, together with the convergence of the constraint error, defined as $\frac{V-\int_\Omega u \ d\mathbf{x}}{V} $, during the training process in Figure~\ref{Fig: G-L model error}. It is worth noting that WANCO employs identical parameters across different $\beta$ values in these results. In contrast, both DRM-P and DRM-AP require individual parameter tuning for each case to obtain satisfactory results (which may still completely fail due to the poor selection of $\beta$).  This implies that the proposed WANCO is relatively insensitive to parameter settings. As shown in Figure~\ref{Fig: G-L model  WANCO-AP-P}, under three different initial penalty terms $\beta$, the proposed WANCO consistently outperforms DRM-P and DRM-AP. For the first case with $\beta=100000$, the profiles of WANCO, DRM-AP, and DRM-P are close to circles, but WANCO and DRM-AP have a clear advantage in the relative error of the constraint. As for the second and third cases with $\beta=10000, 1000$, DRM-P fails, and DRM-AP exhibits relatively inaccurate results (the computed results are not circular) due to the unbalanced terms in the loss, while WANCO trains well. This indicates that relying solely on the adaptive penalty method is still sensitive to parameter selection, and comparable results are achieved only when $\beta=100000$.

In all cases, WANCO can achieve a relative error of the constraint less than $0.005$. As shown in Figure~\ref{Fig: G-L model error}, compared to DRM-AP, WANCO can significantly reduce the training steps due to the adversarial term in~\eqref{WAN_Augmented Lagrangian eq_2d parameter version--A-C equation}. To be more specific, due to the Dirichlet boundary condition and the random initialization of network parameters, the output of the primal network in the initial steps closely approximates $-1$. Consequently, the initial integral of $u$ is also close to $-1$, falling below the constraint of $-0.5$. With the inclusion of the adversarial term $\lambda$, during the early stages of training, the mass of $u$ increases at a faster rate compared to DRM-AP. As the adversarial term $\lambda$ reaches a sufficiently large value, it becomes dominant, causing the minimization step to shift the mass in the positive direction. This corresponds to the initial steps as shown in Figure~\ref{Fig: G-L model error}. Once the error of the constraint becomes positive, $\lambda$ is reduced. Subsequently, when the third term is dominant again, the constraint error decreases rapidly. This corresponds to the stage of rapid descent in Figure~\ref{Fig: G-L model error}, which is more evident in the last two pictures. To conclude, the inclusion of the adversarial term $\lambda$ allows for an adaptive adjustment of its magnitude, enabling the constraint to be rapidly satisfied from both sides. This is in contrast to DRM-AP, where the mass can only decrease rapidly as the penalty term reaches a threshold value. This further demonstrates that WANCO is relatively insensitive to the setting of hyperparameters, as the adversarial term exhibits a notable adaptive weighting and acceleration effect. Moreover, from the training perspective, in WANCO, half of the training steps are on the optimization of the adversarial network, which is usually constructed by a smaller network. This usually results in faster training compared to other methods.

\begin{figure}[ht!]
\centering
\includegraphics[width=0.3\textwidth]{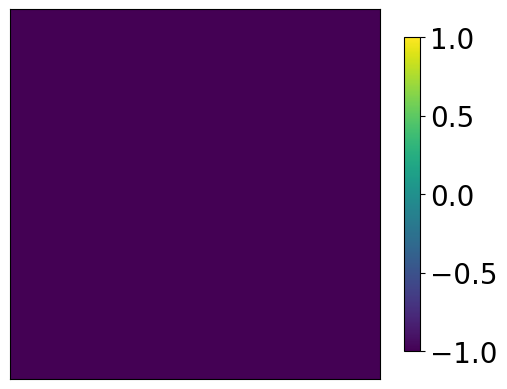} 

\caption{Training solution utilizing adversarial terms only (Lagrange multiplier formula). See Section~\ref{subsection:GLcomparison}.}
\label{Fig: Only adversarial term}
\end{figure}

We also consider the Lagrange multiplier formula without penalty terms, and empirical study shows that it is very difficult to achieve convergence in the adversarial process. As illustrated in Figure~\ref{Fig: Only adversarial term}, using only the Lagrange multiplier method combined with adversarial networks leads to failure. We observe that the Lagrange multiplier only has a corresponding ascent direction during training, and its effectiveness in imposing the constraint is affected by the step size (learning rate). Therefore, relying solely on the Lagrange multiplier may result in convergence difficulties. However, with the augmented Lagrangian method, we observe that the adversarial term adjusts its scale of values based on the current level of constraint satisfaction to accelerate the process of satisfying the constraint, while the penalty term gradually increases to guarantee convergence. As a result, the adversarial process is effectively similar to assigning different weights to the training data, enabling WANCO to impose the constraints more effectively.

\subsection{Choice of activation functions} \label{subsection:activation}

\begin{figure}[htb!]
\centering
\includegraphics[width = 0.32 \textwidth]{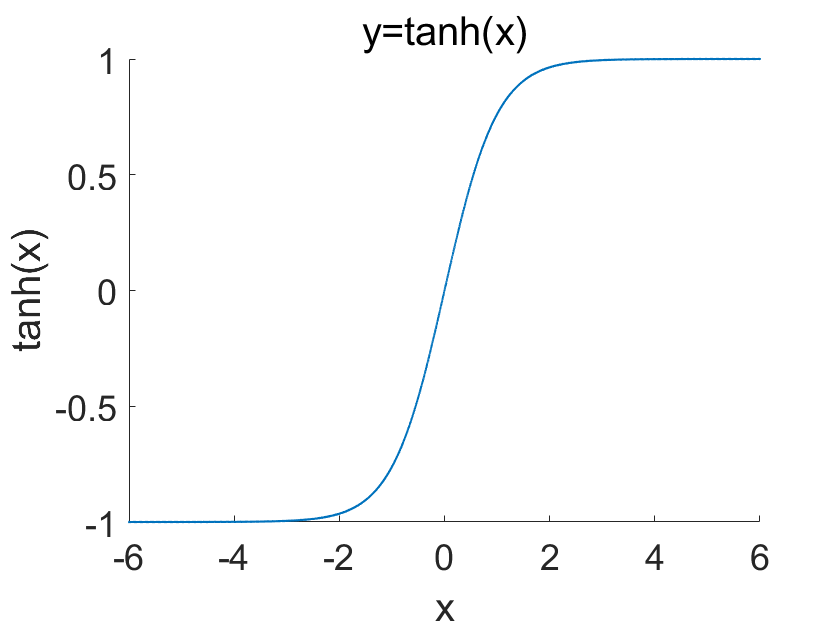} \ \ \ \ \ \ \ \ 
\includegraphics[width = 0.32 \textwidth]{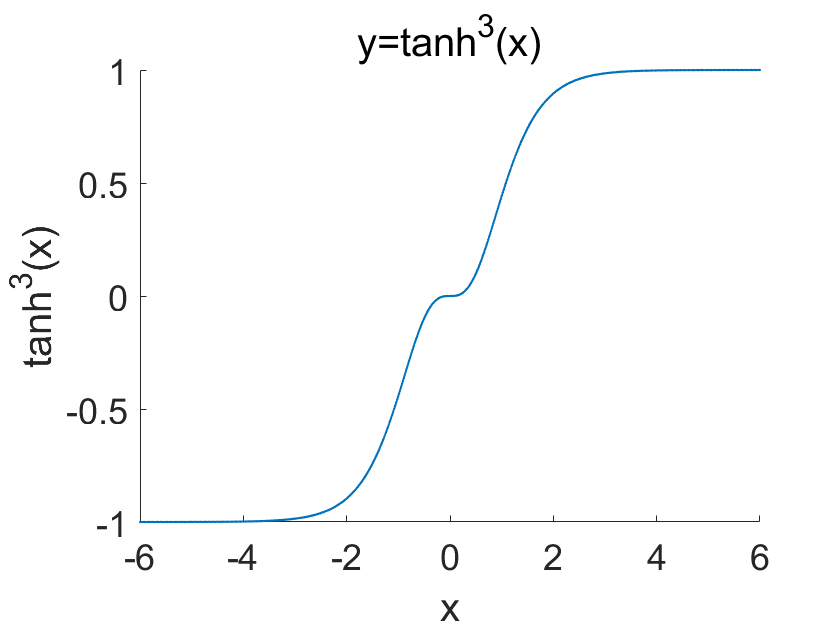}
\caption{Comparison of two activation functions, $\tanh(x)$ and $\tanh^3(x)$, in the range of $[-6,6]$. See Section~\ref{subsection:activation}.}
\label{Fig: activation tanh vs tanh3}
\end{figure}

To illustrate the observations and reasons for the utilization of $\tanh^3$ (see the profile in Figure~\ref{Fig: activation tanh vs tanh3}) in our examples, we present two numerical simulations to validate the advantages of combining ResNet with the activation function $\tanh^3$ in the proposed WANCO. 

For the above problem~\eqref{WAN_Augmented Lagrangian eq_2d parameter version--A-C equation}, we consider: 1. replacing the activation function directly with other activation functions in the model that has been well-tuned for $\tanh^3$, and 2. replacing the activation function with $\tanh^3$ in the model that has been well-tuned for other common activation functions.

\begin{figure}[ht!]
	\centering
	\begin{tabular}{|c|c|c|}
	\hline  $\tanh^3$& $\mathrm{ReLU}$ & $\mathrm{ReLU}^3$ \\ \hline
\includegraphics[width = 0.25 \textwidth]{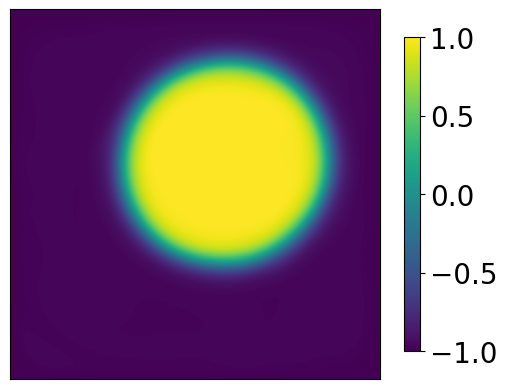} &
\includegraphics[width = 0.25 \textwidth]{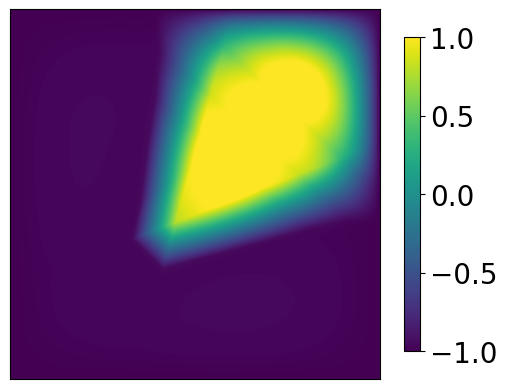} &
\includegraphics[width = 0.25 \textwidth]{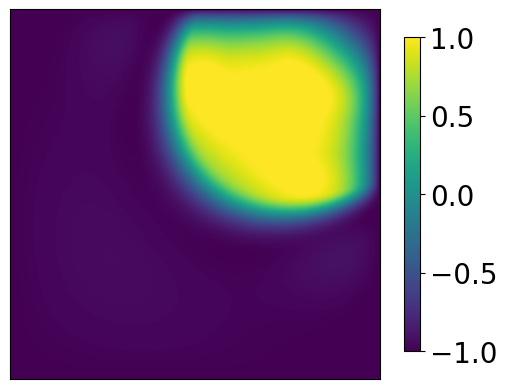}\\ \hline
$\tanh$ & $\mathrm{sigmoid}$ & $\mathrm{softmax}$\\ \hline
\includegraphics[width = 0.25 \textwidth]{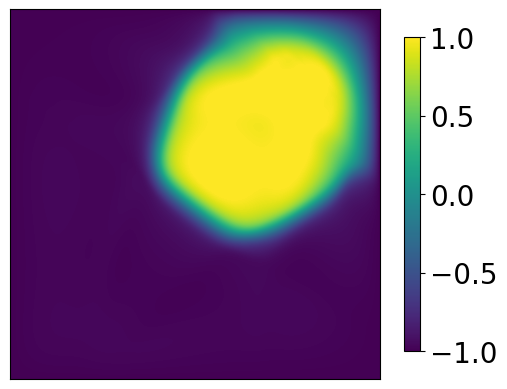}&
\includegraphics[width = 0.25 \textwidth]{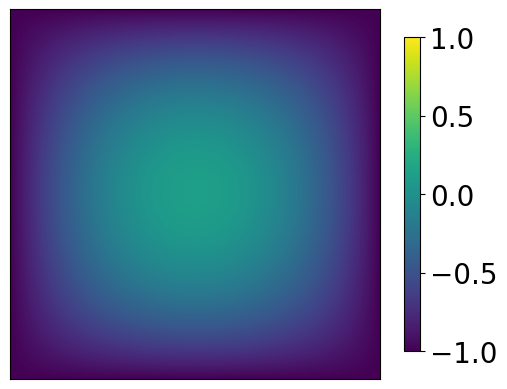} &
\includegraphics[width = 0.25 \textwidth]{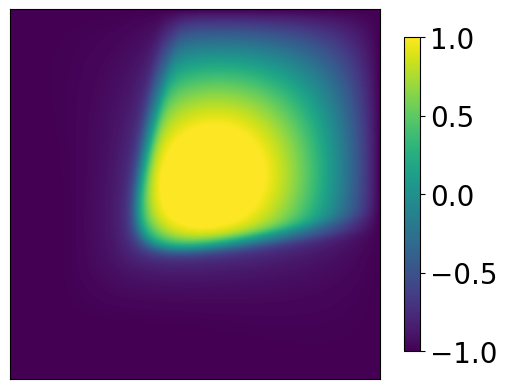} \\ 
 		\hline
	\end{tabular}

\caption{Replacing $\tanh^3$ with other activation functions while keeping the other network parameters unchanged. The training results for the same network parameters with activation functions $\tanh^3$, $\mathrm{ReLU}$, $\mathrm{ReLU}^3$, $\tanh$, $\mathrm{sigmoid}$, and $\mathrm{softmax}$. See Section~\ref{subsection:activation}.}
\label{Fig: Activation test tanh3 to others}
\end{figure} 

\begin{figure}[ht!]
\centering
\begin{tabular}{|m{2.1cm}|m{2.1cm}|m{2.1cm}|m{2.1cm}|m{2.1cm}|}
\hline \centering \footnotesize{$\mathrm{ReLU}$} &
\centering \footnotesize{$\mathrm{ReLU}^3$} &
\centering \footnotesize{$\tanh$} & 
\centering \footnotesize{$\mathrm{sigmoid}$} & \quad\;\;\footnotesize{$\mathrm{softmax}$} \\ \hline
\includegraphics[width = 0.155\textwidth]{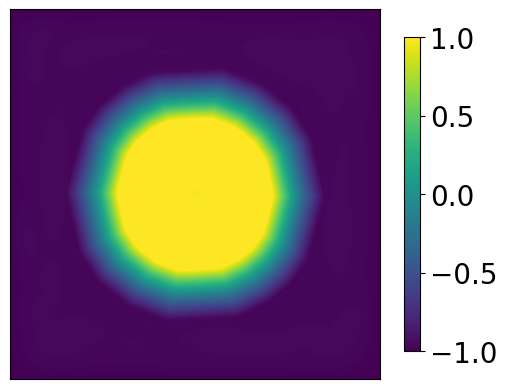} &
\includegraphics[width = 0.155\textwidth]{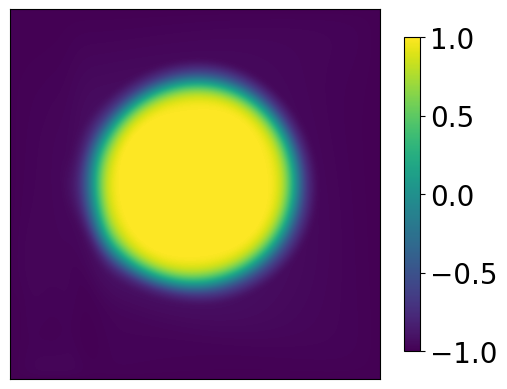} &
\includegraphics[width = 0.155\textwidth]{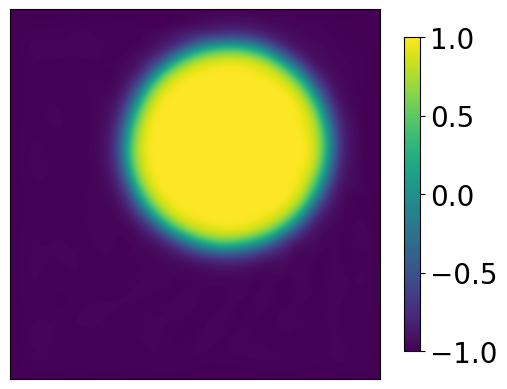} &
\includegraphics[width = 0.155\textwidth]{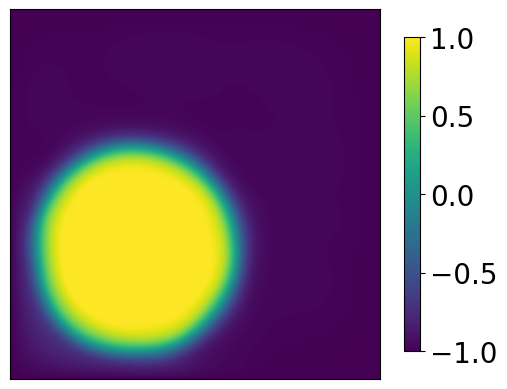} &
\includegraphics[width = 0.155\textwidth]{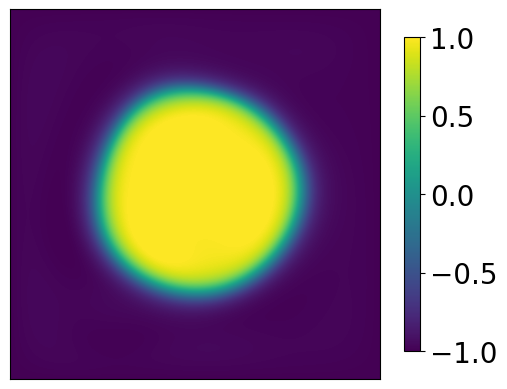} \\ \hline
\centering\footnotesize{$\mathrm{ReLU}$ to $\tanh^3$} & \centering\footnotesize{$\mathrm{ReLU}^3$  to $\tanh^3$}& \centering\footnotesize{$\tanh$  to $\tanh^3$}& \centering\scriptsize{$\mathrm{sigmoid}$  to $\tanh^3$}& \scriptsize{$\mathrm{softmax}$ to $\tanh^3$} \\ \hline
\includegraphics[width = 0.155 \textwidth]{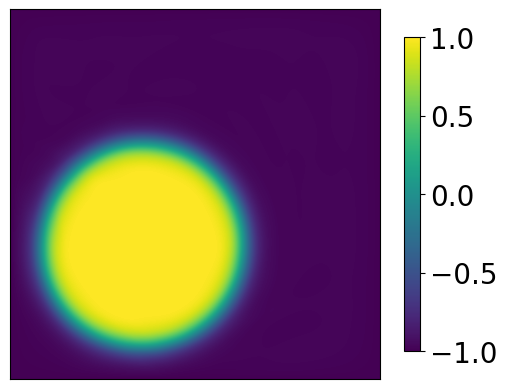} &
\includegraphics[width = 0.155 \textwidth]{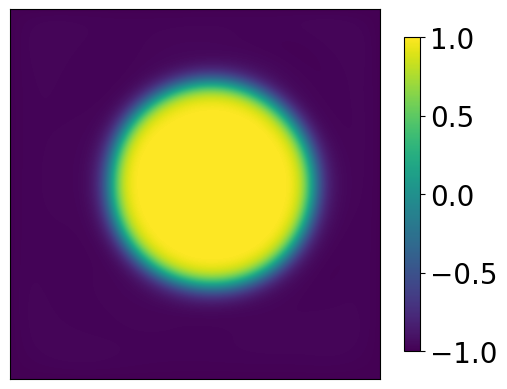}&
\includegraphics[width = 0.155 \textwidth]{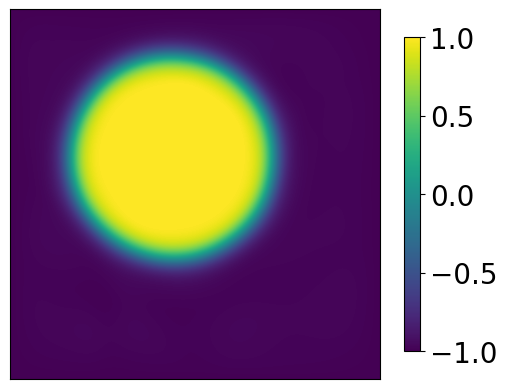} &
\includegraphics[width = 0.155\textwidth]{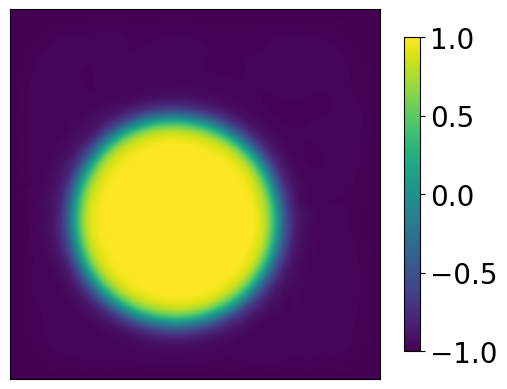} &
\includegraphics[width = 0.155 \textwidth]{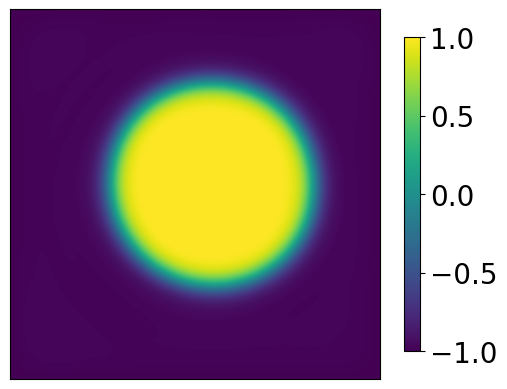} \\
 		\hline
	\end{tabular}
\caption{First row: training results obtained by well-tuned networks with the activation functions $\mathrm{ReLU}$, $\mathrm{ReLU}^3$, $\tanh$, $\mathrm{sigmoid}$, and $\mathrm{softmax}$. Second row: training results after replacing the activation function with $\tanh^3$, while keeping the other parameter settings unchanged throughout the process. See Section~\ref{subsection:activation}.}
\label{Fig: Activation test others to tanh3}
\end{figure}

As shown in Figure~\ref{Fig: Activation test tanh3 to others}, when replacing the activation function $\tanh^3$ in the well-tuned neural network with other activation functions, the training results are often poor or even incorrect. However, as shown in Figure~\ref{Fig: Activation test others to tanh3}, when initially utilizing a well-tuned neural network with a different activation function and subsequently changing it to $\tanh^3$, we observe that similar or even better results can be achieved.

Based on the observations from these two examples, we believe that combining ResNet with $\tanh^3$ exhibits favorable properties and can reduce parameter tuning time for WANCO. This choice is also used and verified in examples for the subsequent problems.

\section{Optimal partition problems}\label{Subsec: Dirichlet partition}

In this section, we consider the optimal partition problem, which minimizes the Dirichlet energy with $L^2$ norm preserving in different dimensions and different boundary conditions \cite{caffarelli2007optimal, Du_2008} as follows.
 \begin{equation}
    \begin{aligned}
                \min_{\mathbf{u}}& \:\int_{\Omega}\left[\frac{\epsilon}{2}\sum_{i}|\nabla u_i|^2 + \frac{1}{\epsilon}f(\mathbf{u})\right]d\mathbf{x}\\
      \st \: & \int_\Omega u_i^2 \ d\mathbf{x} = 1, \qquad i = 1,2,...,n, \\
     & \text{and}\ \mathbf{u} \ \text{satisfies the homogeneous Dirichlet} \\ &\text{or periodic boundary conditions},
    \end{aligned} \label{eq:Dirichlet1}
\end{equation}
where $\Omega=[0,1]^d$, $d$ denotes the dimension, $n$ represents the number of phases, $\mathbf{u} = [u_1,u_2,...,u_n]'$ denotes the phase variables, with each $u_i$ mapping $\mathbf{x}$ to $\mathbb{R}_+$, $\epsilon$ is the model parameter related to the thickness of the transition layer, and $f$ takes the form $f(\mathbf{u}) = \sum_{i, j=1, i\neq j}^n u_i^2u_j^2$.

Following the augmented Lagrangian method, one can obtain the minimax problem:
\begin{equation}\label{WAN_Augmented Lagrangian eq_2d}
\begin{aligned}
    \min_{\mathbf{u}}\max_{\boldsymbol{\lambda}}\: C_0\int_{\Omega}\left[\frac{\epsilon}{2}\sum_{i}|\nabla u_i|^2 + \frac{1}{\epsilon}\sum_{i=1, j=1, i\neq j}^nu_i^2u_j^2\right]d\mathbf{x} \\ 
    -\boldsymbol{\lambda}\cdot\left(\int_\Omega\mathbf{u^2}d\mathbf{x}-\mathbf{1}\right)+\frac{\beta}{2}\left\|\int_\Omega\mathbf{u^2}d\mathbf{x}-\mathbf{1}\right\|_2^2,
\end{aligned}
\end{equation}
where $\boldsymbol{\lambda}$ denotes $n$ Lagrange multipliers for $n$ constraints in \eqref{eq:Dirichlet1}, and $\mathbf{1}$ is the vector in $\mathbb{R}^n$ whose component are all ones. 

Then we construct the primal and adversarial networks to represent $\mathbf{u} = \mathbf{u}(\mathbf{x};\theta)$ and $\boldsymbol{\lambda} =  \boldsymbol{\lambda}(\mathbf{x};\eta)$, respectively. Since the $L^2$ norm preserving constraint is independent of $\mathbf{x}$, we construct the adversarial network as $\boldsymbol{\lambda}(0;\eta)=[\lambda_1(0;\eta),\lambda_2(0;\eta),...,\lambda_n(0;\eta)]'$ using a shallow neural network. Consequently, the problem~\eqref{WAN_Augmented Lagrangian eq_2d} is transformed into the following form:
\begin{equation}\label{WAN_Augmented Lagrangian eq_2d parameter version Dirichlet partition}
\begin{aligned}
           \min_{\theta} \max_{\eta}\:& L_\beta(\mathbf{u}(\mathbf{x};\theta),\boldsymbol{\lambda}(0;\eta))\\
           &= C_0\int_{\Omega}\left[\frac{\epsilon}{2}\sum_{i=1}^n|\nabla u_i(\mathbf{x};\theta)|^2 + \frac{1}{\epsilon}\sum_{i=1, j=1, i\neq j}^n u_i^2(\mathbf{x};\theta)u_j^2(\mathbf{x};\theta)\right]d\mathbf{x}\\
           &-\left[\boldsymbol{\lambda}(\mathbf{0};\eta)\cdot\left(\int_\Omega \mathbf{u}^2(\mathbf{x};\theta)d\mathbf{x}-\mathbf{1}\right)\right]+\frac{\beta}{2}\left\|\int_\Omega\mathbf{u}^2(\mathbf{x};\theta)d\mathbf{x}-\mathbf{1}\right\|_2^2. 
\end{aligned}
\end{equation}

In the follows, we consider the results with the Dirichlet and periodic boundary conditions separately. The parameter settings and training results for the Dirichlet and periodic boundary conditions are given separately. These results in $2$- and $3$-dimensional space are consistent with the results presented in \cite{Wang_2022_DP,Bogosel_2016,Wang_2019}, indicating the effectiveness and robustness of the proposed method. In addition, the $4$-dimensional partition results are displayed to demonstrate the efficiency of the proposed WANCO. 

\subsection{Dirichlet partition with Dirichlet boundary conditions}\label{subsec:Dirichletboundary}

In this section, we consider the Dirichlet partition problem with Dirichlet boundary conditions in the unit square $\Omega=[0,1]\times[0,1]$ with the $L^2$ norm preserved. The Dirichlet boundary condition is imposed as a strong form at the output layer similar to the first numerical example,
\begin{equation}\label{formula: output regularization of u-net}
\mathbf{u}(\mathbf{x}_0;\theta) = \mathrm{ReLU}(\mathbf{W}_{out}f_{N_d}(\mathbf{x}_{N_d}) + \mathbf{b}_{out})\times\Pi_{i=1}^d \left[(x_0^{(i)})\times(1-x_0^{(i)})\right],
\end{equation}
in which $x_0^{(i)}$ corresponds to the $i$-th entry of the input $\mathbf{x}_0\in \Omega$ and $\mathrm{ReLU}$ here is used to enforce the non-negativity of $\mathbf{u}$. 

The parameters for the primal networks are set as follows: $\epsilon=0.05$, $C_0=100$, $\beta=10000$, $\alpha=1.0003$, $N_d=8$, and $N_w=120$. A shallow neural network $\boldsymbol{\lambda}(0;\eta)$ with $1$ hidden layer and width $10$ is used to represent the adversarial network. The training is performed for $N=20000$ steps with inner iteration steps $N_u=N_\lambda=1$, and $N_r=10000$ points are used for training in each iteration. The testing results are obtained by evaluating the trained neural network on a uniform $1000\times 1000$ grid.

\begin{center}
\begin{figure}[h!]
\centering
\includegraphics[width = 0.18 \textwidth]{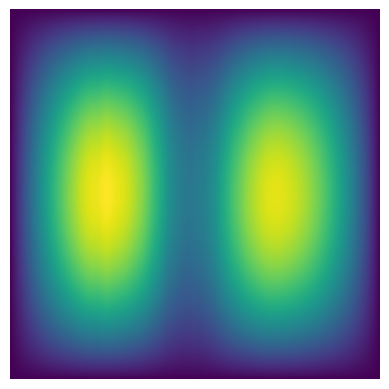}\ \ \
\includegraphics[width = 0.18 \textwidth]{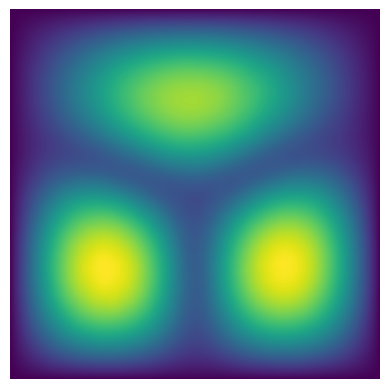}\ \ \
\includegraphics[width = 0.18 \textwidth] {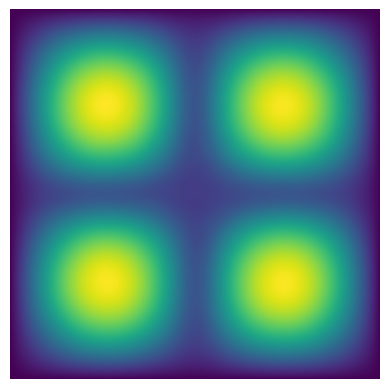}\ \ \
\includegraphics[width = 0.18 \textwidth] {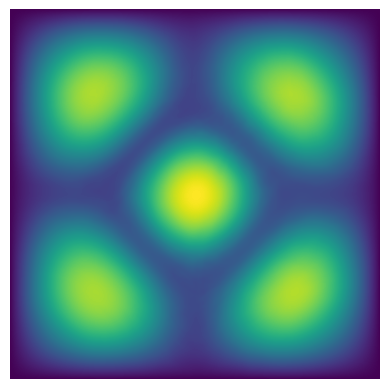}\\
\includegraphics[width = 0.18 \textwidth]{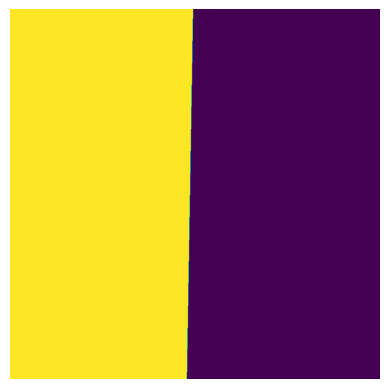}\ \ \
\includegraphics[width = 0.18 \textwidth]{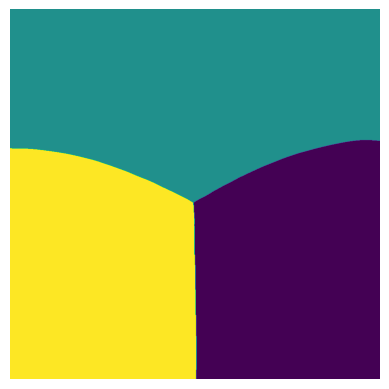}\ \ \
\includegraphics[width = 0.18 \textwidth]{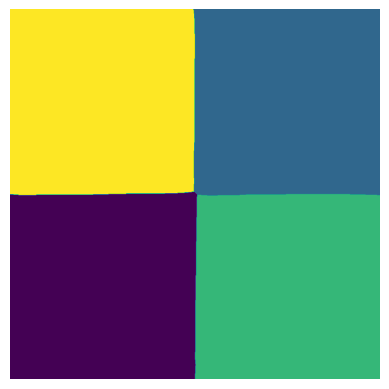}\ \ \
\includegraphics[width = 0.18 \textwidth] {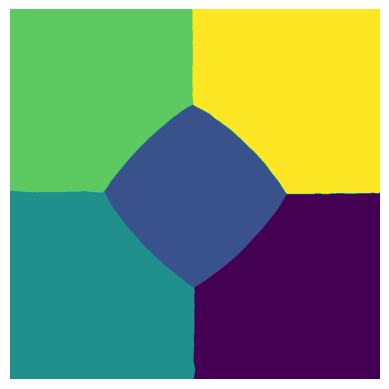}\\
\includegraphics[width = 0.18 \textwidth]{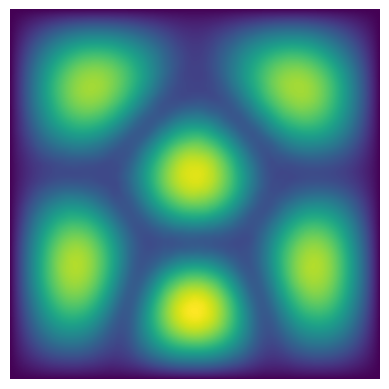}\ \ \
\includegraphics[width = 0.18 \textwidth]{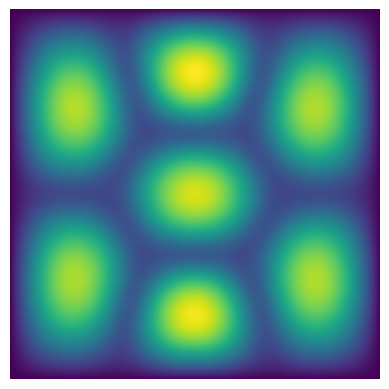}\ \ \
\includegraphics[width = 0.18 \textwidth]{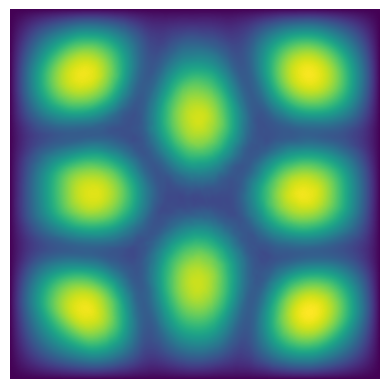}\ \ \
\includegraphics[width = 0.18 \textwidth]{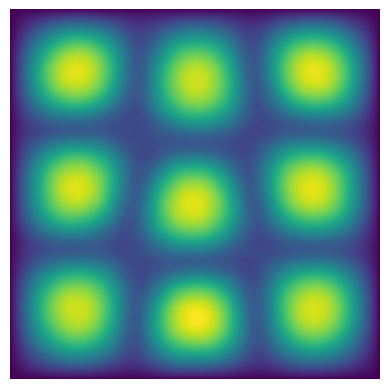}\\
\includegraphics[width = 0.18 \textwidth]{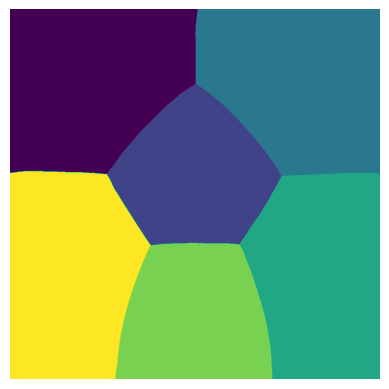}\ \ \
\includegraphics[width = 0.18 \textwidth]{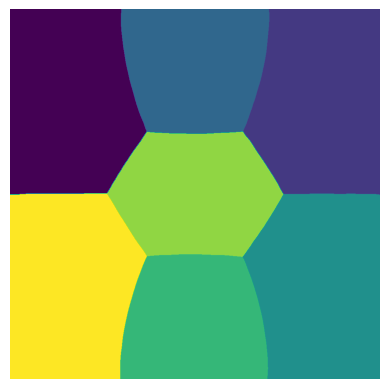}\ \ \ 
\includegraphics[width = 0.18 \textwidth]{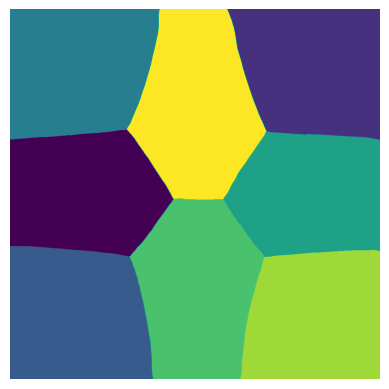}\ \ \
\includegraphics[width = 0.18 \textwidth]{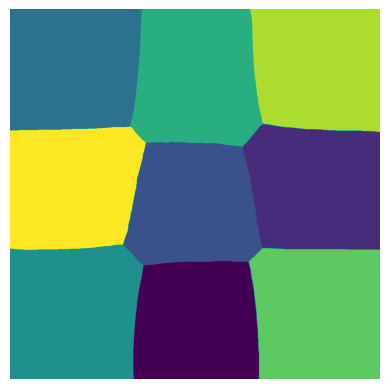}\\
\caption{Results for the Dirichlet partition problem with Dirichlet boundary conditions in the unit square $[0,1]\times[0.1]$. First and third rows: the training solution $\sum_{i=1}^n(u_i(\mathbf{x};\theta))$ for $n=2$-$9$. Second and fourth rows: the corresponding partitions after the projection \eqref{formula: projection 2d n phase} for $2$-$9$ partitions. These results are based on the evaluation of the trained neural network on $1000\times 1000$ grid points. See Section~\ref{subsec:Dirichletboundary}.}
\label{Fig: Dirichlet Partition 2d unit square 2-9 phase}
\end{figure} 
\end{center}

In Figure~\ref{Fig: Dirichlet Partition 2d unit square 2-9 phase}, we display the optimal $n$-partitions with Dirichlet boundary conditions, which are consistent with the results reported in \cite{Wang_2022_DP, Bogosel_2016}, computed with traditional numerical methods. The first and third rows of Figure~\ref{Fig: Dirichlet Partition 2d unit square 2-9 phase} represent the network output after being well trained. The second and fourth rows represent the plot of $n$ partitions based on the following projection,
\begin{equation}\label{formula: projection 2d n phase}
\phi(\mathbf{x}) = m, \quad m=\mathop{\arg\max}\limits_{1\leq k\leq n} u_k(\mathbf{x};\theta), \quad \text{for}\quad \mathbf{x} \in \Omega.
\end{equation}
In the traditional method for this optimization problem, one has to carefully handle the boundary condition by employing concave relaxation or Schr{\"o}dinger operator relaxation and design the optimization scheme to obtain reasonable and regular results. In WANCO, the aforementioned parameters can effectively train the network for all numbers (i.e., $n=2$-$9$) of partitions, which correspond to different dimensions of the network output. This implies the robustness of WANCO to the parameter setting.

\subsection{Dirichlet partition with periodic boundary conditions}\label{sec:periodic}

In this section, we consider the Dirichlet partition problem in $2$-, $3$-, and $4$-dimensional flat tori. The periodic boundary conditions on $\Omega=[0,1]^d$ are imposed by transforming the input layer as follows~\cite{Lyu_2021},
\begin{equation}
\begin{aligned}
f_{in}(\mathbf{x}_0) &= f((x_0^{(1)},x_0^{(2)},...,x_0^{(d)})')\\ = (\cos(2\pi x_0^{(1)}),\sin(2\pi x_0^{(1)}),&\cos(2\pi x_0^{(2)}),\sin(2\pi x_0^{(2)}),...,\cos(2\pi x_0^{(d)}),\sin(2\pi x_0^{(d)}))',
\end{aligned}
\end{equation}
in which $x_0^{(i)}$ corresponds to the $i$-th entry of the input $\mathbf{x}_0\in \Omega$.

The parameters for the primal networks are set as follows: $\epsilon=0.04$, $C_0=2500$, $\beta=100000$, $\alpha=1.0003$,  $N_d=3$, and the width is increasing from $N_w=50$ to $N_w=200$ as the dimension or the number of phases increases. A shallow neural network $\boldsymbol{\lambda}(0;\eta)$ with $1$ hidden layer and width $10$ is used to represent the adversarial network. The training is performed for $N=5000$ steps with inner iteration steps $N_u=N_\lambda=2$. $N_r=10000, 40000, 60000$ points are used for $2$-, $3$-, and $4$-dimensional cases in each iteration, respectively. The sampling points used are fixed as the Hammersley sequence in the corresponding dimension, similar to that mentioned in~\cite{Wu_2023}.  

\begin{figure}[t!]
\centering
\includegraphics[width = 0.18 \textwidth]{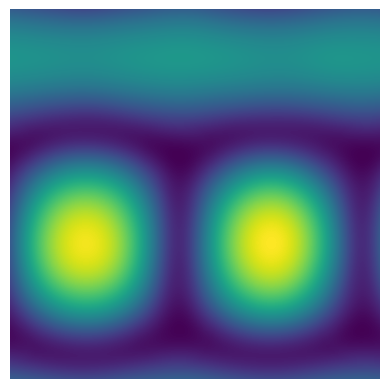}\ \
\includegraphics[width = 0.18 \textwidth]{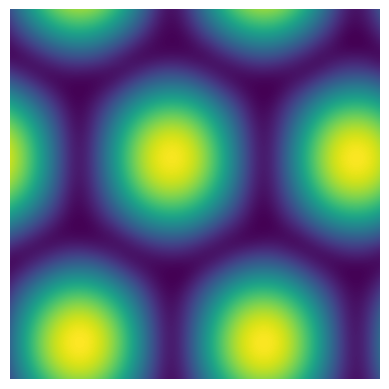}\ \
\includegraphics[width = 0.18 \textwidth] {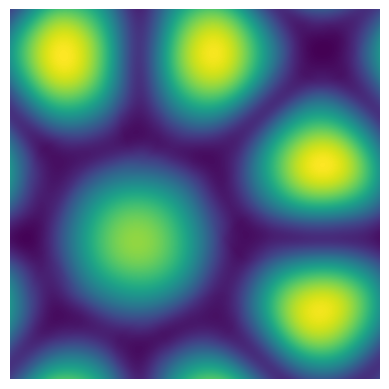}\ \
\includegraphics[width = 0.18 \textwidth] {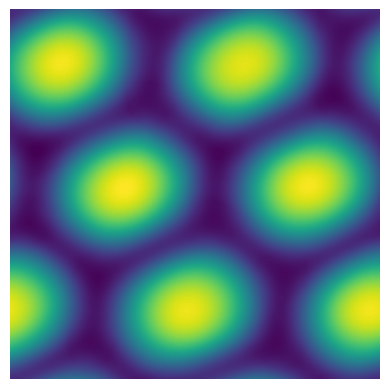}\ \
\includegraphics[width = 0.18 \textwidth]{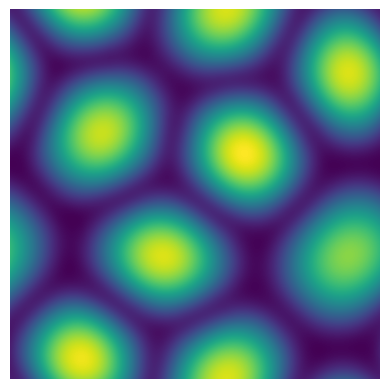}\\
\includegraphics[width = 0.18 \textwidth]{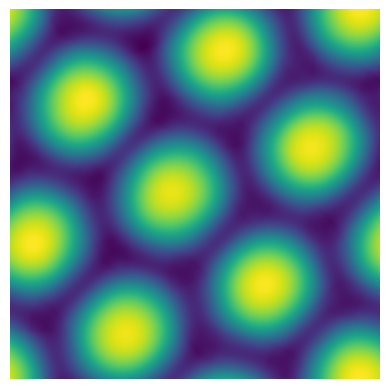}\ \
\includegraphics[width = 0.18 \textwidth]{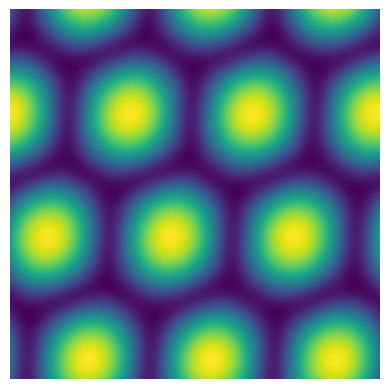}\ \
\includegraphics[width = 0.18 \textwidth] {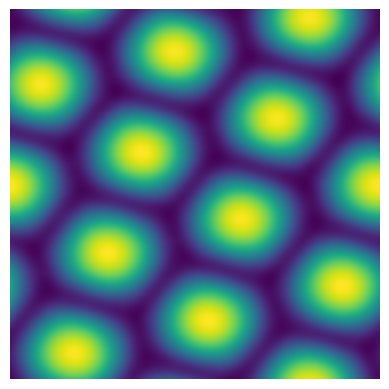}\ \
\includegraphics[width = 0.18 \textwidth] {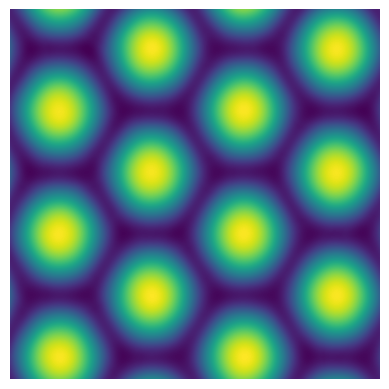}\ \
\includegraphics[width = 0.18 \textwidth]{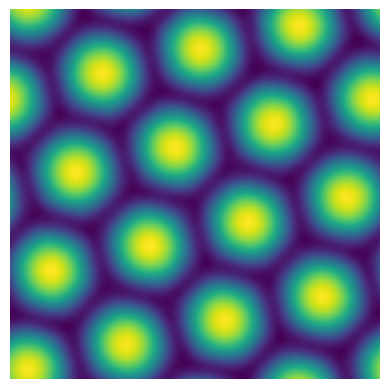}\\
\includegraphics[width = 0.18 \textwidth]{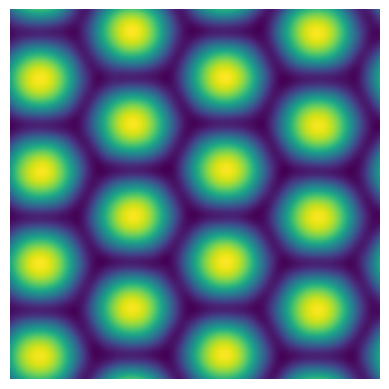}\ \
\includegraphics[width = 0.18 \textwidth]{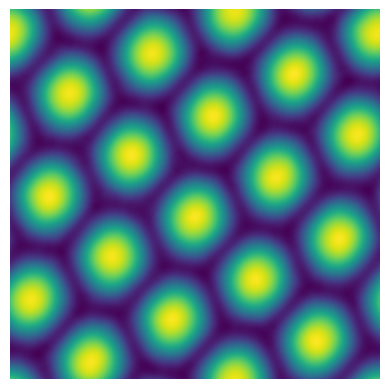}\ \
\includegraphics[width = 0.18 \textwidth] {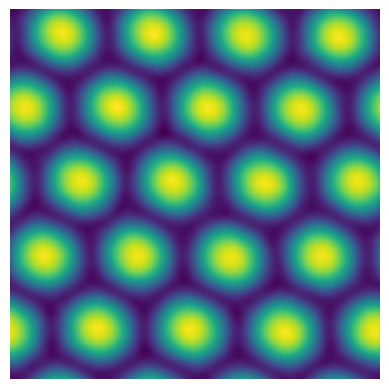}\ \
\includegraphics[width = 0.18 \textwidth] {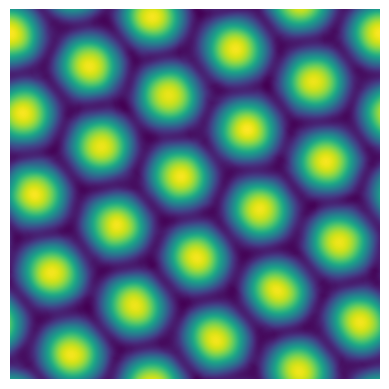}\ \
\includegraphics[width = 0.18 \textwidth]{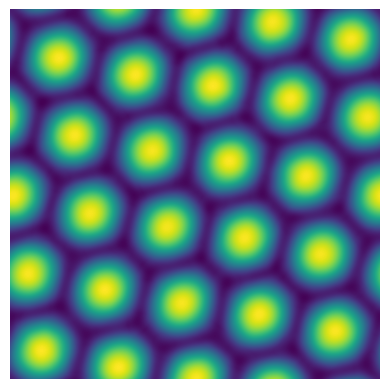}\\
\caption{Trained results for the Dirichlet partition problem with periodic boundary conditions in the $2$-dimensional flat torus for $n=3-9, 11, 12,15,16,18,20,23$, and $24$. See Section~\ref{sec:periodic}.}
\label{Fig: Dirichlet Partition 2d unit square 3-8 phase(periodic bdry)}
\end{figure}

\begin{figure}[h!]
\centering
\includegraphics[width = 0.18 \textwidth]{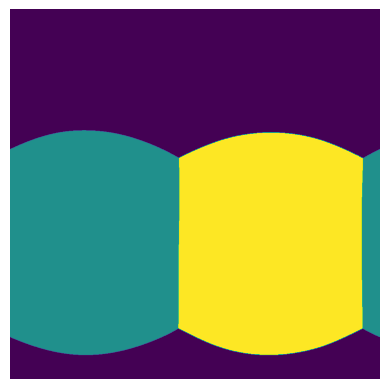}\ \
\includegraphics[width = 0.18 \textwidth]{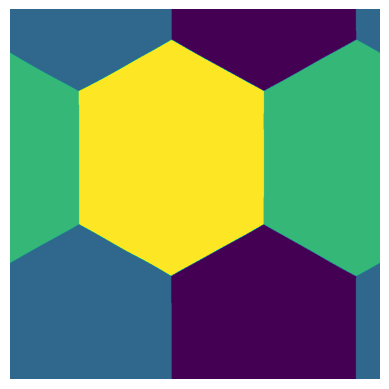}\ \
\includegraphics[width = 0.18 \textwidth] {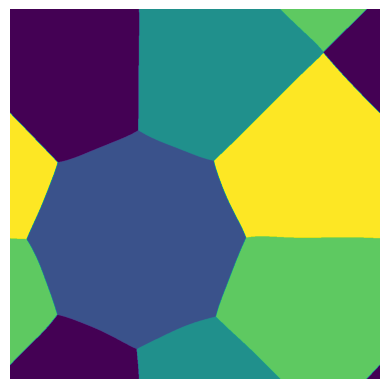}\ \
\includegraphics[width = 0.18 \textwidth] {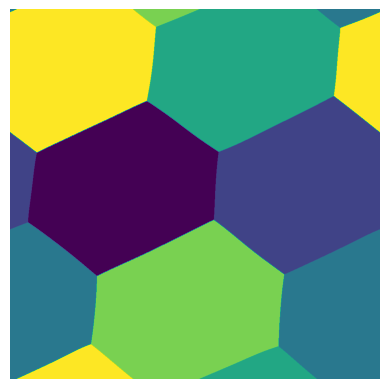}\ \
\includegraphics[width = 0.18 \textwidth]{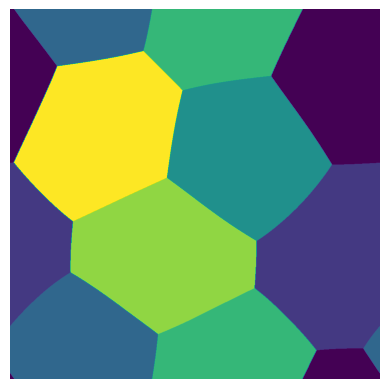} \\
\includegraphics[width = 0.18 \textwidth]{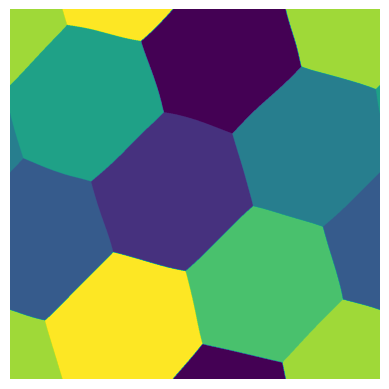}\ \
\includegraphics[width = 0.18 \textwidth]{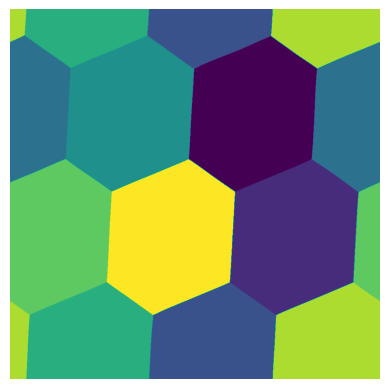}\ \
\includegraphics[width = 0.18 \textwidth] {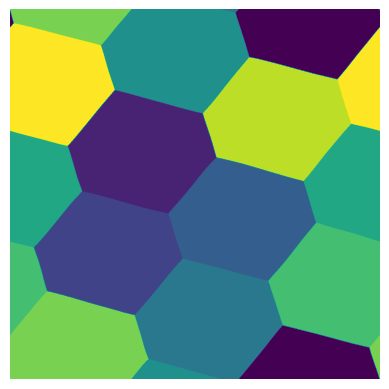}\ \
\includegraphics[width = 0.18 \textwidth] {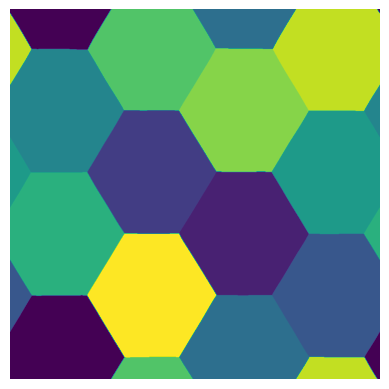}\ \
\includegraphics[width = 0.18 \textwidth]{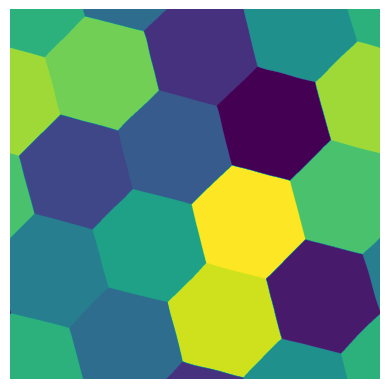} \\
\includegraphics[width = 0.18 \textwidth]{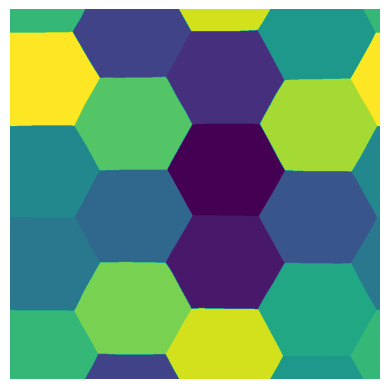}\ \
\includegraphics[width = 0.18 \textwidth]{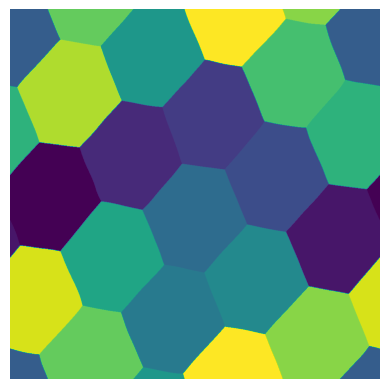}\ \
\includegraphics[width = 0.18 \textwidth] {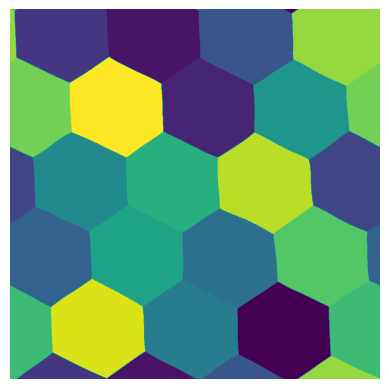}\ \
\includegraphics[width = 0.18 \textwidth] {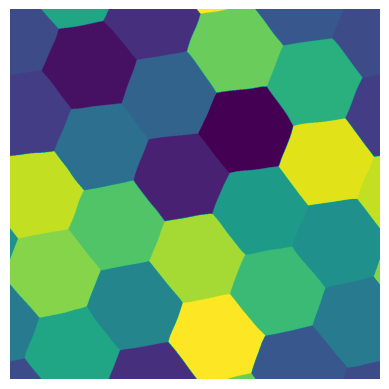}\ \
\includegraphics[width = 0.18 \textwidth]{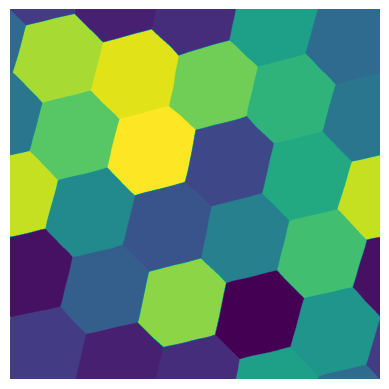}\\
\caption{The projected partitions for the Dirichlet partition problem with periodic boundary conditions in the $2$-dimensional flat torus for $n=3-9, 11, 12,15,16,18,20,23$, and $24$. See Section~\ref{sec:periodic}.}
\label{Fig: Dirichlet Partition 2d unit square 3-8 phase(periodic bdry)2}
\end{figure}

{\bf Two-dimensional flat torus:} In Figures~\ref{Fig: Dirichlet Partition 2d unit square 3-8 phase(periodic bdry)} and \ref{Fig: Dirichlet Partition 2d unit square 3-8 phase(periodic bdry)2}, we display the results of $n$-partition in a $2$-dimensional flat torus trained by the proposed WANCO for $n=3-9, 11, 12,15,16,18,20,23$, and $24$. Again, the results in Figure~\ref{Fig: Dirichlet Partition 2d unit square 3-8 phase(periodic bdry)2} are projections using~\eqref{formula: projection 2d n phase} based on the trained results in Figure~\ref{Fig: Dirichlet Partition 2d unit square 3-8 phase(periodic bdry)}. The testing results are obtained by evaluating the trained neural network on a uniform $1000\times 1000$ grid. To our knowledge, even in the $2$-dimensional case, there are not too many efficient numerical methods that can compute the $n$-partition problems for a large $n$ (\eg, $n=24$). Existing methods \cite{Cheng_2020, Chu_2021} usually require thousands of steps to converge for the problem with a large number of partitions, even in the $2$-dimensional case. In all of the previous studies, hexagons are ubiquitous for large values of $n$. Figures~\ref{Fig: Dirichlet Partition 2d unit square 3-8 phase(periodic bdry)} and \ref{Fig: Dirichlet Partition 2d unit square 3-8 phase(periodic bdry)2} show that WANCO obtains consistent results for a large range of $n$. Since the domain has an aspect ratio equal to one, regular hexagons cannot be used to tile the domain, so the hexagons are slightly distorted. 

\begin{figure}[ht!]
  \centering
    \includegraphics[width=0.18\textwidth, clip, trim = 7cm 5cm 7cm 4.7cm]{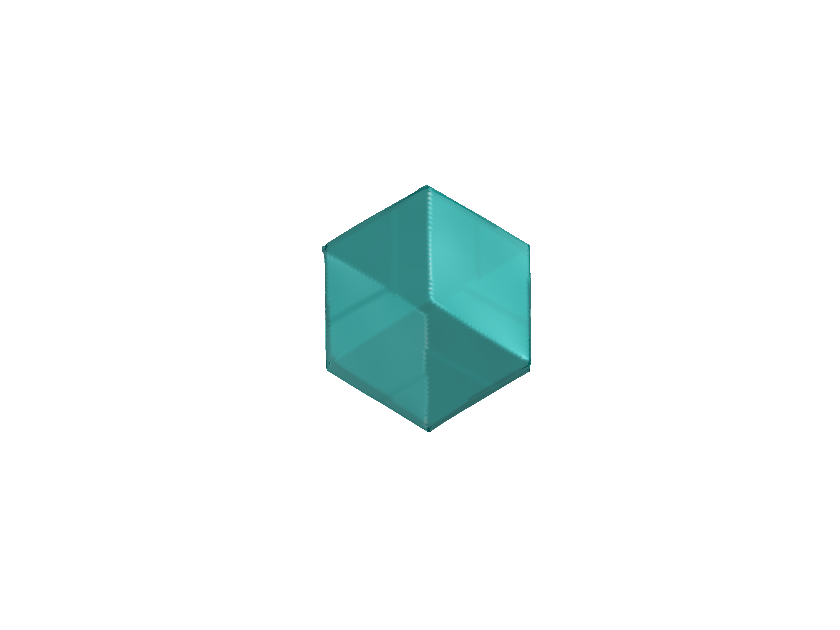}
    \includegraphics[width=0.16\textwidth, clip, trim = 7cm 4cm 5.5cm 3cm]{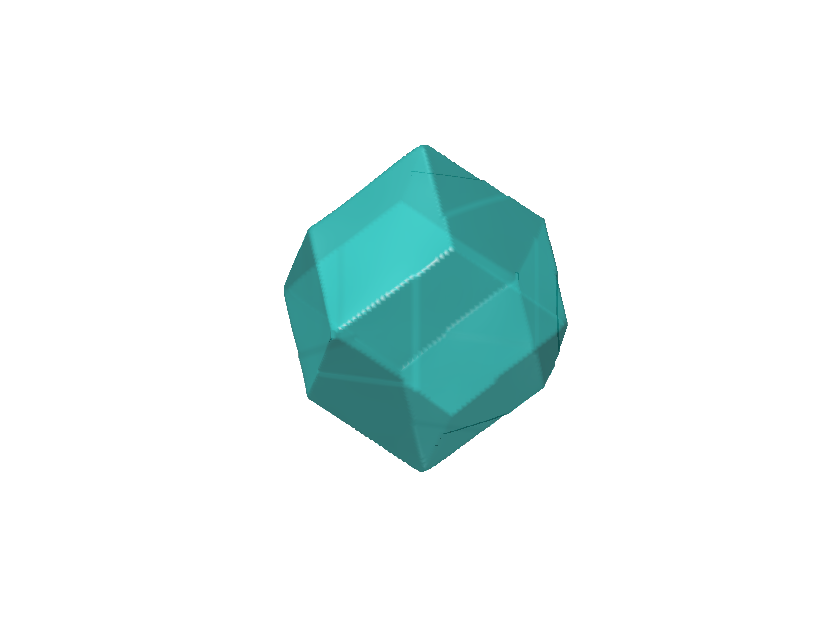}
    \includegraphics[width=0.16\textwidth, clip, trim = 9cm 5.5cm 7cm 5cm]{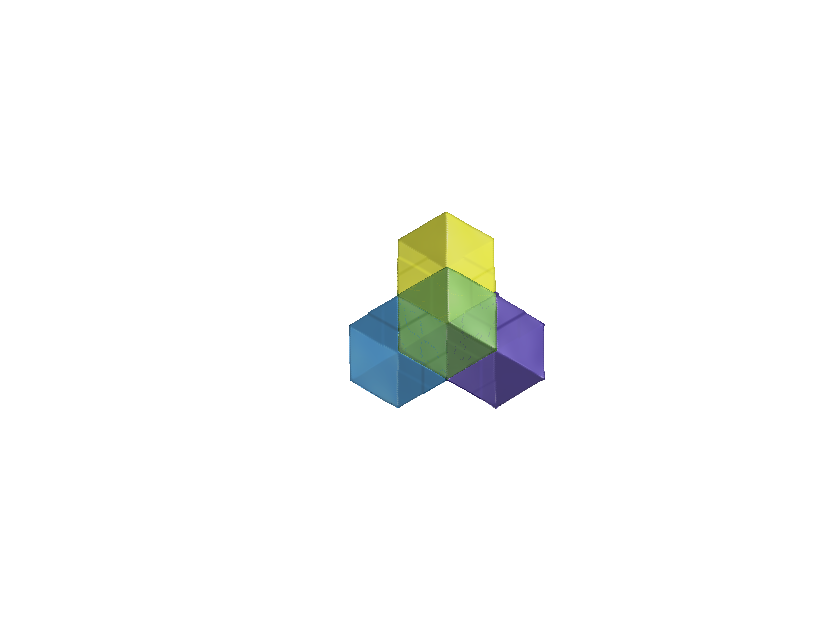}
    \includegraphics[width=0.16\textwidth, clip, trim = 8cm 5cm 7cm 5cm]{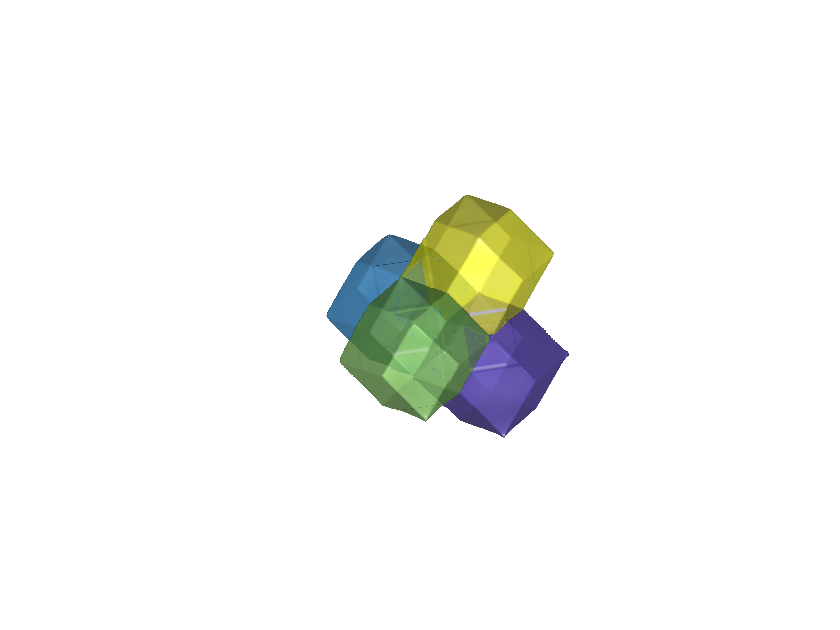}\\
\caption{Partition results for the Dirichlet partition problem with periodic boundary conditions in the $3$-dimensional flat torus with $n=4$. See Section~\ref{sec:periodic}.}
\label{Fig: D-P P Model 3d n4}
\end{figure}

\begin{figure}[ht!]
\centering
\includegraphics[width=0.16\textwidth, clip, trim = 7cm 6cm 7cm 3cm]{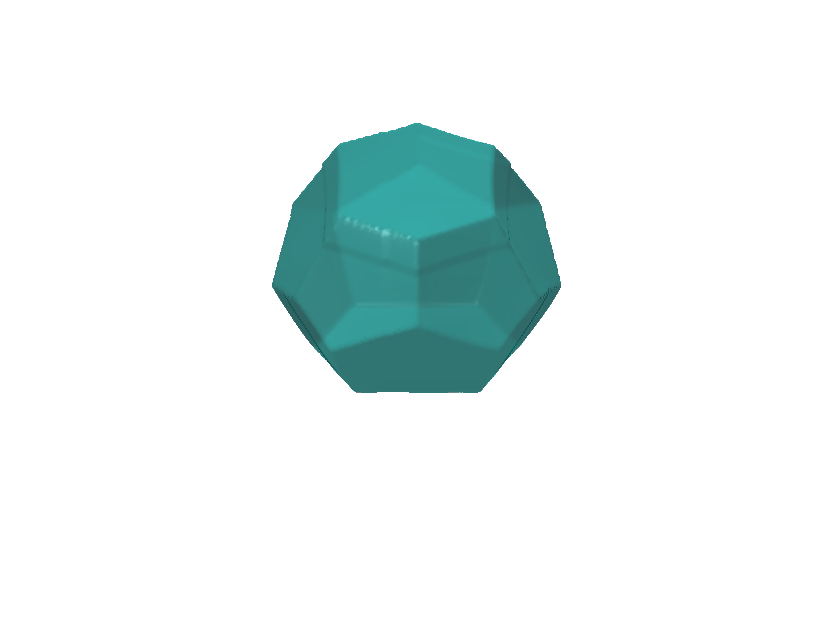}
\includegraphics[width=0.16\textwidth, clip, trim = 5cm 4cm 6cm 2.5cm]{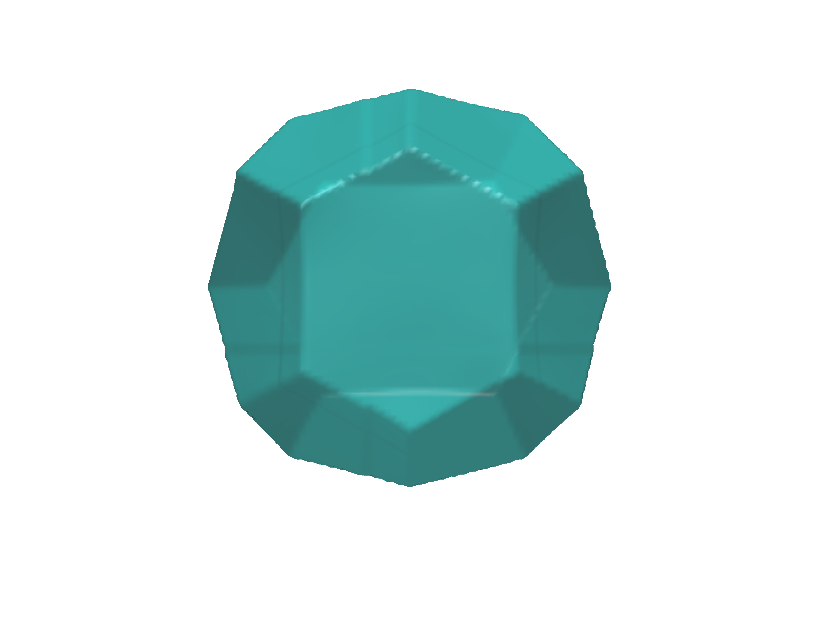}\ 
\includegraphics[width=0.18\textwidth, clip, trim = 5cm 3.5cm 5cm 3cm]{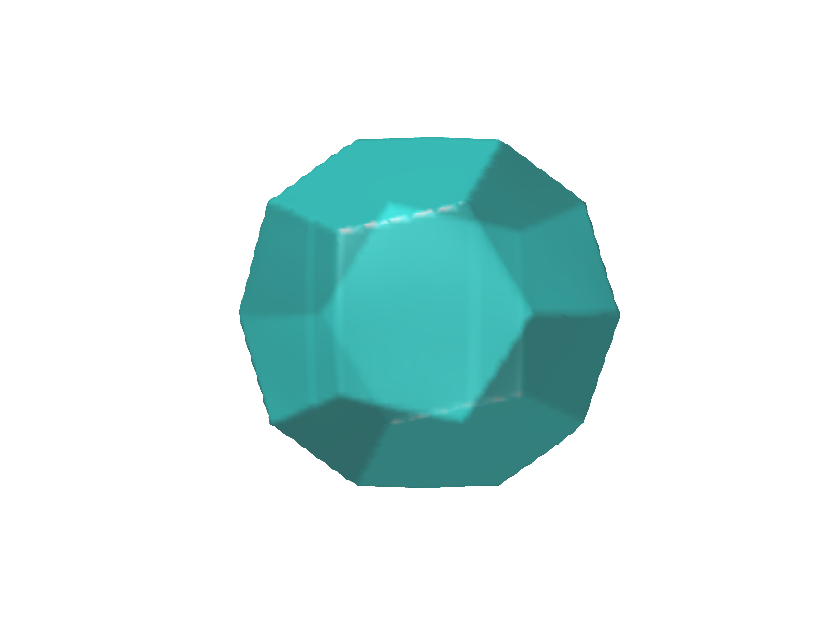}
\includegraphics[width=0.16\textwidth, clip, trim = 5cm 3cm 5cm 3cm]{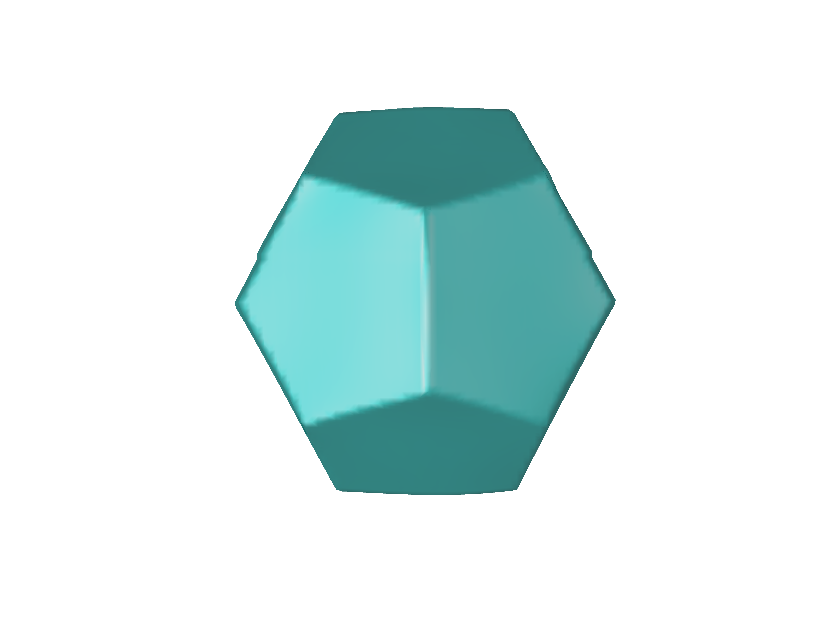}
\includegraphics[width=0.13\textwidth, clip, trim = 7cm 4cm 8cm 2cm]{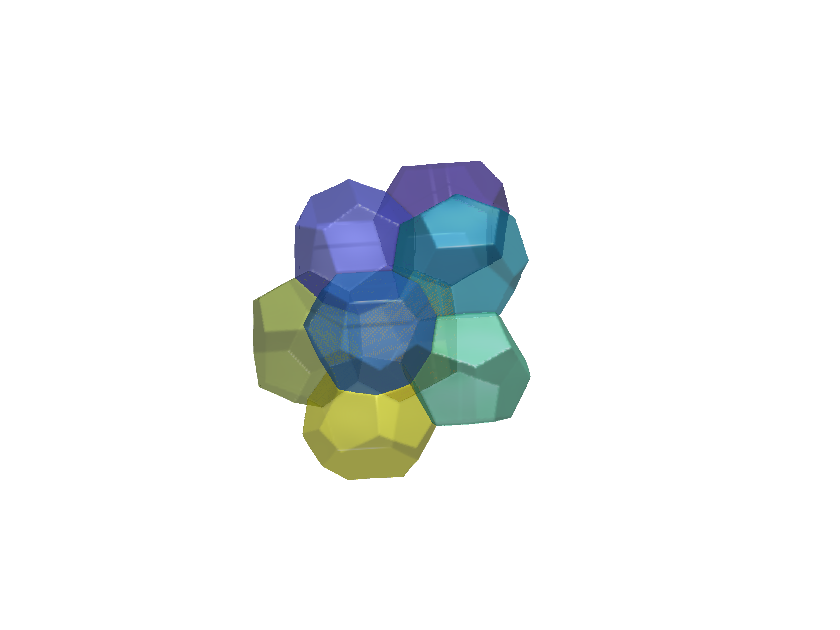}
\includegraphics[width=0.14\textwidth, clip, trim = 5cm 4cm 7cm 2cm]{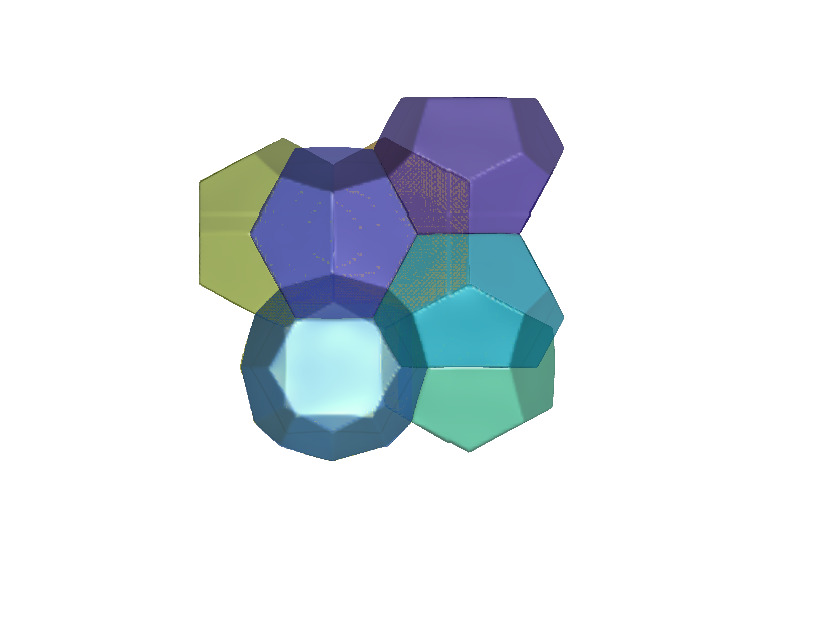}
\caption{Partition results for the Dirichlet partition problem with periodic boundary conditions in the $3$-dimensional flat torus with $n=8$. Left two: different views of the first type Weaire--Phelan structure. Middle two: different views of the second type Weaire--Phelan structure. Right two: different views of the periodic extension of the structure. See Section~\ref{sec:periodic}.}
\label{Fig: D-P P Model 3d n8}
\end{figure}

{\bf Three-dimensional flat torus:} As for the $3$-dimensional flat torus, we consider the situations for $n=4$ and $8$. The results are based on the evaluation of the trained network on a uniform $100^3$ grid. For $n=4$, we obtain a partition of the cube consisting of four identical rhombic dodecahedron structures, as displayed in Figure~\ref{Fig: D-P P Model 3d n4}. For $n=8$, we obtain a partition of the cube that is similar to the Weaire--Phelan structure. Figure~\ref{Fig: D-P P Model 3d n8} shows different views of the first and second type Weaire--Phelan structures, together with different views of a periodic extension of the partition.

\begin{figure}[ht!]
\centering
	\begin{tabular}{|m{1.5cm}|m{2.2cm}|m{2.2cm}|m{2.2cm}|m{2.2cm}|}
	\hline   & \centering{$j=1$} & \centering{$j=2$}  & \centering{$j=3$}  & \quad\quad\;$j=4$ \\ \hline
        $x_j=0$ & \includegraphics[width=0.16\textwidth, clip, trim = 3cm 2cm 3cm 1cm]{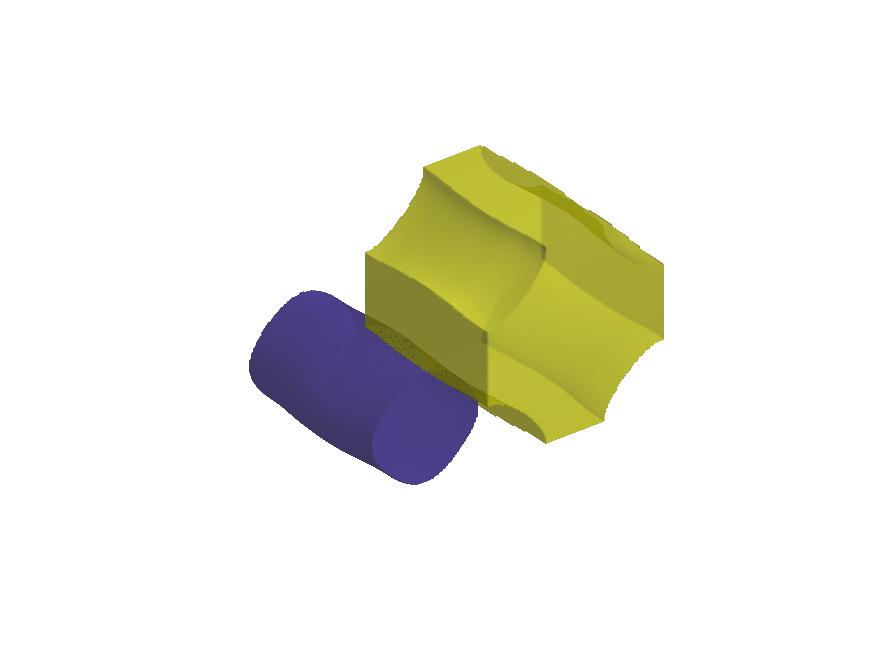}
        & \includegraphics[width=0.16\textwidth, clip, trim = 3cm 2cm 3cm 1cm]{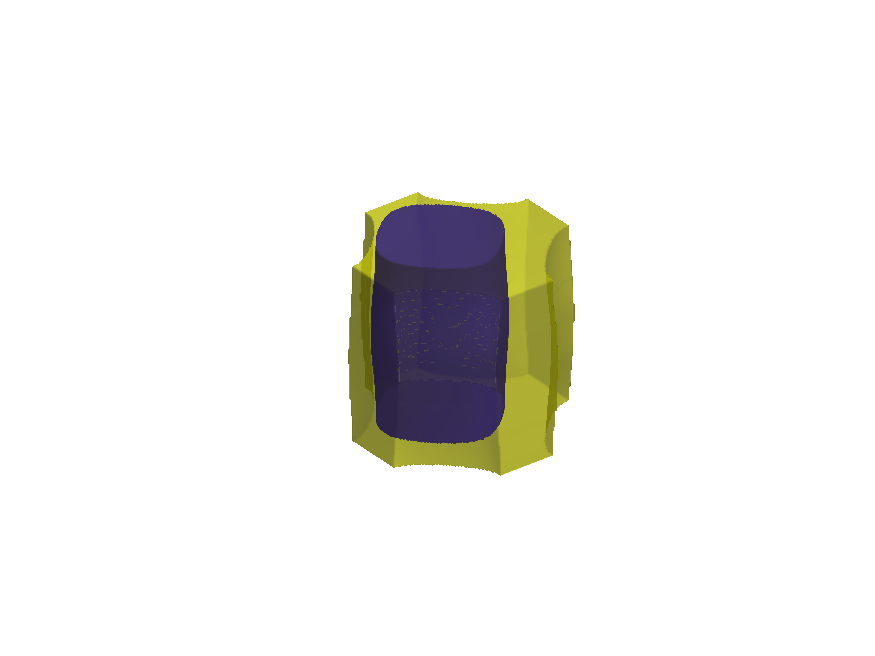}
        & \includegraphics[width=0.16\textwidth, clip, trim = 3cm 2cm 3cm 1cm]{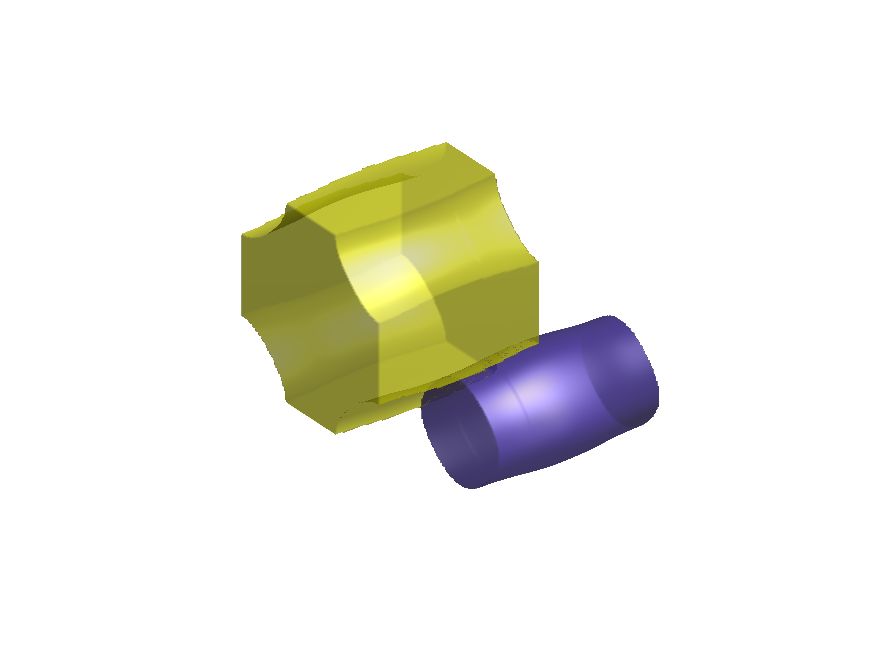}
        & \includegraphics[width=0.16\textwidth, clip, trim = 3cm 2cm 3cm 1cm]{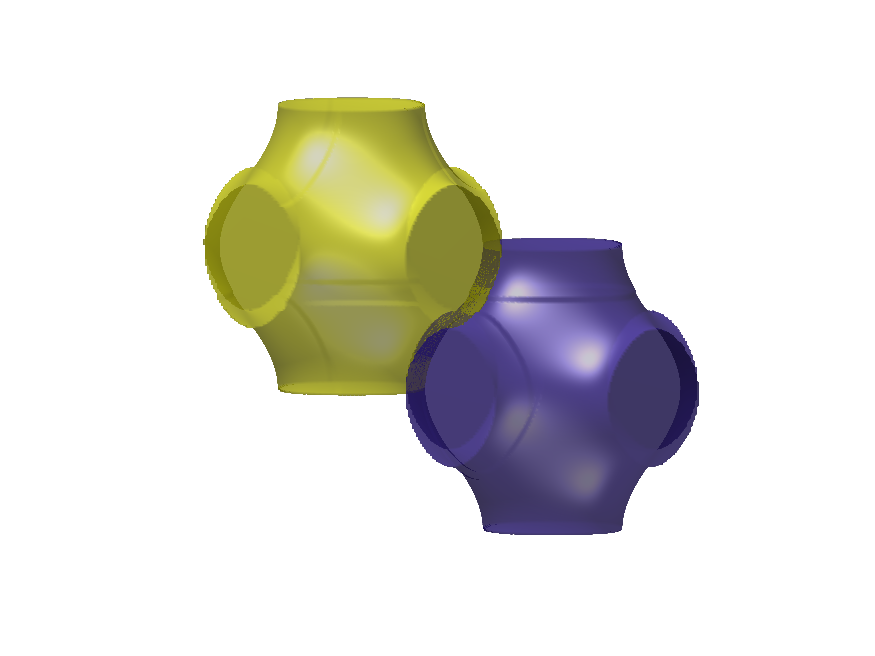} \\ \hline
        $x_j=0.25$ & \includegraphics[width=0.16\textwidth, clip, trim = 3cm 2cm 3cm 1cm]{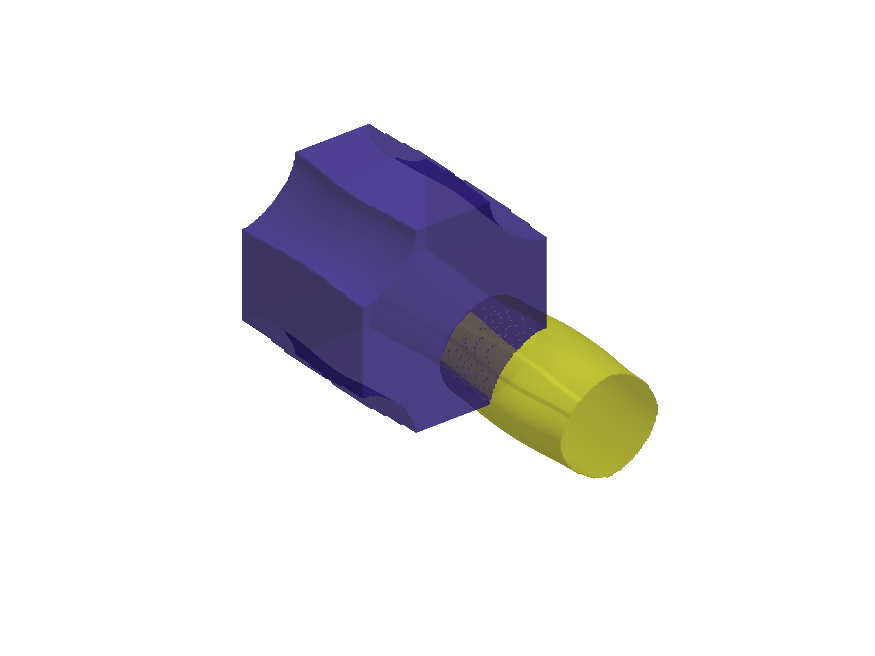}
        & \includegraphics[width=0.16\textwidth, clip, trim = 3cm 2cm 3cm 1cm]{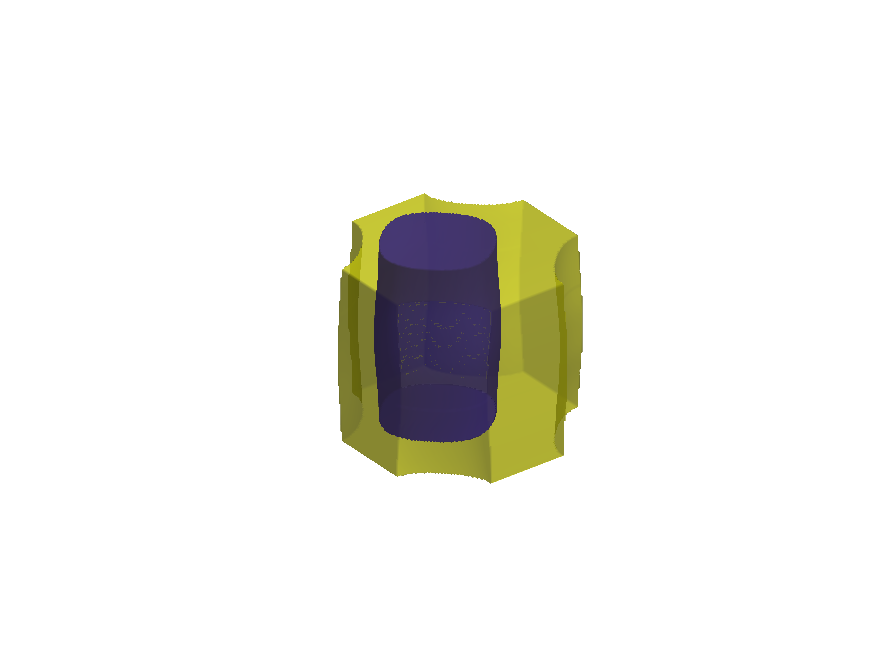}
        &\includegraphics[width=0.16\textwidth, clip, trim = 3cm 2cm 3cm 1cm]{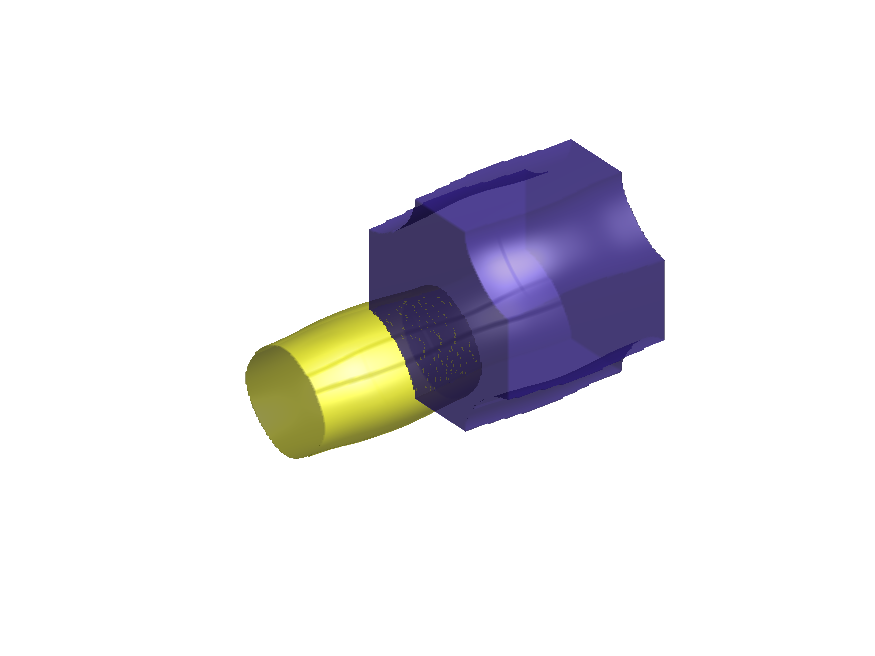}
        & \includegraphics[width=0.16\textwidth, clip, trim = 3cm 2cm 3cm 1cm]{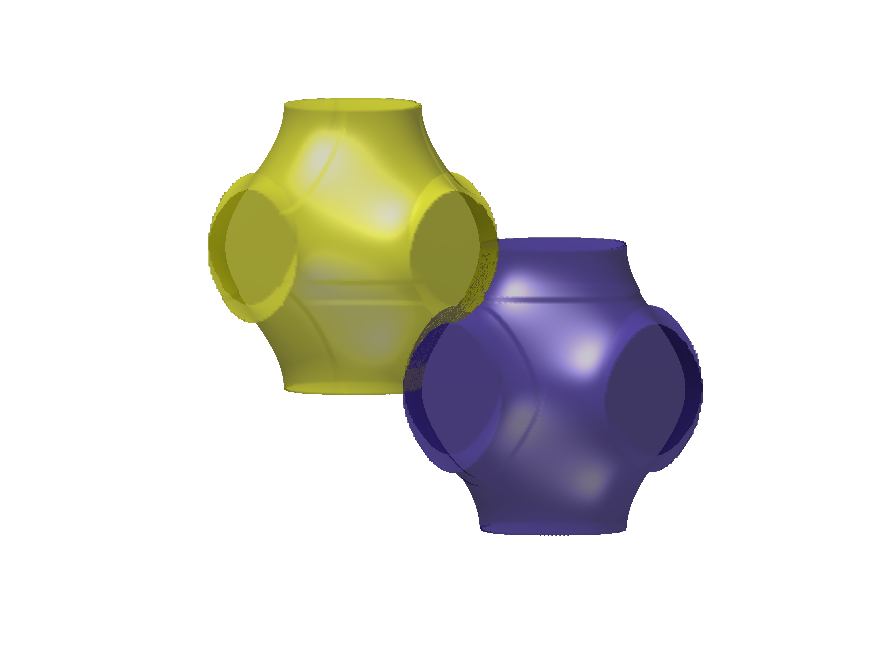}\\ \hline
        $x_j=0.5$ & \includegraphics[width=0.16\textwidth, clip, trim = 3cm 2cm 3cm 1cm]{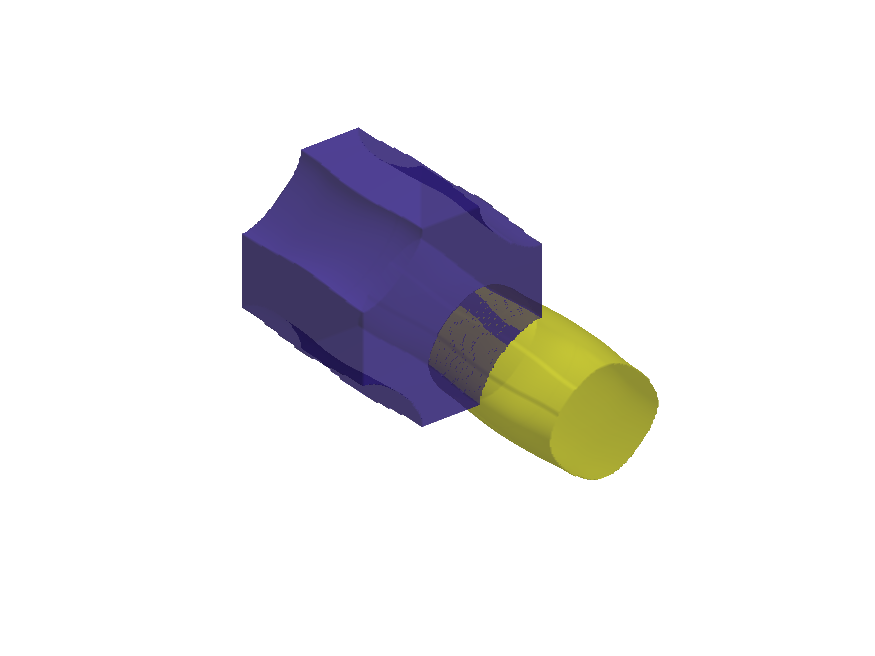}
        & \includegraphics[width=0.16\textwidth, clip, trim = 3cm 2cm 3cm 1cm]{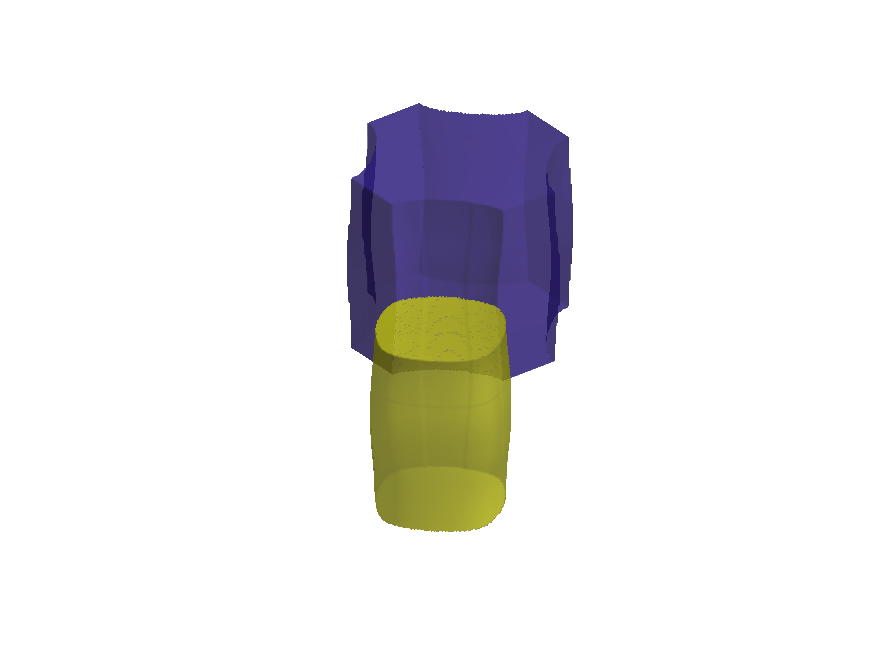}
        & \includegraphics[width=0.16\textwidth, clip, trim = 3cm 2cm 3cm 1cm]{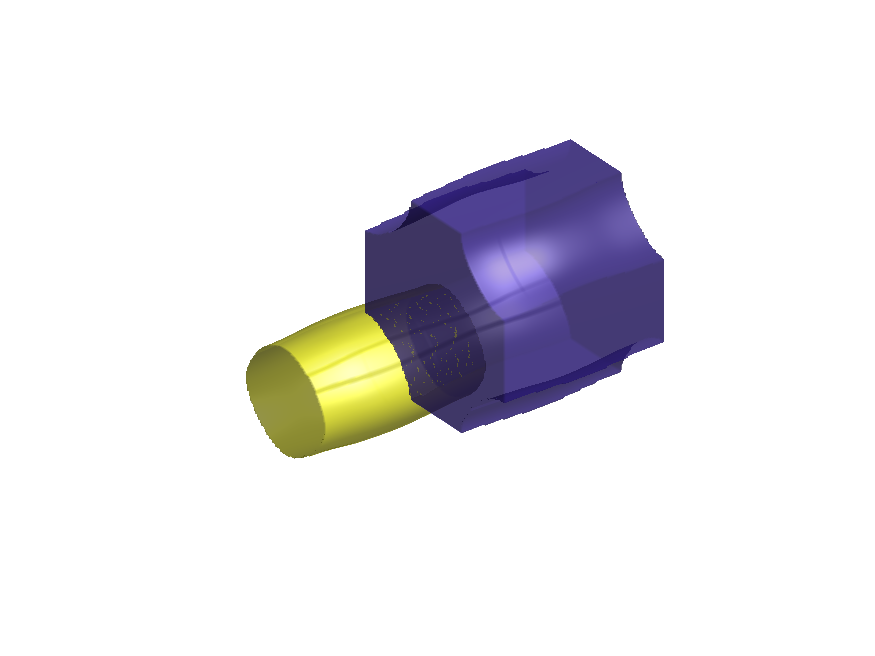}
        & \includegraphics[width=0.16\textwidth, clip, trim = 3cm 2cm 3cm 1cm]{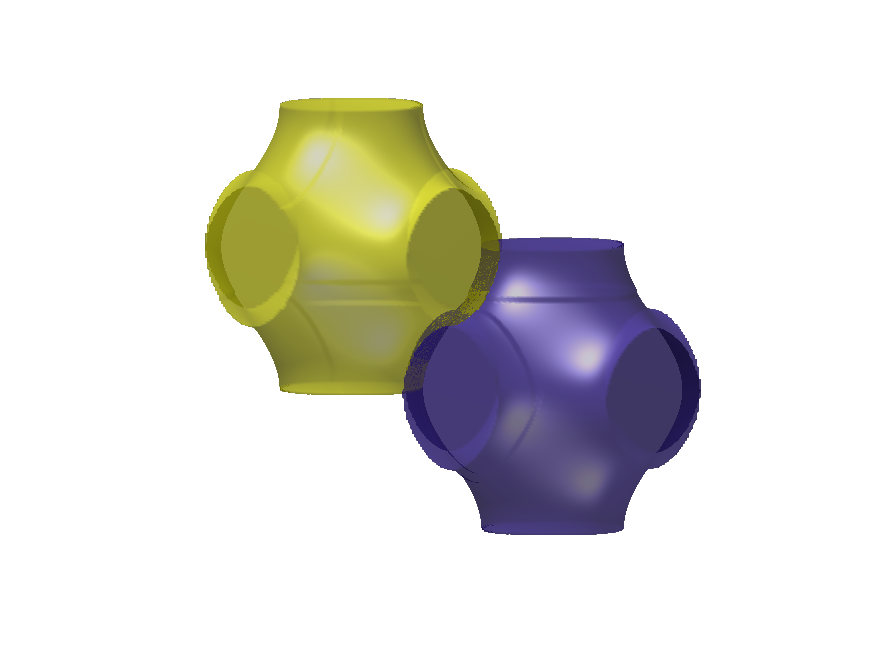}\\ \hline
        $x_j=0.75$ & \includegraphics[width=0.16\textwidth, clip, trim = 3cm 2cm 3cm 1cm]{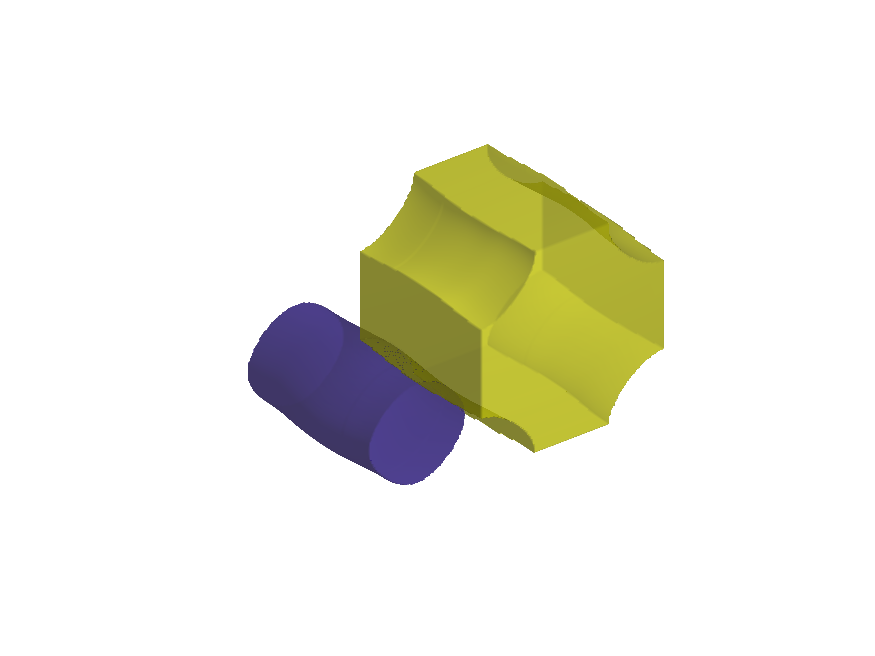}
        & \includegraphics[width=0.16\textwidth, clip, trim = 3cm 2cm 3cm 1cm]{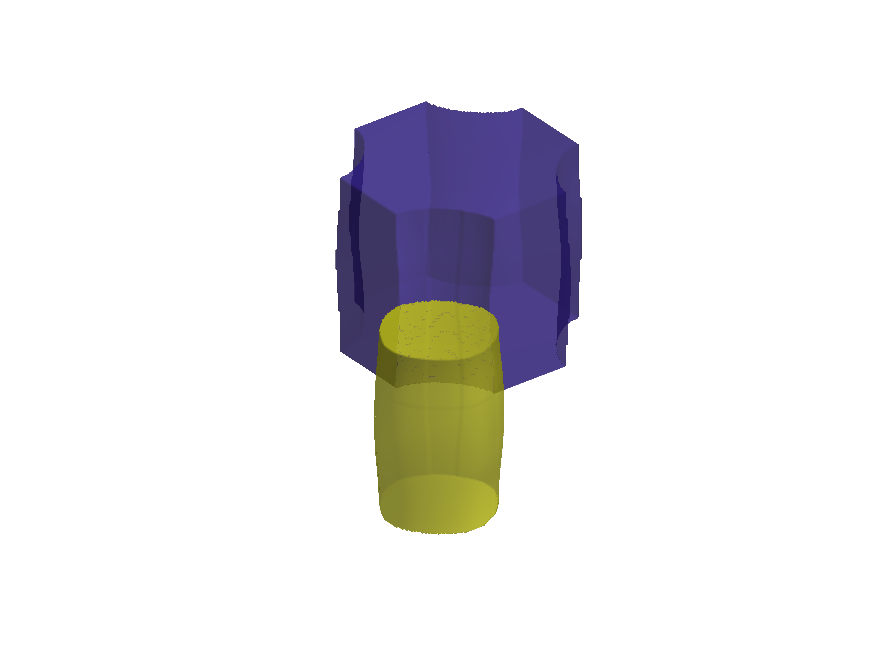}
        & \includegraphics[width=0.16\textwidth, clip, trim = 3cm 2cm 3cm 1cm]{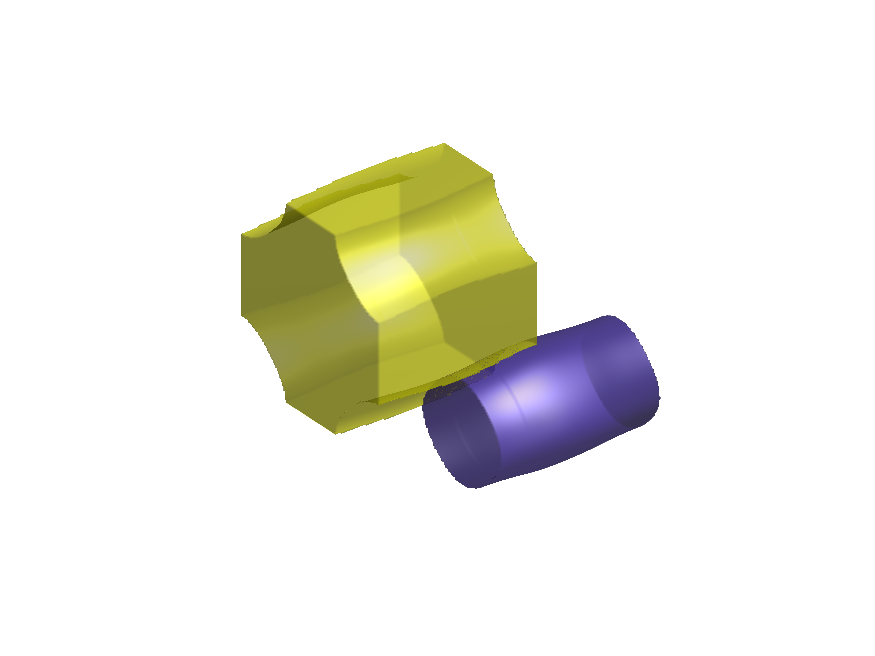}
        & \includegraphics[width=0.16\textwidth, clip, trim = 3cm 2cm 3cm 1cm]{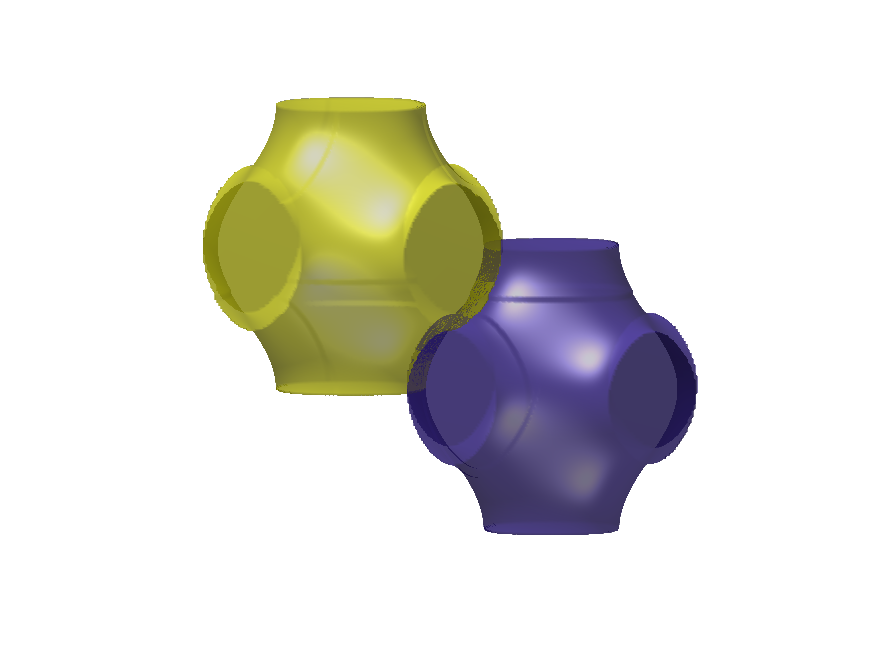}
        \\
        \hline
    \end{tabular}
\caption{A $n=2$ Dirichlet partition problem with periodic boundary conditions in the $4$-dimensional flat torus. The four columns correspond to the slides perpendicular to the $x_1$-, $x_2$-, $x_3$-, and $x_4$-axis respectively. The four rows correspond to the slices at $x_j = 0, 0.25, 0.5, 0.75$, respectively. These partitions are based on testing the trained neural network on a uniform $100^3$ grid. See Section~\ref{sec:periodic}.}
\label{Fig: Dirichlet Partition 4d unit tesseract n=2}
\end{figure}

\begin{figure}[ht!]
\centering
	\begin{tabular}{|m{1.5cm}|m{2.2cm}|m{2.2cm}|m{2.2cm}|m{2.2cm}|}
	\hline   & \centering{$j=1$} & \centering{$j=2$}  & \centering{$j=3$}  & \quad\quad$j=4$ \\ \hline
         $x_j=0$ & \includegraphics[width=0.16\textwidth, clip, trim = 4cm 1.5cm 3cm 1cm]{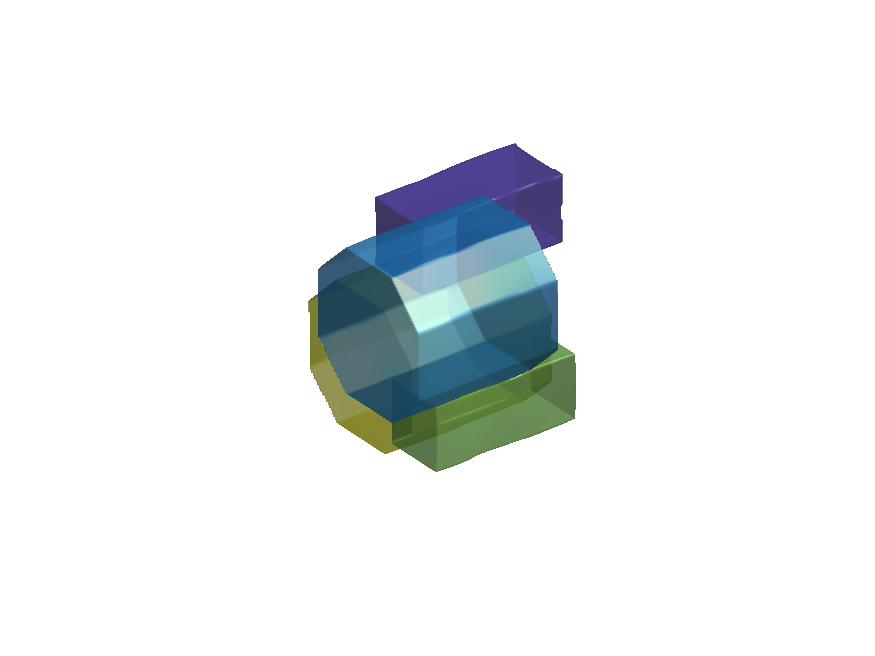} 
         & \includegraphics[width=0.16\textwidth, clip, trim = 4cm 1.5cm 3cm 1cm]{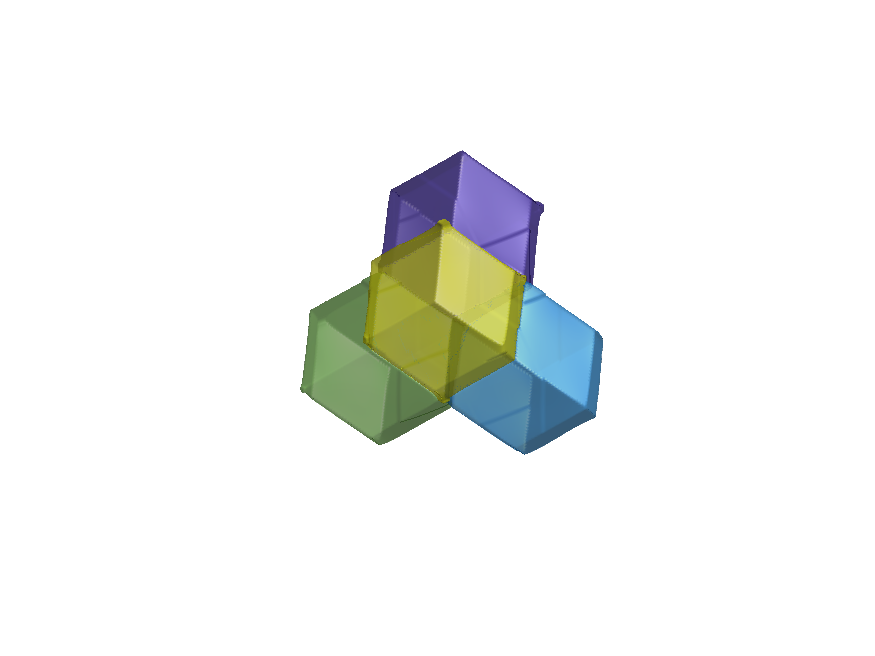} 
         & \includegraphics[width=0.16\textwidth, clip, trim = 4cm 1.5cm 3cm 1cm]{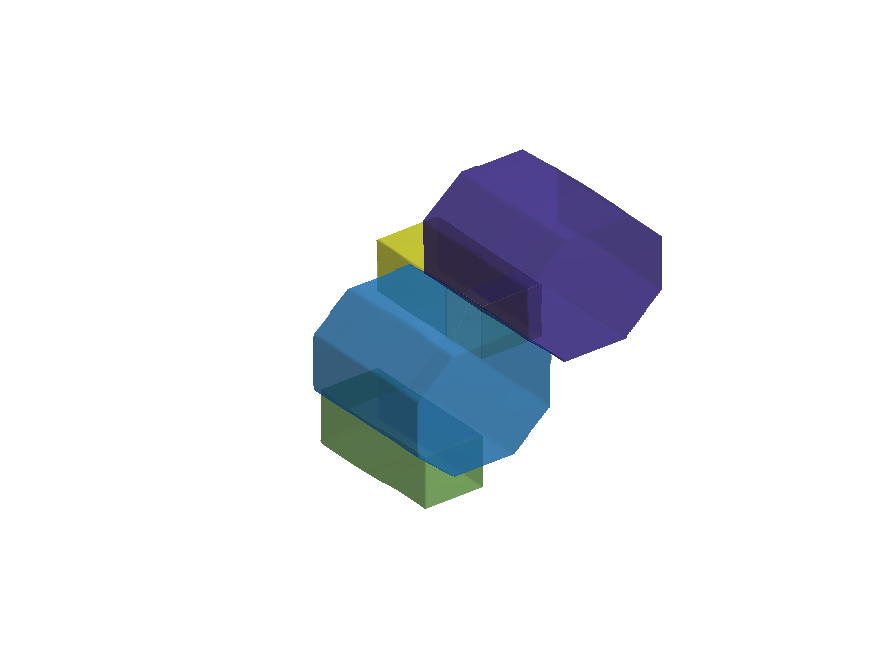}
         & \includegraphics[width=0.16\textwidth, clip, trim = 4cm 1.5cm 3cm 1cm]{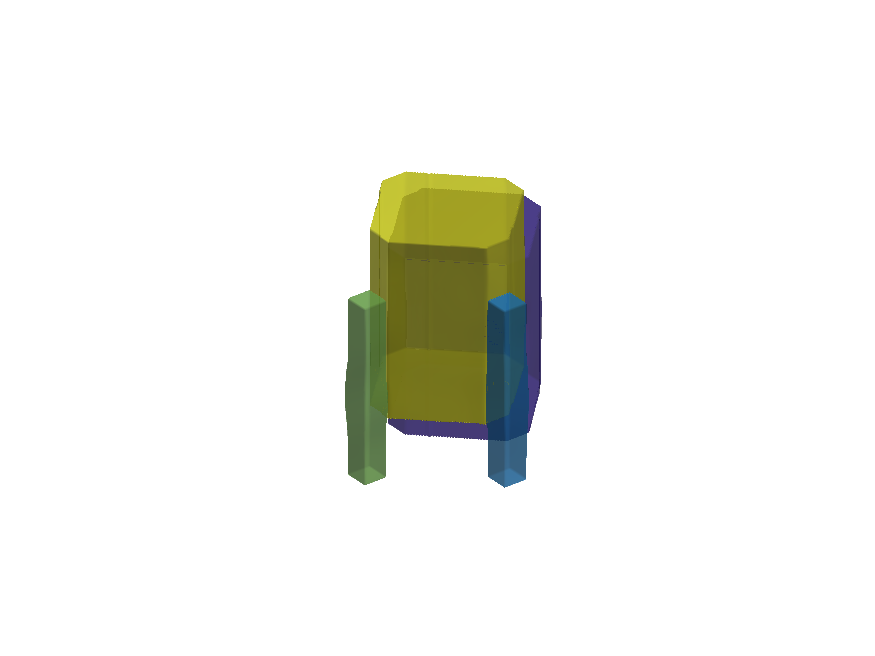}  \\ \hline
         $x_j=0.25$ & \includegraphics[width=0.16\textwidth, clip, trim = 4cm 1.5cm 3cm 1cm]{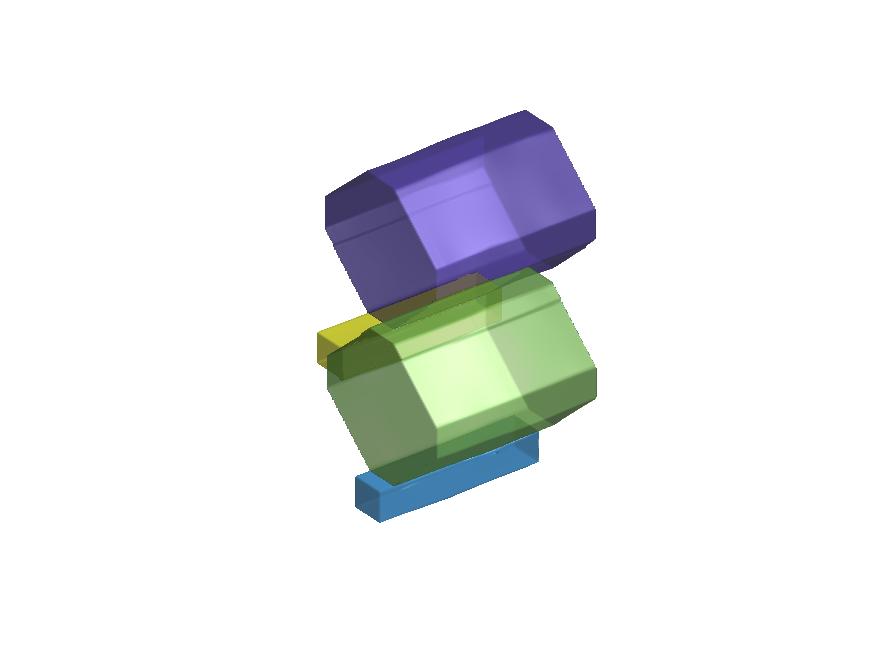}
         & \includegraphics[width=0.16\textwidth, clip, trim = 4cm 1.5cm 3cm 1cm]{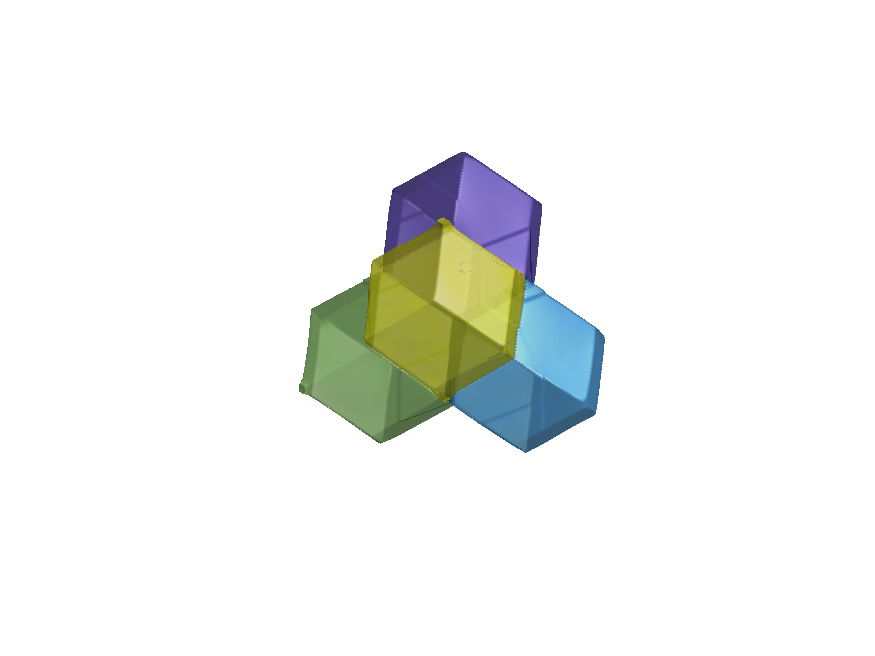}
         & \includegraphics[width=0.16\textwidth, clip, trim = 4cm 1.5cm 3cm 1cm]{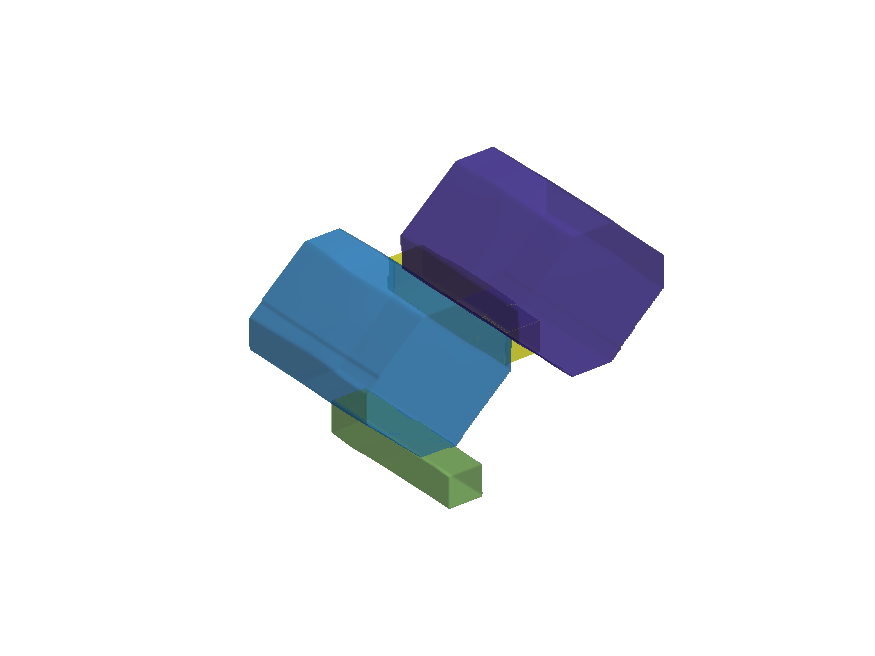}
         & \includegraphics[width=0.16\textwidth, clip, trim = 4cm 1.5cm 3cm 1cm]{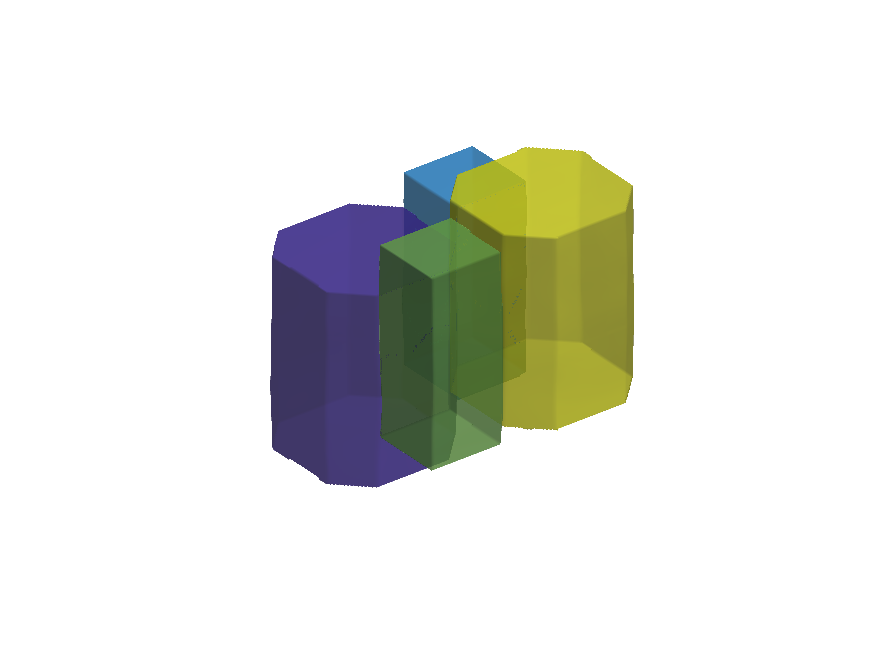}  \\ \hline
         $x_j=0.5$ & \includegraphics[width=0.16\textwidth, clip, trim = 4cm 1.5cm 3cm 1cm]{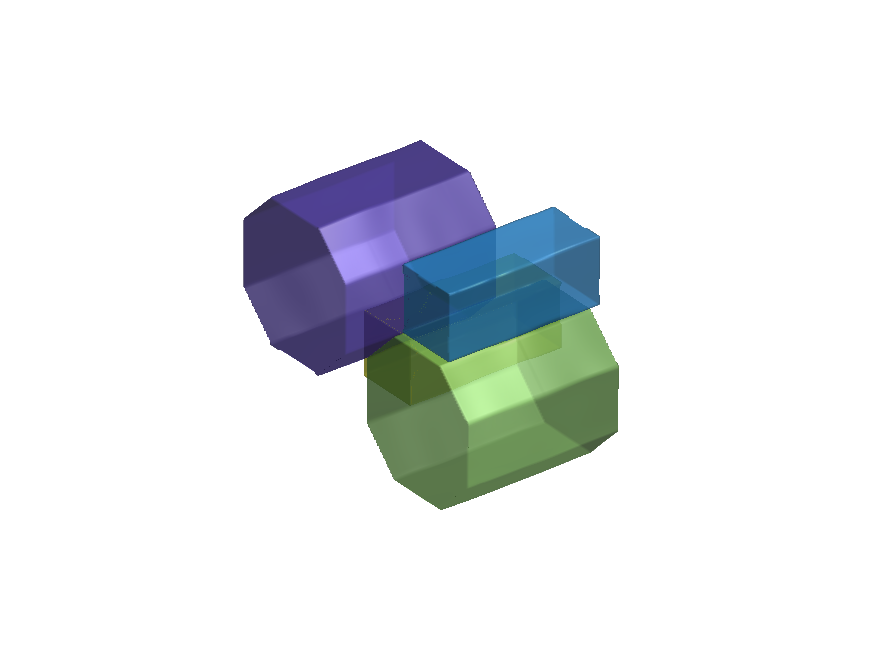}
         & \includegraphics[width=0.16\textwidth, clip, trim = 4cm 1.5cm 3cm 1cm]{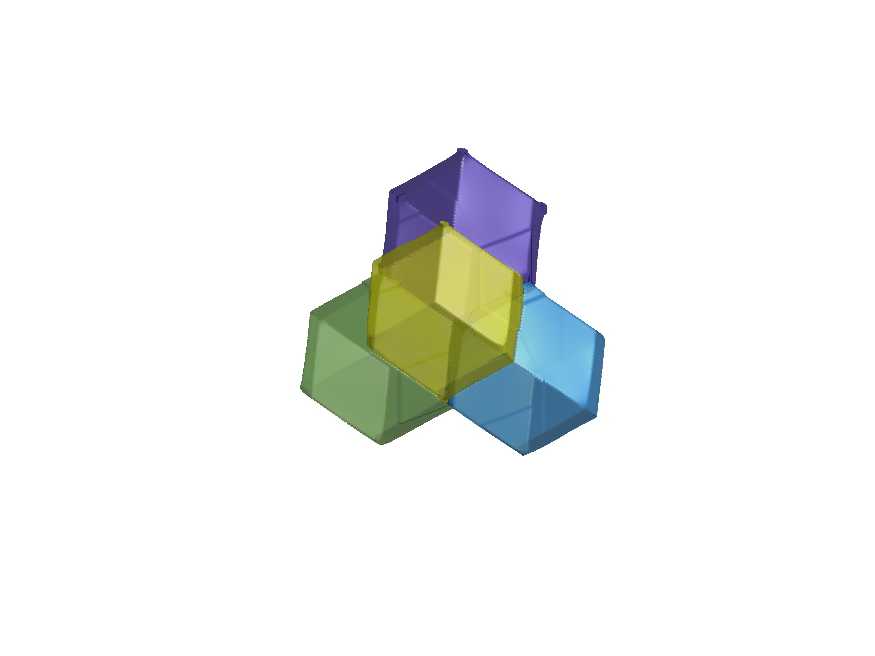}
         & \includegraphics[width=0.16\textwidth, clip, trim = 4cm 1.5cm 3cm 1cm]{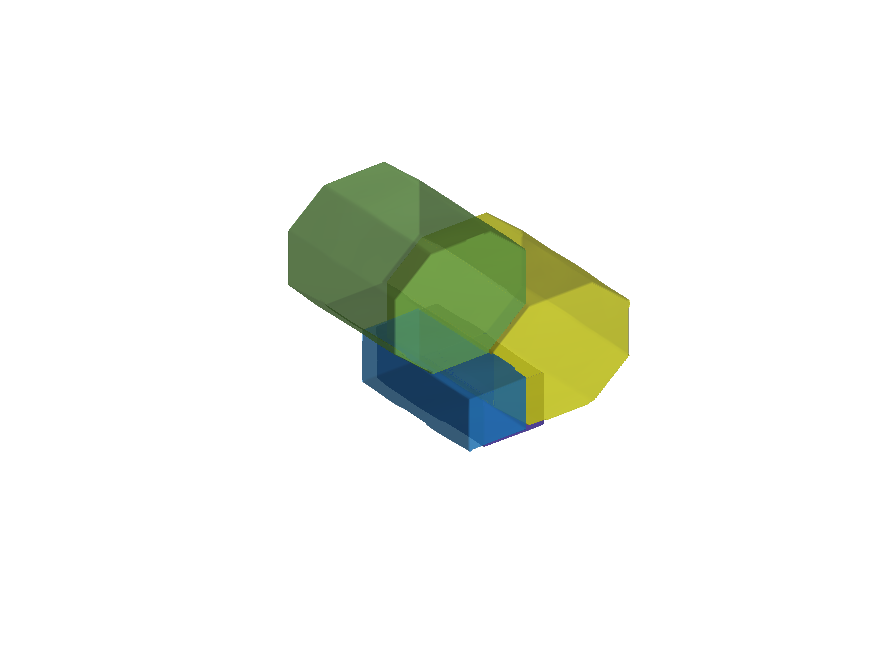}
         & \includegraphics[width=0.16\textwidth, clip, trim = 4cm 1.5cm 3cm 1cm]{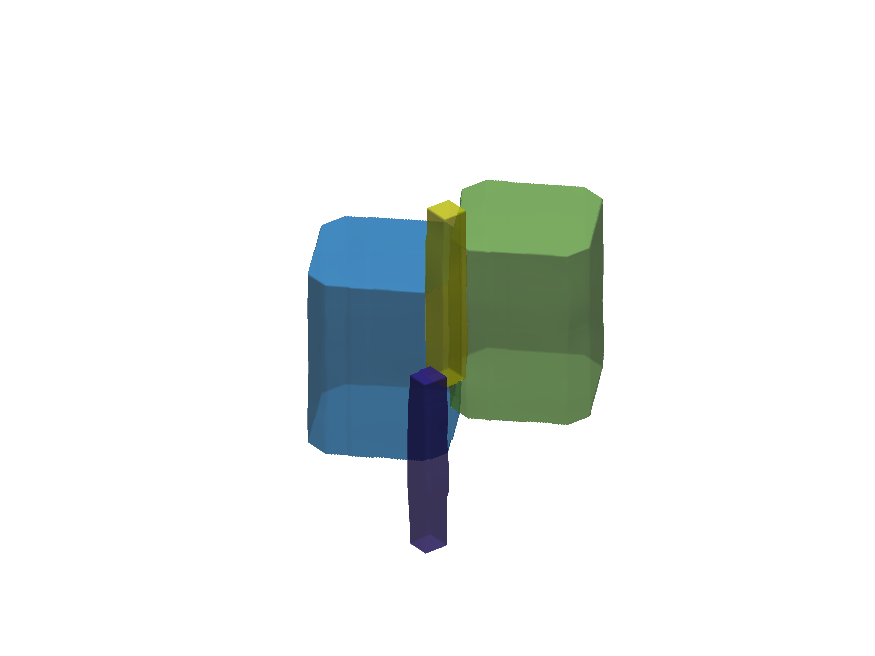}  \\ \hline
         $x_j=0.75$ & \includegraphics[width=0.16\textwidth, clip, trim = 4cm 1.5cm 3cm 1cm]{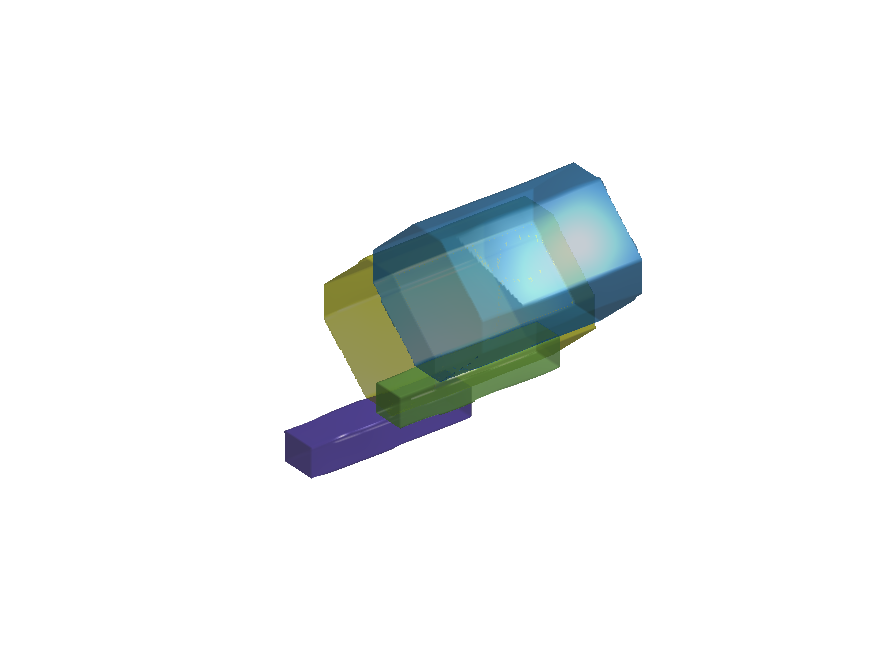}
         & \includegraphics[width=0.16\textwidth, clip, trim = 4cm 1.5cm 3cm 1cm]{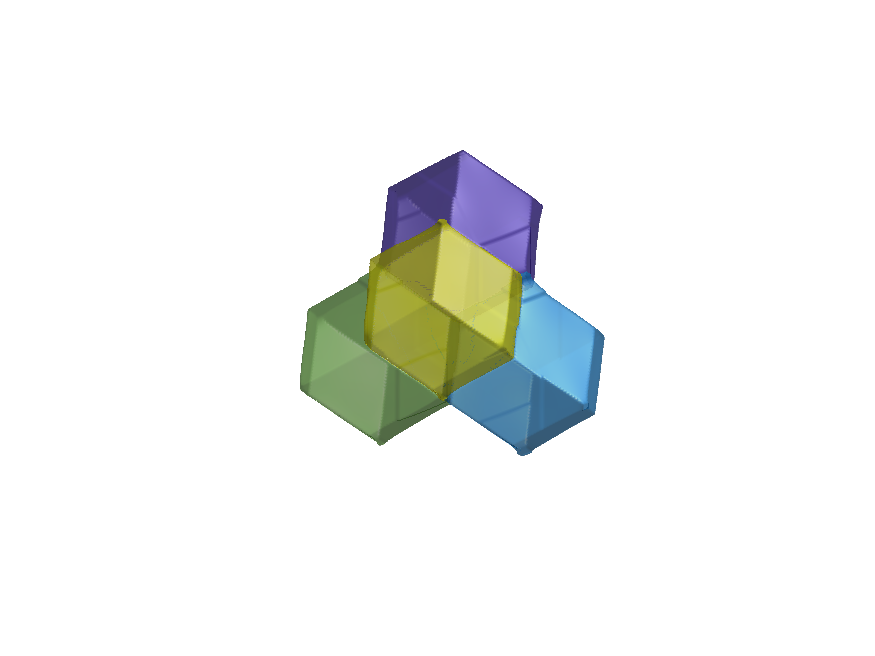}
         & \includegraphics[width=0.16\textwidth, clip, trim = 4cm 1.5cm 3cm 1cm]{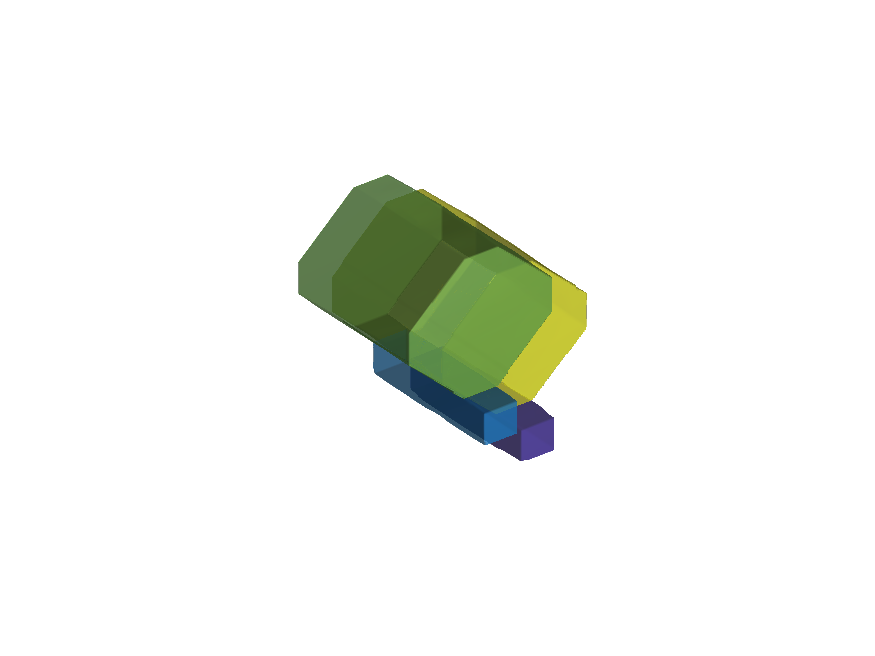}
         & \includegraphics[width=0.16\textwidth, clip, trim = 4cm 1.5cm 3cm 1cm]{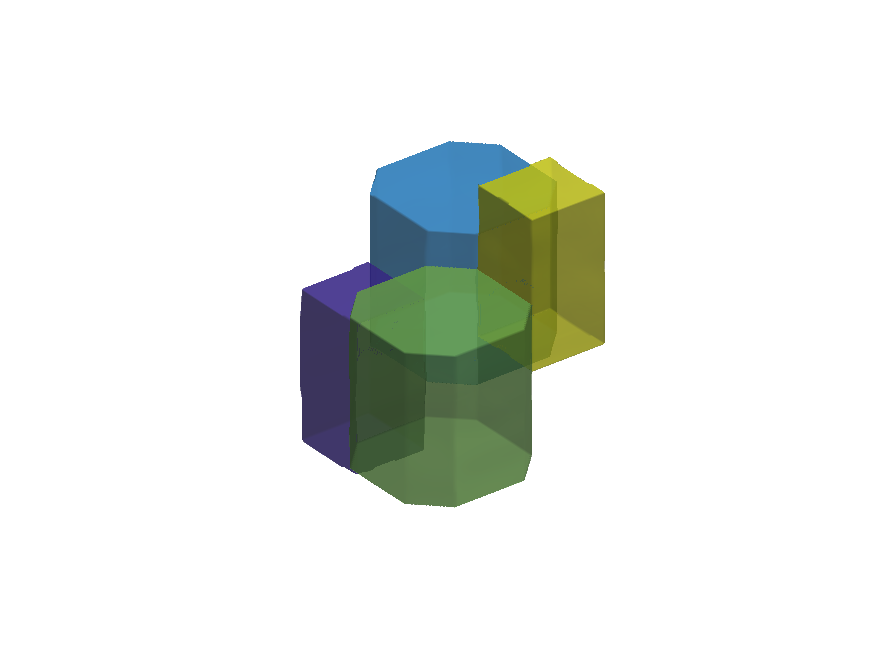}
    \\
     		\hline
    	\end{tabular}
\caption{A $n=4$ Dirichlet partition problem with periodic boundary conditions in the $4$-dimensional flat torus. The four columns correspond to the slides perpendicular to the $x_1$-, $x_2$-, $x_3$-, and $x_4$-axis respectively. The four rows correspond to the slices at $x_j = 0, 0.25, 0.5, 0.75$, respectively. These partitions are based on testing the trained neural network on a uniform $100^3$ grid. See Section~\ref{sec:periodic}.}
\label{Fig: Dirichlet Partition 4d unit tesseract n=4}
\end{figure}

\begin{figure}[ht!]
\centering
	\begin{tabular}{|m{1.5cm}|m{2.2cm}|m{2.2cm}|m{2.2cm}|m{2.2cm}|}
	\hline   & \centering{$j=1$} & \centering{$j=2$}  & \centering{$j=3$}  & \quad\quad\;\;$j=4$ \\ \hline
        $x_j=0$ & \includegraphics[width=0.16\textwidth, clip, trim = 4cm 2cm 4cm 1.5cm]{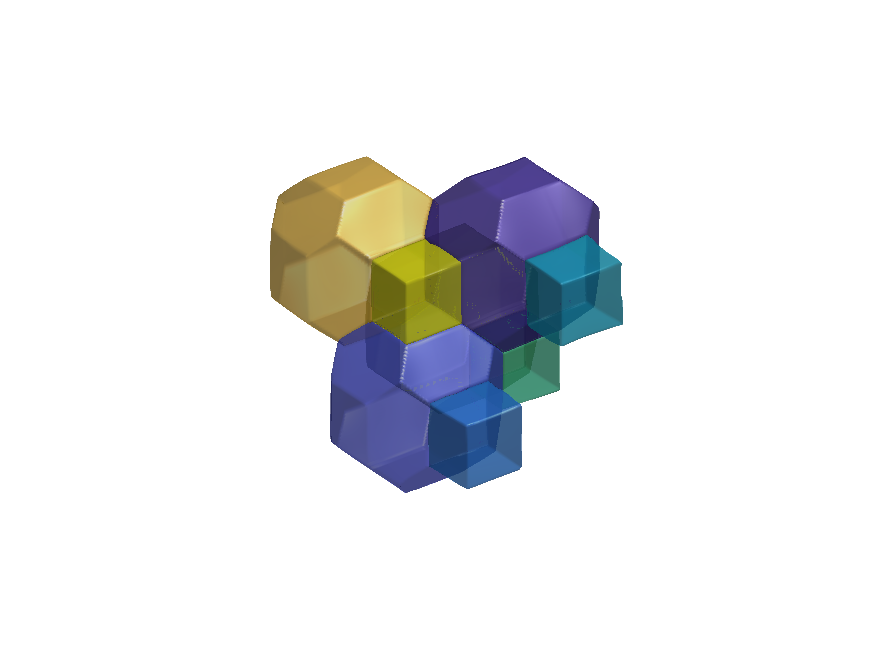}
        & \includegraphics[width=0.16\textwidth, clip, trim = 4cm 2cm 4cm 1.5cm]{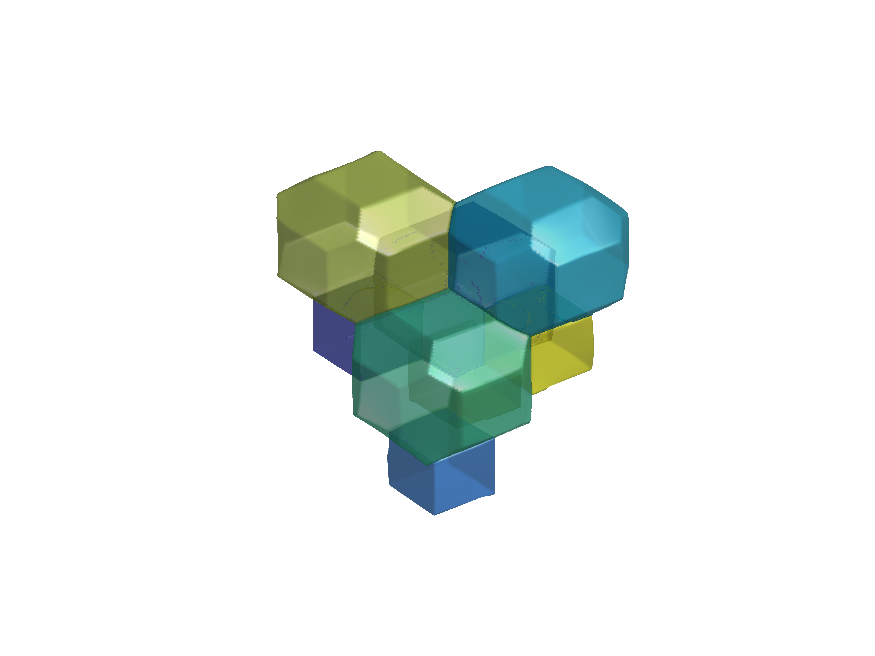}
        & \includegraphics[width=0.16\textwidth, clip, trim = 4cm 2cm 4cm 1.5cm]{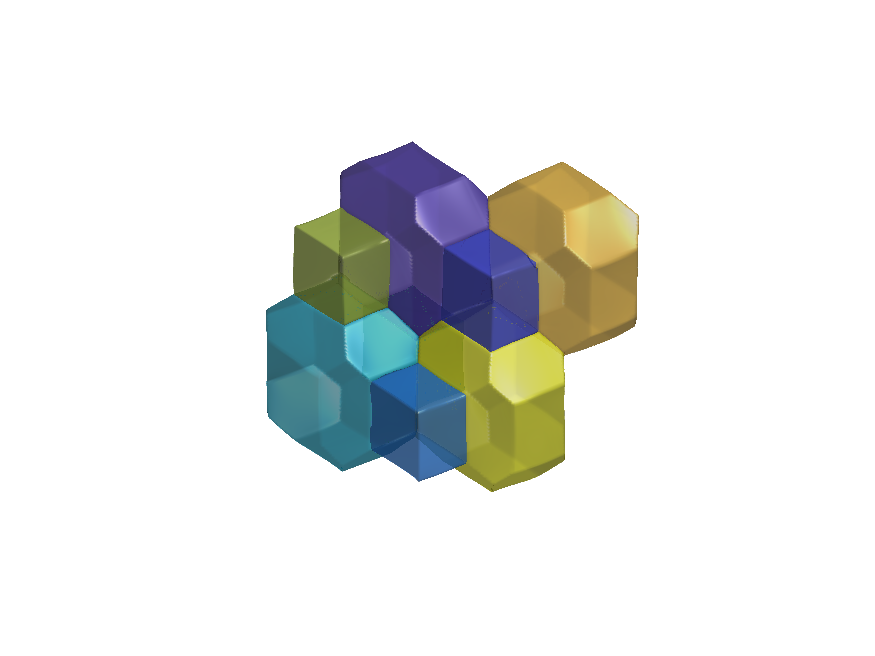}
        & \includegraphics[width=0.16\textwidth, clip, trim = 4cm 2cm 4cm 1.5cm]{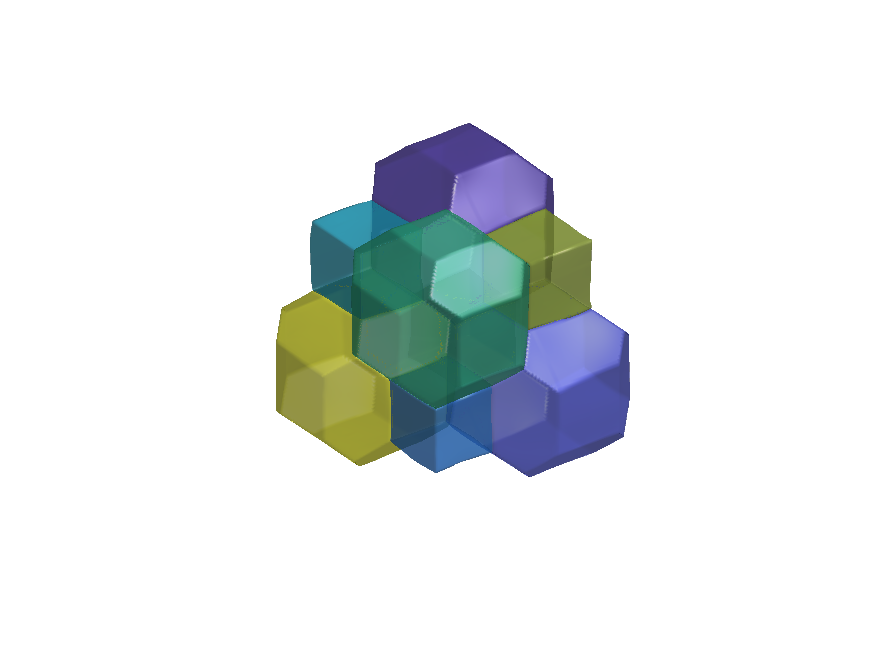}\\ \hline
        $x_j=0.25$ & \includegraphics[width=0.16\textwidth, clip, trim = 4cm 2cm 4cm 1.5cm]{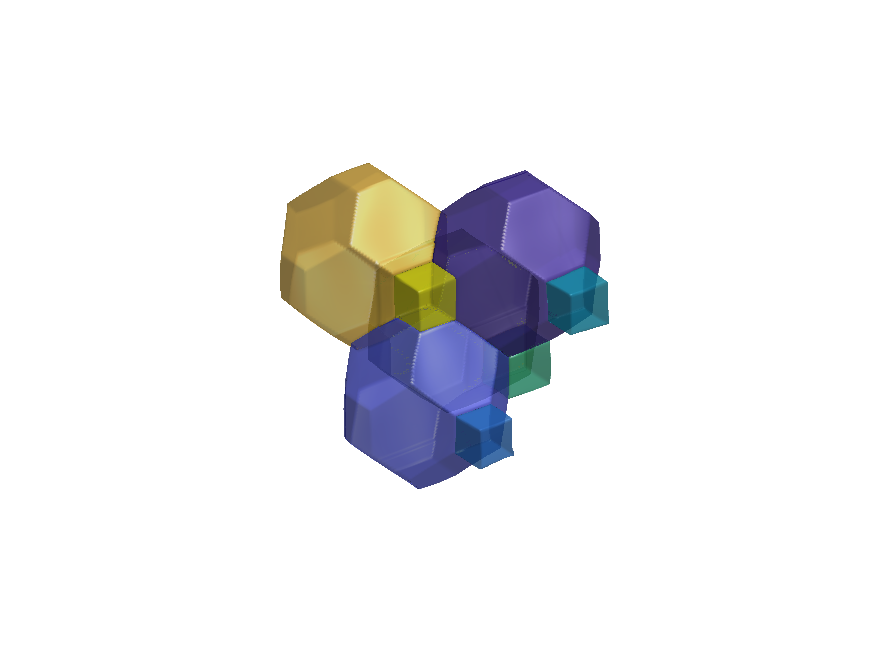}
        & \includegraphics[width=0.16\textwidth, clip, trim = 4cm 2cm 4cm 1.5cm]{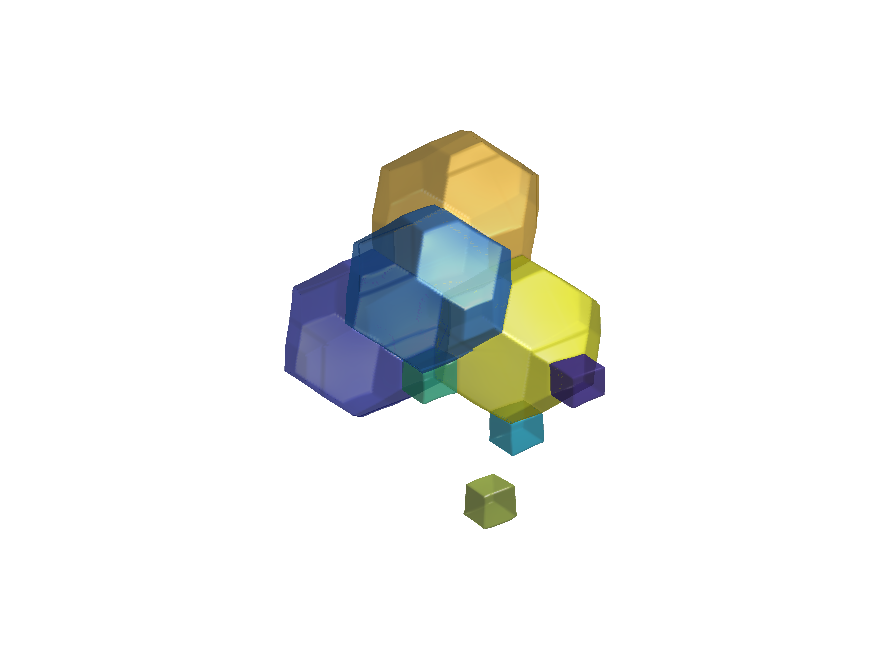}
        & \includegraphics[width=0.16\textwidth, clip, trim = 4cm 2cm 4cm 1.5cm]{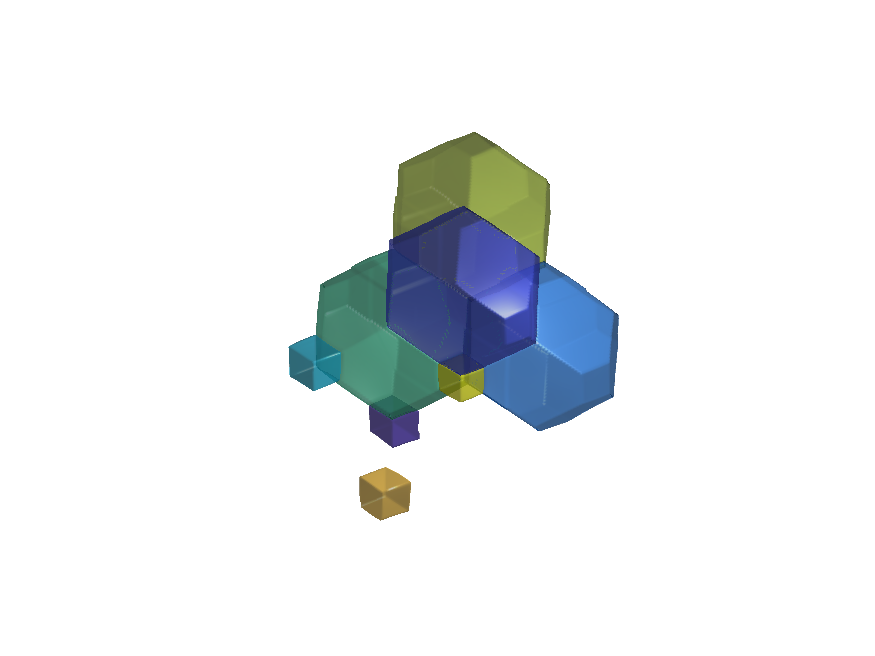}
        & \includegraphics[width=0.16\textwidth, clip, trim = 4cm 2cm 4cm 1.5cm]{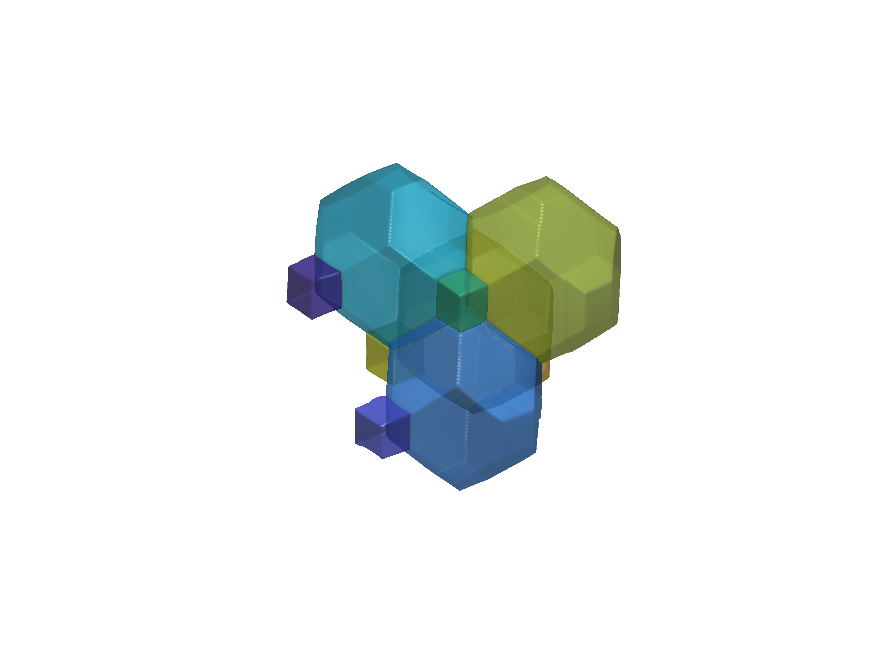}\\ \hline
        $x_j=0.5$ & \includegraphics[width=0.16\textwidth, clip, trim = 4cm 2cm 4cm 1.5cm]{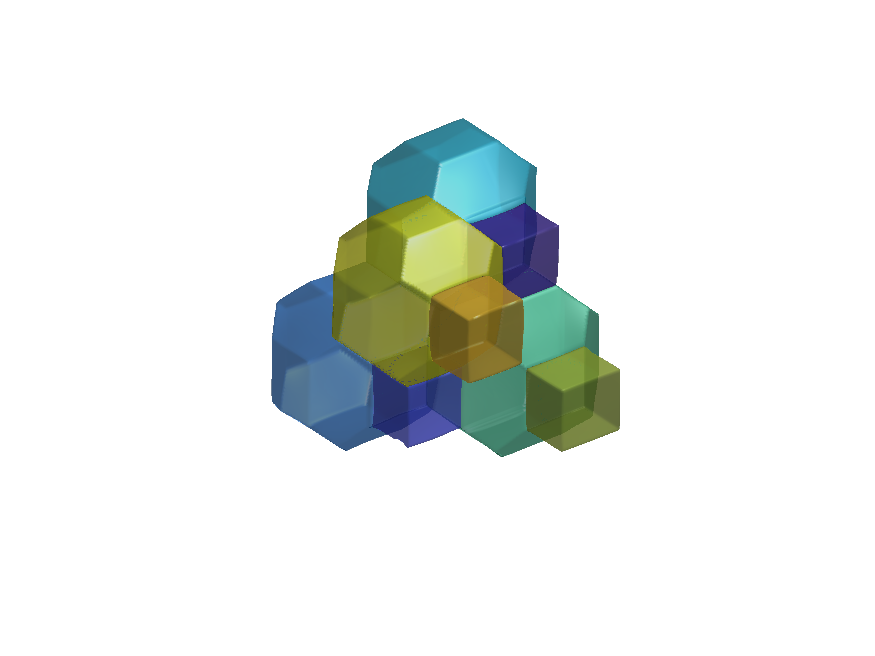}
        & \includegraphics[width=0.16\textwidth, clip, trim = 4cm 2cm 4cm 1.5cm]{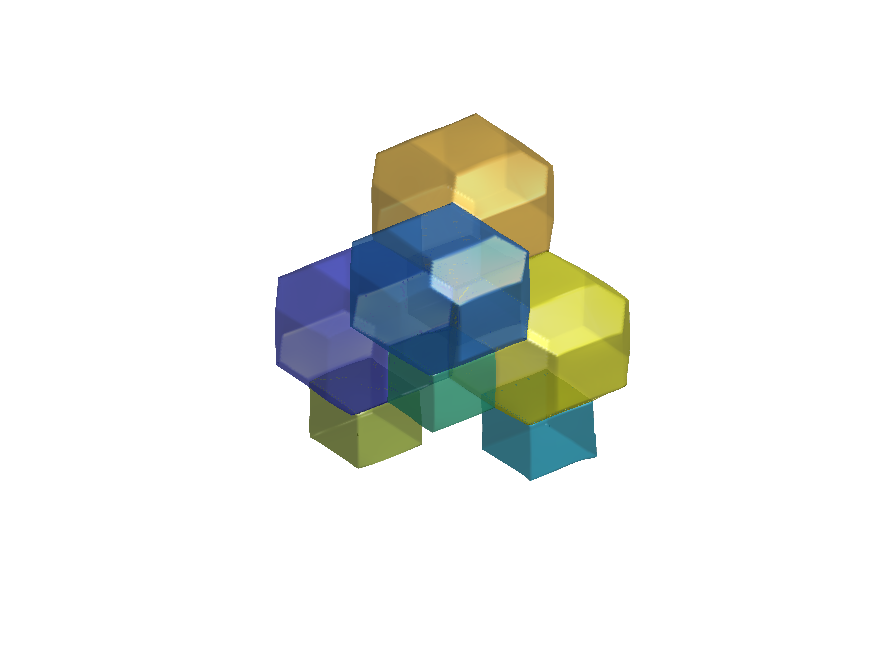}
        & \includegraphics[width=0.16\textwidth, clip, trim = 4cm 2cm 4cm 1.5cm]{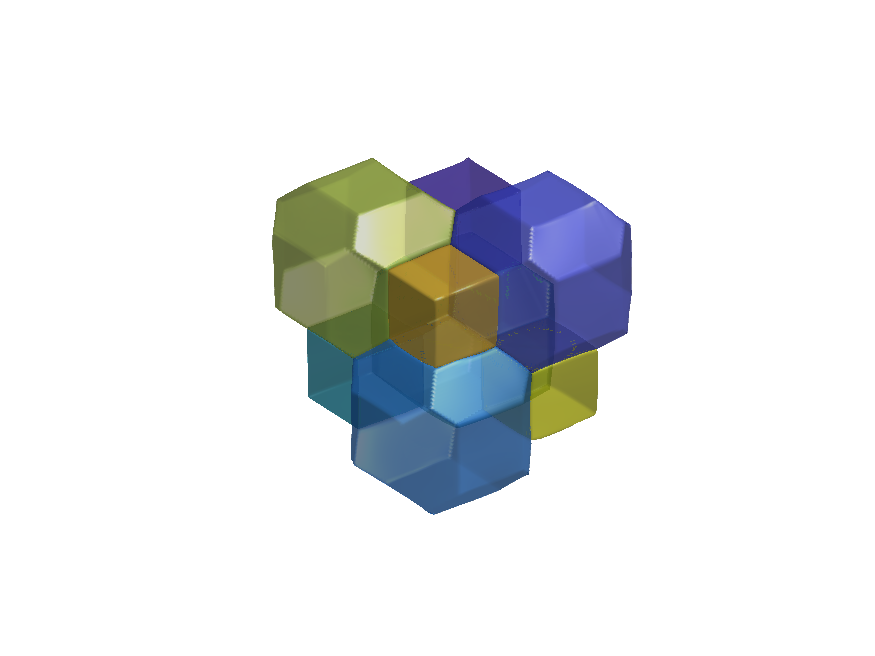}
        & \includegraphics[width=0.16\textwidth, clip, trim = 4cm 2cm 4cm 1.5cm]{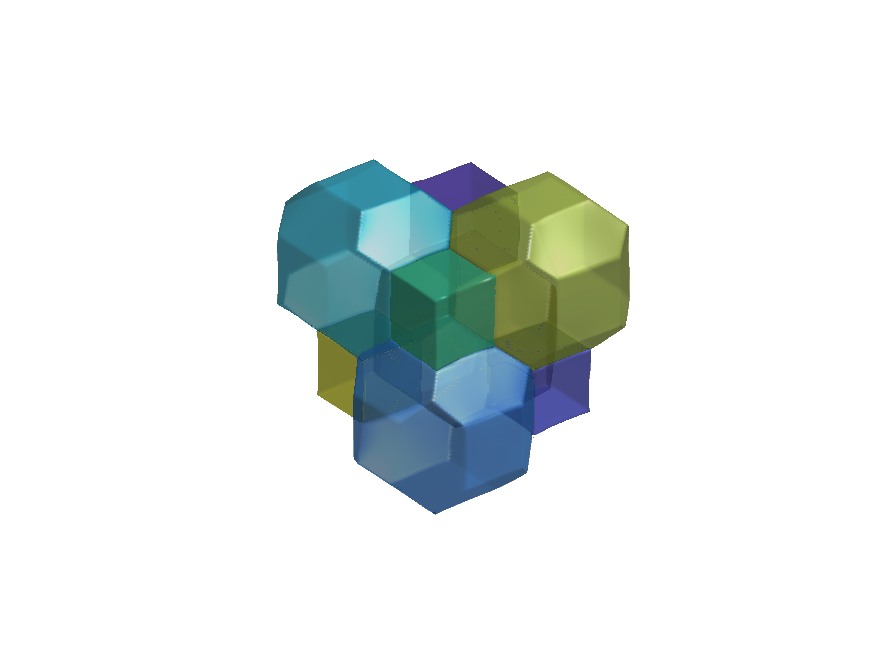}\\ \hline
        $x_j=0.75$ & \includegraphics[width=0.16\textwidth, clip, trim = 4cm 2cm 4cm 1.5cm]{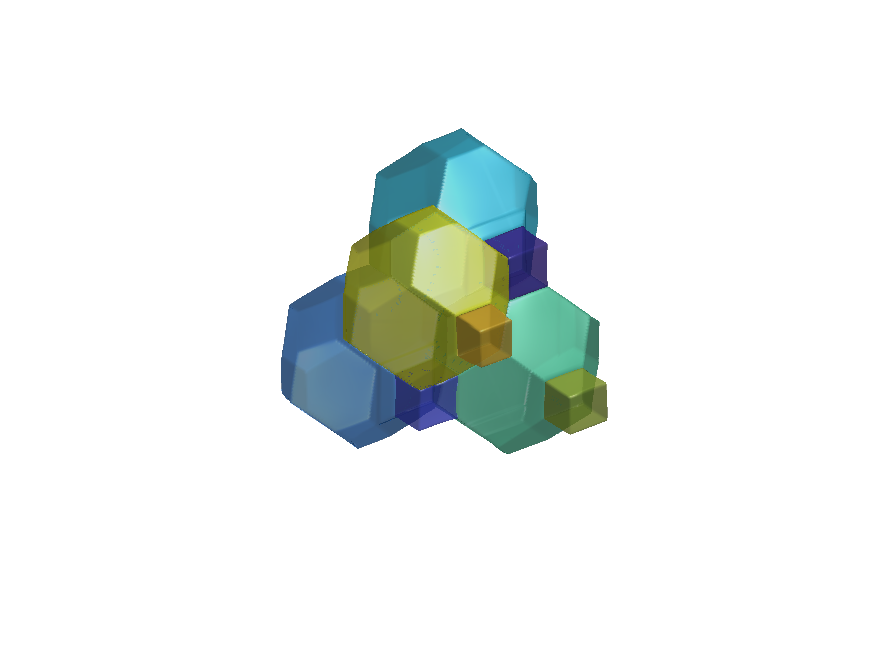}
        & \includegraphics[width=0.16\textwidth, clip, trim = 4cm 2cm 4cm 1.5cm]{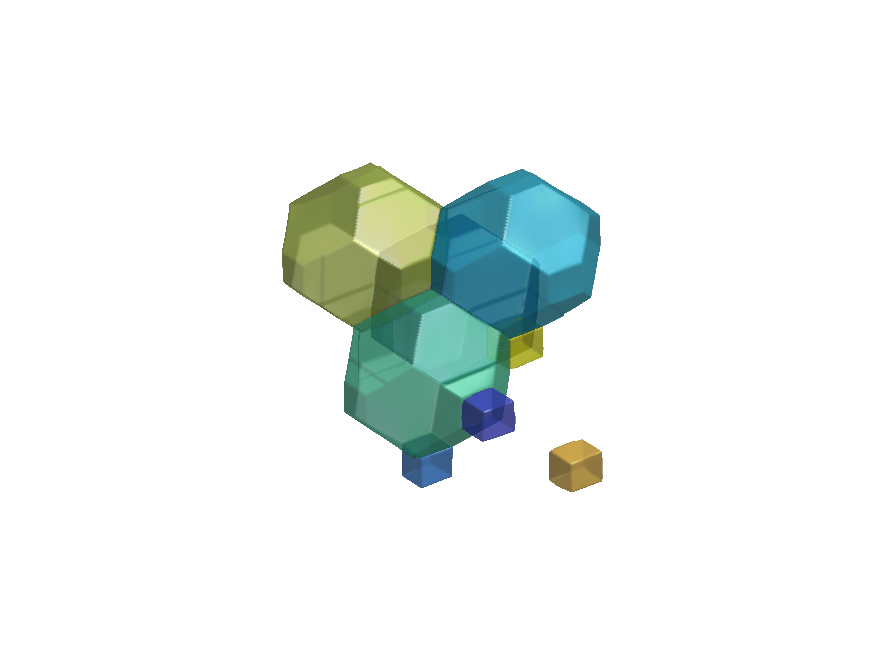}
        & \includegraphics[width=0.16\textwidth, clip, trim = 4cm 2cm 4cm 1.5cm]{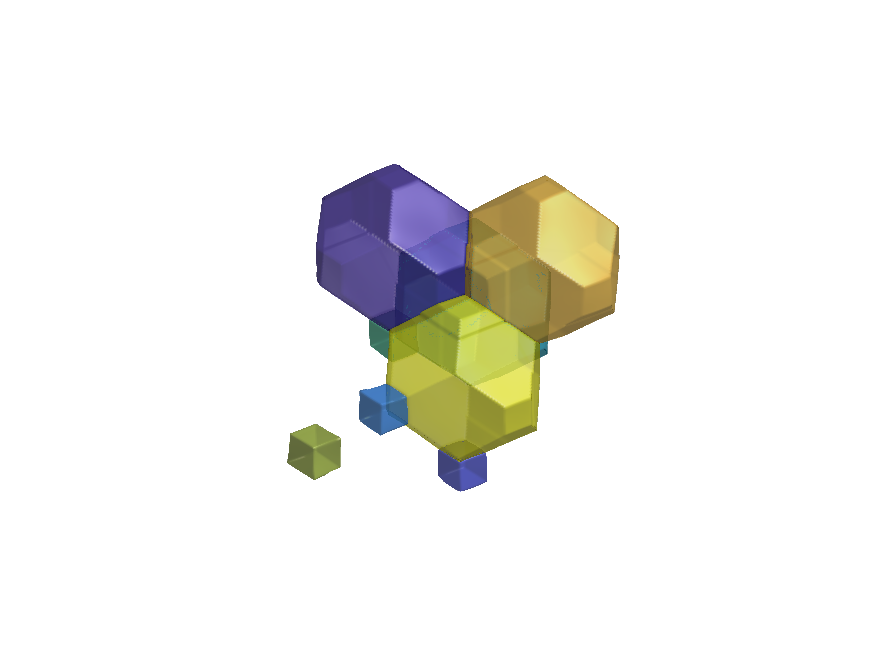}
        & \includegraphics[width=0.16\textwidth, clip, trim = 4cm 2cm 4cm 1.5cm]{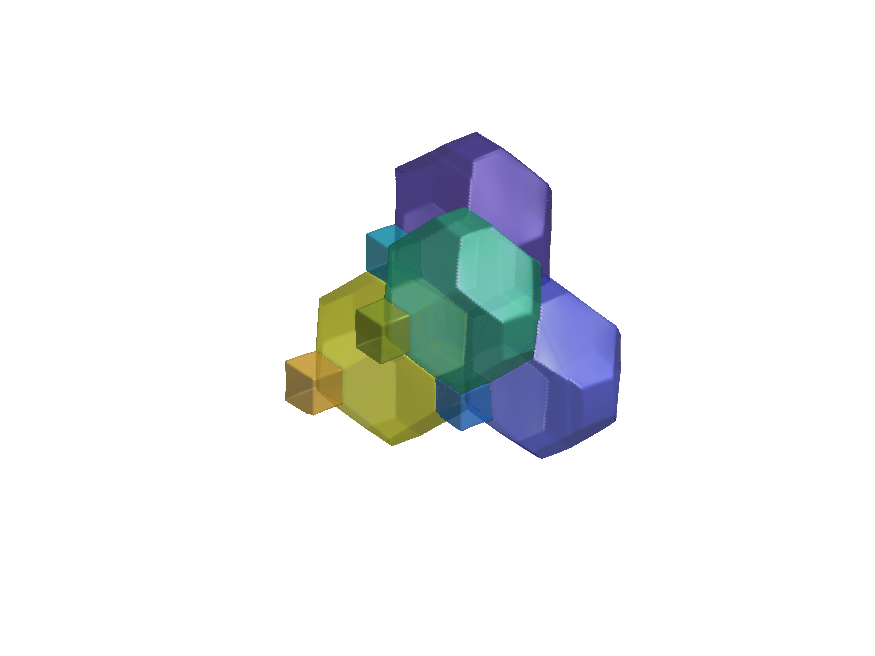}
            \\
     		\hline
    	\end{tabular}
\caption{A $n=8$ Dirichlet partition problem with periodic boundary conditions in the $4$-dimensional flat torus. The four columns correspond to the slides perpendicular to the $x_1$-, $x_2$-, $x_3$-, and $x_4$-axis respectively. The four rows correspond to the slices at $x_j = 0, 0.25, 0.5, 0.75$, respectively. These partitions are based on testing the trained neural network on a uniform $100^3$ grid. See Section~\ref{sec:periodic}.}
\label{Fig: Dirichlet Partition 4d unit tesseract n=8}
\end{figure}

{\bf Four-dimensional flat torus:} We then apply the proposed WANCO to $4$-dimensional partitioning problems. In Figures~\ref{Fig: Dirichlet Partition 4d unit tesseract n=2},~\ref{Fig: Dirichlet Partition 4d unit tesseract n=4} and~\ref{Fig: Dirichlet Partition 4d unit tesseract n=8}, the four columns correspond to slices perpendicular to the $x_1$-, $x_2$-, $x_3$-, and $x_4$-axes, respectively. The four rows correspond to the slices at $x_j = 0, 0.25, 0.5, 0.75$, respectively. These partitions are based on testing the trained neural network on a uniform $100^3$ grid at corresponding $3$-dimensional slices.

For $n=2$, we obtain a constant extension of the structure along the fourth direction, similar to the Schwarz P-surface, as displayed in Figure~\ref{Fig: Dirichlet Partition 4d unit tesseract n=2}. Traditional methods usually result in slab partitions for $2$ partitions in $3$-dimensional flat torus with a random initialization of the partition. However, based on the evidence from the numerical study in~\cite{Wang_2019},  the Schwarz P-surface also appears to be a local minimizer for the 3-dimensional case. In this $4$-dimensional flat torus, the random initialization represented by the neural network can find different local minimizers or stationary solutions more easily.

For $n=4$, we obtain a partition structure similar to a constant extension along the fourth direction of the rhombic dodecahedron as displayed in Figure~\ref{Fig: Dirichlet Partition 4d unit tesseract n=4}, which is the $4$-partition in the $3$-dimensional flat torus as displayed in Figure~\ref{Fig: D-P P Model 3d n4}.

For $n=8$, we obtain a partition that is symmetric along all four directions, as shown in Figure~\ref{Fig: Dirichlet Partition 4d unit tesseract n=8}. The obtained partition is similar to the structure known as a 24-cell honeycomb, which is a tessellation by regular 24-cells.

\section{Fluid-solid optimization}\label{Subsec: Fluid-solid optimization}

Topology optimization is a technique that aims to find the optimal layout or distribution of materials to maximize the objective within a specified design space while satisfying a set of constraints. Typically, these constraints are usually described by coupled partial differential equations. Such problems can be formulated as interface related optimization problems with PDE constraints. In the follows, we consider the application of WANCO to solve the fluid-solid optimization problem as an example.

To be specific, following \cite{Garcke_2015, Chen_2022}, we consider an optimization model to find the optimal fluid-solid channel inside a designated region with given inflow and outflow boundary conditions, where the interface is represented by the phase field order parameter $\phi$, and the objective is the total dissipation of the system in the following model. 

\begin{equation}
\begin{aligned}\label{Eq: topology optimization}
&\min_{\phi} \: C_\alpha J_\alpha(\phi,\mathbf{u})+C_\epsilon J_\epsilon(\phi) \\
\st \: &\begin{cases}-\Delta \mathbf{u}+\nabla p+\alpha(\phi) \mathbf{u}=\mathbf{0}, & \mathbf{x} \in D, \\ \nabla \cdot \mathbf{u}=0, & \mathbf{x} \in D,\\
\int_D\phi \ d\mathbf{x} = C_V|D|, & \mathbf{x} \in D,\\
\mathbf{u} = \mathbf{g}, & \mathbf{x} \in \partial D.
\end{cases}\end{aligned}\end{equation}Here,
\begin{equation}\label{Eq: J alpha}
    J_\alpha(\phi, \mathbf{u})=\int_D \frac{1}{2}|\nabla \mathbf{u}|^2 d \mathbf{x}+\int_D \frac{1}{2} \alpha(\phi)|\mathbf{u}|^2 d \mathbf{x},
\end{equation}
\begin{equation}
J_\epsilon(\phi)=\int_D \frac{\epsilon}{2}|\nabla \phi|^2+\frac{1}{\epsilon} F(\phi) d \mathbf{x},
\end{equation}
$\alpha(\phi)$ = $\alpha_0(1-\phi)^2$ is a penalized function. $\alpha_0$ is a constant and $F(\phi) = \frac{1}{4}\phi^2(1-\phi)^2$ is a double well potential. $C_\alpha$ and $C_\epsilon$ represent the weights of the two functional. $|D|$ represents the volume of the computational domain $D$ and $C_V$ denotes the volume fraction. $p$ is the fluid pressure and $\mathbf{g}$ is a given function defined on the boundary $\partial D$.

To simplify the above problem, we interpret the pressure $p$ as the Lagrange multiplier for the divergence-free constraint and treat $\mathbf{u}$ as the decision variable to minimize. Then, we arrive at the following problem:
\begin{equation}
    \begin{aligned}\label{Eq: topology optimization AL}
        &\min_{\mathbf{u},\phi} \: C_\alpha J_\alpha(\phi,\mathbf{u})+C_\epsilon J_\epsilon(\phi) \\
      \st \: &
      \begin{cases}
      \nabla \cdot \mathbf{u}=0, & \mathbf{x} \in D,\\
      \int_D\phi \ d\mathbf{x} = C_V|D|, & \mathbf{x} \in D,\\
      \mathbf{u} = \mathbf{g}, & \mathbf{x} \in \partial D.
      \end{cases}
    \end{aligned}
\end{equation}


The proposed WANCO is then applied to transform the model~\eqref{Eq: topology optimization AL} into the augmented Lagrangian formula. It is straightforward to check that the minimizer of~\eqref{Eq: topology optimization AL} also satisfies~\eqref{Eq: topology optimization}. For this specific problem, we introduce three Lagrange multipliers represented by three individual neural networks: $p(\mathbf{x};\eta_0), \lambda_1(0;\eta_1), \lambda_2(\mathbf{x};\eta_2)$, to incorporate the divergence-free constraint, boundary condition constraint and volume preserving constraint, respectively. We then arrive at the following minimax problem:
\begin{equation}\label{Eq: topology optimization parameter form}
    \begin{aligned}
    \min_{\theta} \max_{\eta_0,\eta_1,\eta_2}\:L_\beta(\mathbf{u}(\mathbf{x};\theta),\phi(\mathbf{x};\theta),&p(\mathbf{x};\eta_0),\lambda_1(0;\eta_1),\lambda_2(\mathbf{x};\eta_2))
    \end{aligned}
\end{equation}
where 
\begin{align*}
& L_\beta(\mathbf{u}(\mathbf{x};\theta),\phi(\mathbf{x};\theta),p(\mathbf{x};\eta_0),\lambda_1(0;\eta_1),\lambda_2(\mathbf{x};\eta_2)) =  C_\alpha J_\alpha(\phi(\mathbf{x};\theta),\mathbf{u}(\mathbf{x};\theta))\\
&+ C_\epsilon J_\epsilon(\phi(\mathbf{x};\theta)) - \int_Dp(\mathbf{x};\eta_0)\nabla \cdot \mathbf{u}(\mathbf{x};\theta)d\mathbf{x} + \frac{\beta_u}{2} \int_D(\nabla \cdot \mathbf{u}(\mathbf{x};\theta))^2d\mathbf{x} \\
    & - \lambda_1(0;\eta_1)(\int_D\phi(\mathbf{x};\theta) d\mathbf{x}-C_V|D|) + \frac{\beta_\phi}{2}(\int_D\phi(\mathbf{x};\theta) d\mathbf{x}-C_V|D|)^2\\
    & - \int_{\partial D}\lambda_2(\mathbf{x};\eta_2)\left(\mathbf{u}(\mathbf{x};\theta) - \mathbf{g}(\mathbf{x})\right)\ ds + \frac{\beta_b}{2}\int_{\partial D} (\mathbf{u}(\mathbf{x};\theta) - \mathbf{g}(\mathbf{x}))^2 \ ds.
\end{align*}
Here, the primal network has an output dimension of $3$, which includes the two entries in $\mathbf{u}(\mathbf{x};\theta)$ and the phase variable $\phi(\mathbf{x};\theta)$. The adversarial network consists of three individual networks, with one involving the Lagrange multipliers corresponding to $\nabla \cdot \mathbf{u}$ as pressure $p(\mathbf{x};\eta_0)$, and the other two involving the multipliers $\lambda_1(0;\eta_1)$ and $\lambda_2(\mathbf{x};\eta_2)$ corresponding to the volume preserving and boundary condition constraints.

In the follows, we consider two examples of the fluid-solid optimization problem using WANCO. The parameters are set as follows: $\epsilon=0.01$, $\alpha_0 = 250000$, $C_\alpha=100$, $C_\epsilon=10$, $\beta_u=100$, $\beta_\phi=100$, $\beta_b = 1000$, and $\alpha_u=\alpha_\phi=\alpha_b=1.0003$. In Example 1, the primal network is set to $N_d=4$ and $N_w=50$, while in Example 2, $N_d=6$ and $N_w=80$. The constructions of three adversarial neural networks are as follows. For the Lagrange multiplier $p$, we use a ResNet with the same depth and width as the primal network to represent it as $p(\mathbf{x};\eta_0)$ for each example. For the Lagrange multiplier $\lambda_1$, which enforces the area constraint, we use a shallow neural network with $1$ hidden layer and width $10$, denoted as $\lambda_1(0;\eta_1)$ to represent it. For the Lagrange multiplier $\lambda_2(\mathbf{x};\eta_2)$, which enforces the boundary conditions, we use a ResNet with $4$ hidden layers and width $50$ to represent it. The training is performed for $N=5000$ steps with inner iteration steps $N_u=3$, $N_p=1=N_\phi=N_b=1$, and $N_r=40000$ points are used for training in each iteration. The testing results are obtained by evaluating the trained neural network on a uniform $1000\times 1000$ grid. 

\begin{figure}[ht!]
\centering
{\includegraphics[width = 0.31 \textwidth]{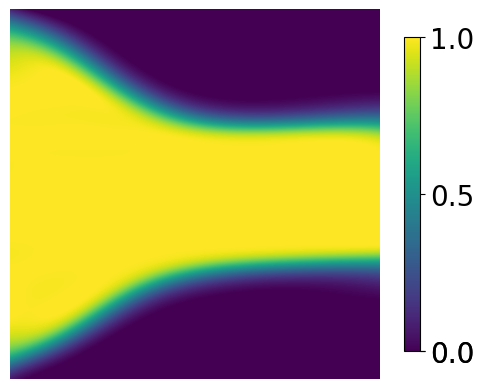}}
{\includegraphics[width = 0.25 \textwidth]{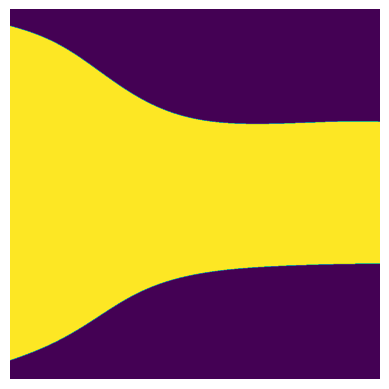}}
{\includegraphics[width = 0.27 \textwidth]{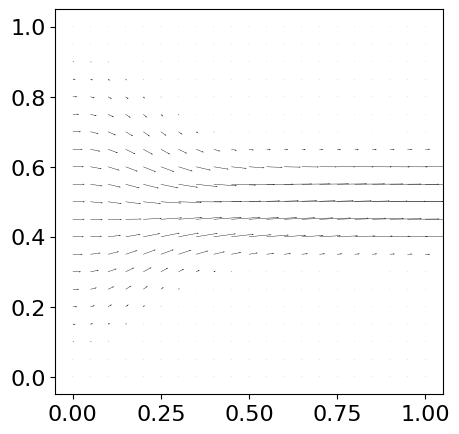}}
\caption{The phase variable (left), projection fluid-solid region (middle), and the velocity field (right) for Example $1$. See Section~\ref{Subsec: Fluid-solid optimization}}
\label{Fig: Fluid-Solid Optimization example 1}
\end{figure}

{\bf Example $1$:} In Figure~\ref{Fig: Fluid-Solid Optimization example 1}, we display the results of fluid-solid optimization in $D=[0,1]\times[0,1]$ with $C_v = \frac{1}{2}$ and the following boundary conditions:

\begin{equation}\label{Eq: case1 boundary condition top opt}
    \begin{aligned}
 &
      \begin{cases}u_1(0,y)=\frac{1}{2}\sin(y\pi), & y \in [0,1], \\ u_1(1,y)=\frac{3}{2}\sin((3y-1)\pi), & y \in [\frac{1}{3},\frac{2}{3}],\\
      u_1(x,y) = u_2(x,y) = 0, & \text{elsewhere on }\partial D.

      \end{cases}
    \end{aligned}
\end{equation}
The training results, including the phase variable $\phi(\mathbf{x};\theta)$ and the velocity field $\mathbf{u}(\mathbf{x};\theta)$, are presented in Figure~\ref{Fig: Fluid-Solid Optimization example 1}. The training result exhibits a smooth transition from a large opening on the left to a smaller opening on the right. The second column shows the phase variable based on the projection given below:
\begin{equation}\label{formula: projection F-S opt}
    \begin{aligned}
 &
      \begin{cases} \Tilde{\phi}(\mathbf{x};\theta)=1 & \text{if} \quad \phi(\mathbf{x};\theta)\geq\frac{1}{2}, \\ \Tilde{\phi}(\mathbf{x};\theta)=0 & \text{if} \quad \phi(\mathbf{x};\theta)<\frac{1}{2}.
      \end{cases}
    \end{aligned}
\end{equation}

\begin{figure}[h!]
  \centering

{\includegraphics[width = 0.31 \textwidth]{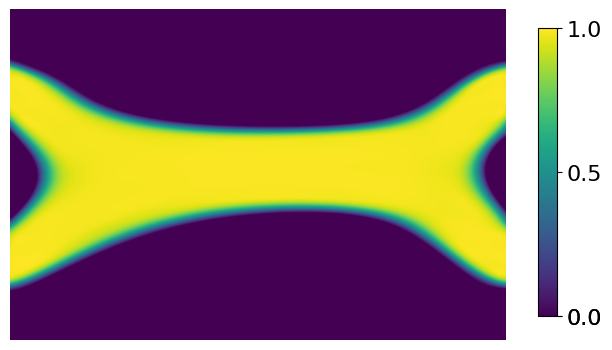}}
{\includegraphics[width = 0.253 \textwidth]{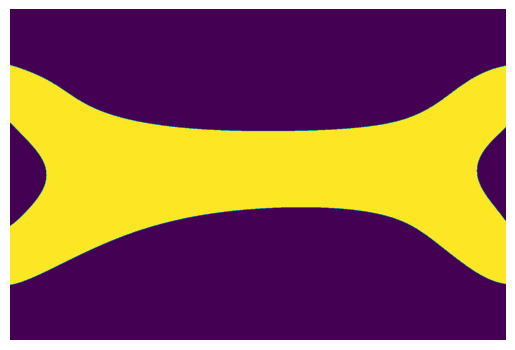}}
{\includegraphics[width = 0.262 \textwidth]{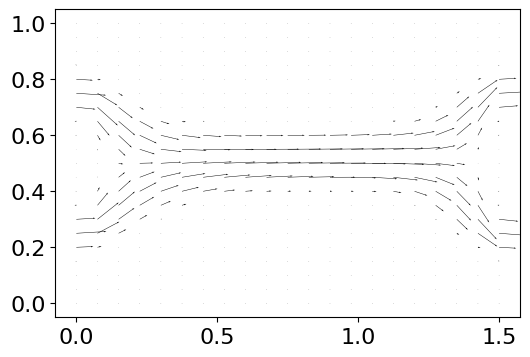}}
{\includegraphics[width = 0.33 \textwidth]{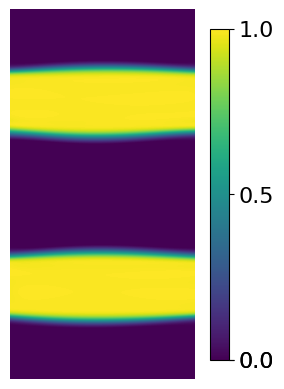}}
{\includegraphics[width = 0.238\textwidth]{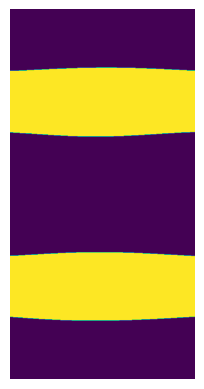}}
{\includegraphics[width = 0.26 \textwidth]{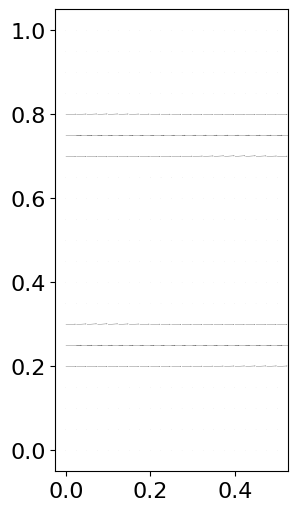}}
\caption{The phase variable (left), projection fluid-solid region (middle), and the velocity field (right) for Example $2$. See Section~\ref{Subsec: Fluid-solid optimization}.}
\label{Fig: Fluid-Solid Optimization example 2}
\end{figure}

{\bf Example $2$:} In Figure~\ref{Fig: Fluid-Solid Optimization example 2}, we display the results of fluid-solid optimization problems with the same boundary conditions but with different optimization regions, denoted as $D = [0,l]\times[0,1]$, where $l$ denotes the length of the desired domain taking two values: $\frac{3}{2}$ and $\frac{1}{2}$. The boundary conditions are given below and $C_v = \frac{1}{3}$. 

\begin{equation}\label{Eq: case3 boundary condition top opt}
    \begin{aligned}
 &
      \begin{cases}u_1(0,y)=u_1(l,y) = 1-(12y-3)^2, & y \in [\frac{1}{6},\frac{2}{6}], \\ u_1(0,y)=u_1(l,y) = 1-(12y-9)^2, & y \in [\frac{4}{6},\frac{5}{6}],\\
      u_1(x,y) = u_2(x,y) = 0, & \text{elsewhere on }\partial D.

      \end{cases}
    \end{aligned}
\end{equation}
The training results are presented in Figure~\ref{Fig: Fluid-Solid Optimization example 2}. The first and second rows correspond to the results of $l=\frac{3}{2}$ and $l=\frac{1}{2}$, respectively. It is noteworthy that in this example, the boundary conditions, network architecture, parameter settings, and other configurations are the same, with the only difference being the optimization region, \ie, different values of $l$. By employing the identical WANCO procedure, we can obtain the corresponding optimal channels. Specifically, when $l=\frac{1}{2}$, two independent and disconnected channels are formed. When $l=\frac{3}{2}$, the channels first merge and then separate, forming an X-shaped channel. These training results agree well with the results in~\cite{Chen_2022}, further highlighting the capabilities of WANCO.

\section{Obstacle Problems}\label{Subsec: Obstacle problem}
In this section, we consider the following obstacle problem that aims to minimize an objective functional subject to inequality constraints that represent obstacle conditions, as well as some boundary conditions,
\begin{equation}\label{Eq: Constrained Optimization Original Form obstacle problems}
    \begin{aligned}
    \min\limits_{u}\:  & \int_{\Omega}|\nabla u|^2 d \mathbf{x} \\
      \st \: &
      \begin{cases}
      u \geq \psi, & \mathbf{x} \in \Omega,\\
      u = g, & \mathbf{x} \in \partial \Omega.\\
      \end{cases}
      \end{aligned}
\end{equation}
This problem is naturally the constrained optimization framework with inequality constraints. By following the augmented Lagrangian method and incorporating the concept of adversarial networks, one can obtain the minimax problem,
\begin{equation}\label{Eq: WANCO 1d obstacle problem parameter form}
\begin{aligned}
           \min_{\theta}& \max_{\eta}\: L_\beta(u(\mathbf{x};\theta),\lambda(\mathbf{x};\eta))= C_0\int_{\Omega}|\nabla u(\mathbf{x};\theta)|^2 d\mathbf{x}\\
           &-\int_\Omega\lambda(\mathbf{x};\eta)(\psi(\mathbf{x})-u(\mathbf{x};\theta))d\mathbf{x} + \frac{\beta}{2}\int(\mathrm{ReLU}(\psi(\mathbf{x})-u(\mathbf{x};\theta)))^2d\mathbf{x},
\end{aligned}
\end{equation}
where the Lagrange multiplier $\lambda$ is non-positive.

We consider some $1$-dimensional obstacle problems in $\Omega = [0,1]$ with different obstacles $\psi$ and boundary conditions as follows, 
\begin{equation}
\psi_1(x)= \begin{cases}100 x^2, & \text { for } 0 \leq x \leq 0.25, \\ 100 x(1-x)-12.5, & \text { for } 0.25 \leq x \leq 0.5,\qquad u(0) = u(1)=0,\\ \psi_1(1-x), & \text { for } 0.5 \leq x \leq 1.\end{cases}
\end{equation}

\begin{equation}
\psi_2(x)= \begin{cases}10 \sin (2 \pi x), & \text { for } 0 \leq x \leq 0.25, \\ 5 \cos (\pi(4 x-1))+5, & \text { for } 0.25 \leq x \leq 0.5,\qquad u(0) = u(1)=0, \\ \psi_2(1-x), & \text { for } 0.5 \leq x \leq 1.\end{cases}
\end{equation}

\begin{equation}
\psi_3(x)=10 \sin ^2\left(\pi(x+1)^2\right), \quad\quad \text { for } 0 \leq x \leq 1\text{,}\qquad\quad  u(0) = 5\text{, }u(1)=10.
\end{equation}
In these examples, we set $C_0=100$, $\beta=1000$ and penalty amplifier $\alpha=1.0003 $. The primal network is a ResNet with $6$ hidden layers and width $80$. We impose the non-positivity of the Lagrange multiplier by modifying the output of the adversarial network $\lambda(x;\eta)=-\mathrm{ReLU}(-\lambda(x;\eta))$, where the adversarial network has the same structure as the primal network. Regarding the boundary conditions, we similarly enforce them strictly using the construction method as used in~\eqref{formula: output modification of u-net G-L energy}. The training is performed for $N=5000$ with inner iteration steps $N_u=N_\lambda=2$, and $N_r=2000$ points are used for training in each iteration. The testing results are obtained by evaluating the trained neural network on a uniform grid of $1000$ points.

\begin{center}
\begin{figure}[ht!]
\centering
{\includegraphics[width = 0.32 \textwidth]{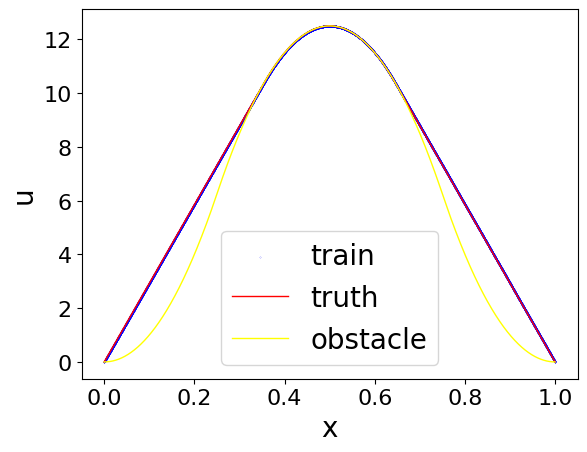}}
{\includegraphics[width = 0.32 \textwidth]{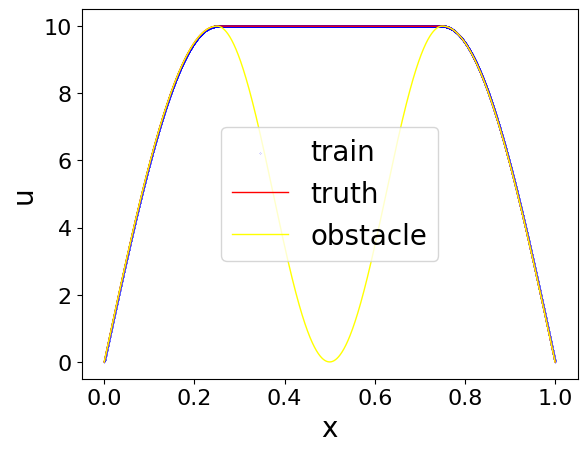}}
{\includegraphics[width = 0.32 \textwidth]{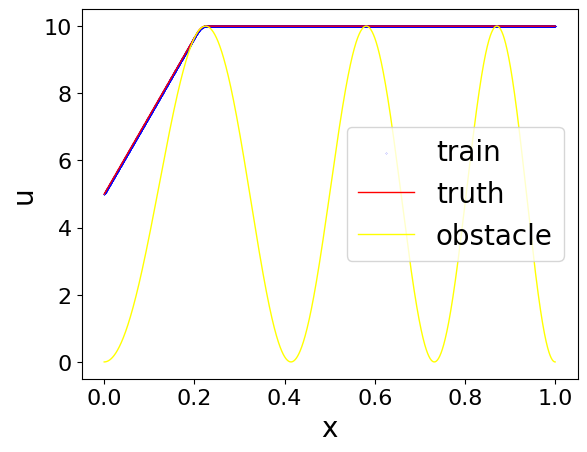}}\\
{\includegraphics[width = 0.32 \textwidth]{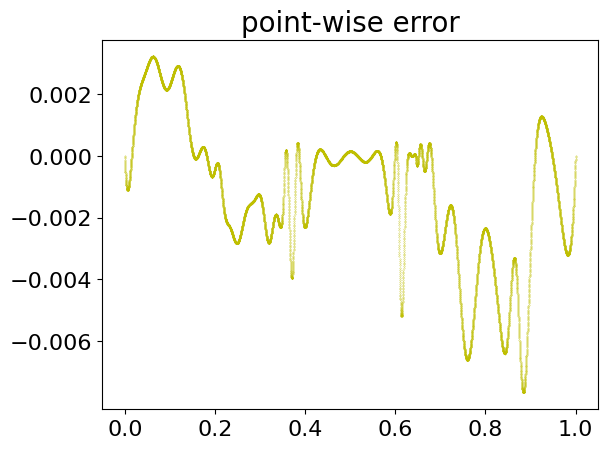}}
{\includegraphics[width = 0.32 \textwidth]{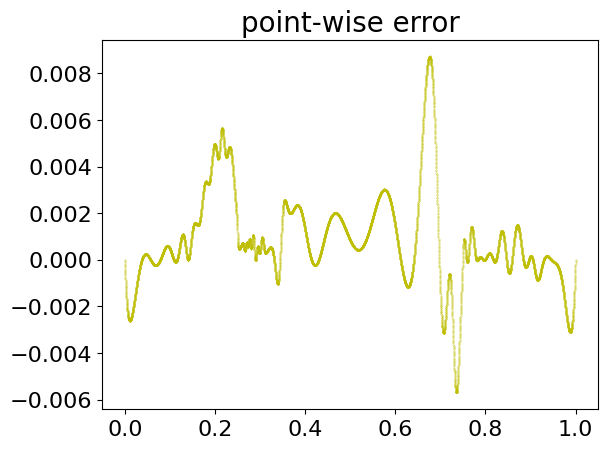}}
{\includegraphics[width = 0.32 \textwidth]{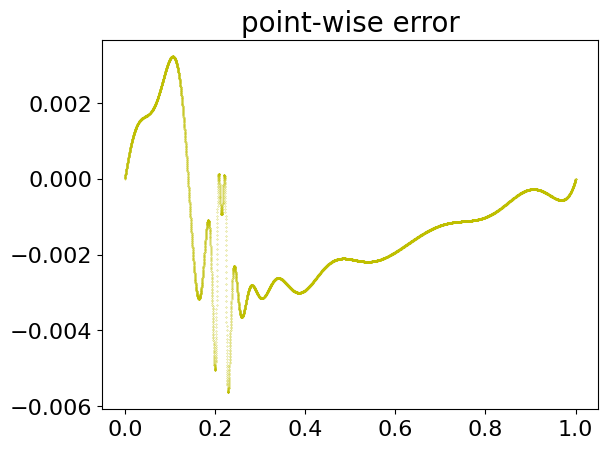}}
\caption{Training results for obstacle problems under three different obstacles and boundary conditions. The first row shows the obstacles, true solutions for the corresponding obstacles, and the training results for each case. The second row represents the point-wise error between the true solution and the trained solution. The left, middle, and right columns correspond to obstacles $\psi_1$, $\psi_2$, and $\psi_3$, respectively. See Section~\ref{Subsec: Obstacle problem}.}
\label{Fig: obstacle problem}
\end{figure} 
\end{center}
As shown in Figure~\ref{Fig: obstacle problem}, WANCO can effectively simulate the obstacle problems and accurately preserve the constraints. This can be observed from the contact regions of the obstacles in the trained solution and the boundary constraints are strictly preserved. In comparison to recent papers~\cite{bahja2023physics, Cheng_2023} that utilize neural networks for solving obstacle problems, WANCO demonstrates certain improvements in accuracy. For instance, when considering the obstacle $\psi_1$, the infinity norm of the point-wise error is approximately $0.012$ as reported in~\cite{bahja2023physics}, whereas it is around $0.09$ as reported in~\cite{Cheng_2023}.

\section{Conclusion and discussion}\label{Sec: Conclusions}

In this paper, we introduced WANCO (Weak Adversarial Networks for Constrained Optimization), a novel approach inspired by WAN~\cite{Zang_2020, Bao_2020}. WANCO utilizes the augmented Lagrangian formulation and separate networks for decision variables and Lagrange multipliers. It exhibits robustness to parameter tuning and adaptive training through the adversarial term. WANCO is a versatile framework that can enhance training efficiency and achieve improved results when integrated with other methods. Our experiments demonstrated the effectiveness of combining the $\tanh^3$ activation function with ResNet for WANCO. We presented WANCO's performance across a range of constrained optimization problems, encompassing scalar linear/nonlinear constraints, PDE constraints, and inequality constraints. The results highlight WANCO as a stable and adaptable solution to constrained optimization problems, demonstrating promising potential for wider applications.

\bibliographystyle{siamplain}
\bibliography{references}

\begin{thebibliography}{10}

\bibitem{bahja2023physics}
{\sc H.~E. Bahja, J.~C. Hauffen, P.~Jung, B.~Bah, and I.~Karambal}, {\em A
  physics-informed neural network framework for modeling obstacle-related
  equations}, arXiv:2304.03552,  (2023).

\bibitem{bao2022mathematical}
{\sc G.~Bao}, {\em Mathematical analysis and numerical methods for inverse
  scattering problems}, in Proc. Int. Cong. Math, vol.~7, 2022, pp.~5034--5055,
  \url{https://doi.org/10.4171/ICM2022/5}.

\bibitem{bao2001mathematical}
{\sc G.~Bao, L.~Cowsar, and W.~Masters}, {\em Mathematical modeling in optical
  science}, SIAM, 2001.

\bibitem{Bao_2015}
{\sc G.~Bao, P.~Li, J.~Lin, and F.~Triki}, {\em Inverse scattering problems
  with multi-frequencies}, Inverse Problems, 31 (2015), p.~093001,
  \url{https://doi.org/10.1088/0266-5611/31/9/093001}.

\bibitem{Bao_2020}
{\sc G.~Bao, X.~Ye, Y.~Zang, and H.~Zhou}, {\em Numerical solution of inverse
  problems by weak adversarial networks}, Inverse Problems, 36 (2020),
  p.~115003, \url{https://doi.org/10.1088/1361-6420/abb447}.

\bibitem{bendsoe2013topology}
{\sc M.~P. Bendsoe and O.~Sigmund}, {\em Topology optimization: theory,
  methods, and applications}, Springer Science \& Business Media, 2013.

\bibitem{bertsekas2014constrained}
{\sc D.~P. Bertsekas}, {\em Constrained optimization and Lagrange multiplier
  methods}, Academic press, 2014.

\bibitem{Bogosel_2016}
{\sc B.~Bogosel and B.~Velichkov}, {\em A multiphase shape optimization problem
  for eigenvalues: Qualitative study and numerical results}, SIAM Journal on
  Numerical Analysis, 54 (2016), p.~210–241,
  \url{https://doi.org/10.1137/140976406}.

\bibitem{Burman_2023}
{\sc E.~Burman, P.~Hansbo, and M.~G. Larson}, {\em The augmented {L}agrangian
  method as a framework for stabilised methods in computational mechanics},
  Archives of Computational Methods in Engineering, 30 (2023), p.~2579–2604,
  \url{https://doi.org/10.1007/s11831-022-09878-6}.

\bibitem{caffarelli2007optimal}
{\sc L.~A. Cafferelli and F.~H. Lin}, {\em An optimal partition problem for
  eigenvalues}, J. Sci. Comp., 31 (2007), pp.~5--18,
  \url{https://doi.org/10.1007/s10915-006-9114-8}.

\bibitem{Chen_2022}
{\sc H.~Chen, H.~Leng, D.~Wang, and X.-P. Wang}, {\em An efficient threshold
  dynamics method for topology optimization for fluids}, CSIAM Transactions on
  Applied Mathematics, 3 (2022), p.~26–56,
  \url{https://doi.org/10.4208/csiam-am.so-2021-0007}.

\bibitem{Chen_2022_RF}
{\sc J.~Chen, X.~Chi, W.~E, and Z.~Yang}, {\em Bridging traditional and machine
  learning-based algorithms for solving {PDE}s: The random feature method},
  Journal of Machine Learning, 1 (2022), p.~268–298,
  \url{https://doi.org/10.4208/jml.220726}.

\bibitem{Cheng_2020}
{\sc Q.~Cheng and J.~Shen}, {\em Global constraints preserving scalar auxiliary
  variable schemes for gradient flows}, SIAM Journal on Scientific Computing,
  42 (2020), p.~A2489–A2513, \url{https://doi.org/10.1137/19m1306221}.

\bibitem{Cheng_2023}
{\sc X.~Cheng, X.~Shen, X.~Wang, and K.~Liang}, {\em A deep neural
  network-based method for solving obstacle problems}, Nonlinear Analysis: Real
  World Applications, 72 (2023), p.~103864,
  \url{https://doi.org/10.1016/j.nonrwa.2023.103864}.

\bibitem{Chu_2021}
{\sc K.~Chu and S.~Leung}, {\em A level set method for the {D}irichlet
  k-partition problem}, Journal of Scientific Computing, 86 (2021),
  \url{https://doi.org/10.1007/s10915-020-01368-w}.

\bibitem{Cuomo_2022}
{\sc S.~Cuomo, V.~S. Di~Cola, F.~Giampaolo, G.~Rozza, M.~Raissi, and
  F.~Piccialli}, {\em Scientific machine learning through physics–informed
  neural networks: Where we are and what’s next}, Journal of Scientific
  Computing, 92 (2022), \url{https://doi.org/10.1007/s10915-022-01939-z}.

\bibitem{Cybenko_1989}
{\sc G.~Cybenko}, {\em Approximation by superpositions of a sigmoidal
  function}, Mathematics of Control, Signals, and Systems, 2 (1989),
  p.~303–314, \url{https://doi.org/10.1007/bf02551274}.

\bibitem{Dong_2021}
{\sc S.~Dong and Z.~Li}, {\em Local extreme learning machines and domain
  decomposition for solving linear and nonlinear partial differential
  equations}, Computer Methods in Applied Mechanics and Engineering, 387
  (2021), p.~114129, \url{https://doi.org/10.1016/j.cma.2021.114129}.

\bibitem{Du_2008}
{\sc Q.~Du and F.~Lin}, {\em Numerical approximations of a norm-preserving
  gradient flow and applications to an optimal partition problem},
  Nonlinearity, 22 (2008), pp.~67--83,
  \url{https://doi.org/10.1088/0951-7715/22/1/005}.

\bibitem{E_2018}
{\sc W.~E and B.~Yu}, {\em The deep {R}itz method: A deep learning-based
  numerical algorithm for solving variational problems}, Communications in
  Mathematics and Statistics, 6 (2018), p.~1–12,
  \url{https://doi.org/10.1007/s40304-018-0127-z}.

\bibitem{Garcke_2015}
{\sc H.~Garcke, C.~Hecht, M.~Hinze, and C.~Kahle}, {\em Numerical approximation
  of phase field based shape and topology optimization for fluids}, SIAM
  Journal on Scientific Computing, 37 (2015), p.~A1846–A1871,
  \url{https://doi.org/10.1137/140969269}.

\bibitem{glowinski1989augmented}
{\sc R.~Glowinski and P.~Le~Tallec}, {\em Augmented Lagrangian and
  operator-splitting methods in nonlinear mechanics}, SIAM, 1989.

\bibitem{He_2016}
{\sc K.~He, X.~Zhang, S.~Ren, and J.~Sun}, {\em Deep residual learning for
  image recognition}, in 2016 IEEE Conference on Computer Vision and Pattern
  Recognition (CVPR), IEEE, June 2016,
  \url{https://doi.org/10.1109/cvpr.2016.90}.

\bibitem{Huang_2022}
{\sc J.~Huang, H.~Wang, and T.~Zhou}, {\em An augmented {L}agrangian deep
  learning method for variational problems with essential boundary conditions},
  Communications in Computational Physics, 31 (2022), p.~966–986,
  \url{https://doi.org/10.4208/cicp.oa-2021-0176}.

\bibitem{Hwang_2022}
{\sc H.~J. Hwang and H.~Son}, {\em Lagrangian dual framework for conservative
  neural network solutions of kinetic equations}, Kinetic and Related Models,
  15 (2022), p.~551, \url{https://doi.org/10.3934/krm.2021046}.

\bibitem{Karniadakis_2021}
{\sc G.~E. Karniadakis, I.~G. Kevrekidis, L.~Lu, P.~Perdikaris, S.~Wang, and
  L.~Yang}, {\em Physics-informed machine learning}, Nature Reviews Physics, 3
  (2021), p.~422–440, \url{https://doi.org/10.1038/s42254-021-00314-5}.

\bibitem{kharazmi2019variational}
{\sc E.~Kharazmi, Z.~Zhang, and G.~E. Karniadakis}, {\em Variational
  physics-informed neural networks for solving partial differential equations},
  arXiv:1912.00873,  (2019).

\bibitem{kingma2014adam}
{\sc D.~P. Kingma and J.~Ba}, {\em Adam: A method for stochastic optimization},
  arXiv:1412.6980,  (2014).

\bibitem{Li_2021}
{\sc J.~Li, D.~Fridovich-Keil, S.~Sojoudi, and C.~J. Tomlin}, {\em Augmented
  {L}agrangian method for instantaneously constrained reinforcement learning
  problems}, in 2021 60th IEEE Conference on Decision and Control (CDC), IEEE,
  Dec. 2021, \url{https://doi.org/10.1109/cdc45484.2021.9683088}.

\bibitem{Linka_2022}
{\sc K.~Linka, A.~Schäfer, X.~Meng, Z.~Zou, G.~E. Karniadakis, and E.~Kuhl},
  {\em Bayesian physics informed neural networks for real-world nonlinear
  dynamical systems}, Computer Methods in Applied Mechanics and Engineering,
  402 (2022), p.~115346, \url{https://doi.org/10.1016/j.cma.2022.115346}.

\bibitem{Lu_2021}
{\sc L.~Lu, X.~Meng, Z.~Mao, and G.~E. Karniadakis}, {\em Deep{XDE}: A deep
  learning library for solving differential equations}, SIAM Review, 63 (2021),
  p.~208–228, \url{https://doi.org/10.1137/19m1274067}.

\bibitem{Lu_2021hard}
{\sc L.~Lu, R.~Pestourie, W.~Yao, Z.~Wang, F.~Verdugo, and S.~G. Johnson}, {\em
  Physics-informed neural networks with hard constraints for inverse design},
  SIAM Journal on Scientific Computing, 43 (2021), p.~B1105–B1132,
  \url{https://doi.org/10.1137/21m1397908}.

\bibitem{Lyu_2021}
{\sc L.~Lyu, K.~Wu, R.~Du, and J.~Chen}, {\em Enforcing exact boundary and
  initial conditions in the deep mixed residual method}, CSIAM Transactions on
  Applied Mathematics, 2 (2021), p.~748–775,
  \url{https://doi.org/10.4208/csiam-am.so-2021-0011}.

\bibitem{McClenny_2022}
{\sc L.~McClenny and U.~Braga-Neto}, {\em Self-adaptive physics-informed neural
  networks}, SSRN Electronic Journal,  (2022),
  \url{https://doi.org/10.2139/ssrn.4086448}.

\bibitem{nandwani2019primal}
{\sc Y.~Nandwani, A.~Pathak, and P.~Singla}, {\em A primal dual formulation for
  deep learning with constraints}, Advances in Neural Information Processing
  Systems, 32 (2019).

\bibitem{Qi_2022}
{\sc Q.~Qi, W.~Lin, B.~Guo, J.~Chen, C.~Deng, G.~Lin, X.~Sun, and Y.~Chen},
  {\em Augmented {L}agrangian-based reinforcement learning for network slicing
  in {IIoT}}, Electronics, 11 (2022), p.~3385,
  \url{https://doi.org/10.3390/electronics11203385}.

\bibitem{Raissi_2019}
{\sc M.~Raissi, P.~Perdikaris, and G.~Karniadakis}, {\em Physics-informed
  neural networks: A deep learning framework for solving forward and inverse
  problems involving nonlinear partial differential equations}, Journal of
  Computational Physics, 378 (2019), p.~686–707,
  \url{https://doi.org/10.1016/j.jcp.2018.10.045}.

\bibitem{sangalli2021constrained}
{\sc S.~Sangalli, E.~Erdil, A.~H{\"o}tker, O.~Donati, and E.~Konukoglu}, {\em
  Constrained optimization to train neural networks on critical and
  under-represented classes}, Advances in neural information processing
  systems, 34 (2021), pp.~25400--25411.

\bibitem{Sharma_2020}
{\sc S.~Sharma, S.~Sharma, and A.~Athaiya}, {\em Activation functions in neural
  networks}, International Journal of Engineering Applied Sciences and
  Technology, 04 (2020), p.~310–316,
  \url{https://doi.org/10.33564/ijeast.2020.v04i12.054}.

\bibitem{Sirignano_2018}
{\sc J.~Sirignano and K.~Spiliopoulos}, {\em {DGM}: A deep learning algorithm
  for solving partial differential equations}, Journal of Computational
  Physics, 375 (2018), p.~1339–1364,
  \url{https://doi.org/10.1016/j.jcp.2018.08.029}.

\bibitem{Son_2023}
{\sc H.~Son, S.~W. Cho, and H.~J. Hwang}, {\em Enhanced physics-informed neural
  networks with augmented {L}agrangian relaxation method ({AL-PINNs})},
  Neurocomputing, 548 (2023), p.~126424,
  \url{https://doi.org/10.1016/j.neucom.2023.126424}.

\bibitem{Sukumar_2022}
{\sc N.~Sukumar and A.~Srivastava}, {\em Exact imposition of boundary
  conditions with distance functions in physics-informed deep neural networks},
  Computer Methods in Applied Mechanics and Engineering, 389 (2022), p.~114333,
  \url{https://doi.org/10.1016/j.cma.2021.114333}.

\bibitem{Wang_2022_DP}
{\sc D.~Wang}, {\em An efficient unconditionally stable method for dirichlet
  partitions in arbitrary domains}, SIAM Journal on Scientific Computing, 44
  (2022), p.~A2061–A2088, \url{https://doi.org/10.1137/21m1443406}.

\bibitem{Wang_2019}
{\sc D.~Wang and B.~Osting}, {\em A diffusion generated method for computing
  {D}irichlet partitions}, Journal of Computational and Applied Mathematics,
  351 (2019), p.~302–316, \url{https://doi.org/10.1016/j.cam.2018.11.015}.

\bibitem{Wang_2024}
{\sc S.~Wang, S.~Sankaran, and P.~Perdikaris}, {\em Respecting causality for
  training physics-informed neural networks}, Computer Methods in Applied
  Mechanics and Engineering, 421 (2024), p.~116813,
  \url{https://doi.org/10.1016/j.cma.2024.116813}.

\bibitem{wright1997primal}
{\sc S.~J. Wright}, {\em Primal-dual interior-point methods}, SIAM, 1997.

\bibitem{Wu_2023}
{\sc C.~Wu, M.~Zhu, Q.~Tan, Y.~Kartha, and L.~Lu}, {\em A comprehensive study
  of non-adaptive and residual-based adaptive sampling for physics-informed
  neural networks}, Computer Methods in Applied Mechanics and Engineering, 403
  (2023), p.~115671, \url{https://doi.org/10.1016/j.cma.2022.115671}.

\bibitem{Yu_2021}
{\sc H.~Yu, X.~Tian, W.~E, and Q.~Li}, {\em Onsager{N}et: Learning stable and
  interpretable dynamics using a generalized {O}nsager principle}, Physical
  Review Fluids, 6 (2021),
  \url{https://doi.org/10.1103/physrevfluids.6.114402}.

\bibitem{Yu_2022}
{\sc J.~Yu, L.~Lu, X.~Meng, and G.~E. Karniadakis}, {\em Gradient-enhanced
  physics-informed neural networks for forward and inverse {PDE} problems},
  Computer Methods in Applied Mechanics and Engineering, 393 (2022), p.~114823,
  \url{https://doi.org/10.1016/j.cma.2022.114823}.

\bibitem{zang2023particlewnn}
{\sc Y.~Zang and G.~Bao}, {\em Particle{WNN}: a novel neural networks framework
  for solving partial differential equations}, arXiv:2305.12433,  (2023).

\bibitem{Zang_2020}
{\sc Y.~Zang, G.~Bao, X.~Ye, and H.~Zhou}, {\em Weak adversarial networks for
  high-dimensional partial differential equations}, Journal of Computational
  Physics, 411 (2020), p.~109409,
  \url{https://doi.org/10.1016/j.jcp.2020.109409}.

\end{thebibliography}
\end{document}